\documentclass{amsart}
\usepackage{graphicx} 
\usepackage{amsmath,amssymb,amsthm,latexsym,amsfonts,mathtools,tikz,mathrsfs} 
\usepackage[top=30truemm,bottom=30truemm,left=25truemm,right=25truemm]{geometry}
\usepackage{url}

\usepackage{hyperref}

\newtheorem{thm}{Theorem}[section]
\newtheorem{lem}[thm]{Lemma}

\newtheorem{rem}[thm]{Remark}
\newtheorem{cor}[thm]{Corollary}

\theoremstyle{definition}

\makeatletter
\@addtoreset{equation}{section}

\makeatother

\def\Wick#1{\mathord{{:}{#1}{:}}}

\allowdisplaybreaks

\title[Maximal dissipativity for the Dirichlet operator associated to the dynamical $\Phi^4_2$ model]{Maximal dissipativity for the Dirichlet operator associated to the dynamical $\Phi^4$ model on the plane}

\author{Hirotatsu Nagoji}
\address{Graduate School of Advanced Science and Engineering, Hiroshima University, 1-3-1 Kagamiyama, Higashi-Hiroshima City, Hiroshima, 739-8526, Japan}
\email{nagojihiro@hiroshima-u.ac.jp}

\date{}
\keywords{Euclidean quantum field theory, stochastic quantization, singular SPDE, Dirichlet operator, Kolmogorov operator}
\subjclass[2020]{81S20, 60H15, 35Q40, 35R60}

\begin{document}

\begin{abstract}
We consider the Dirichlet/Kolmogorov operator associated with the stochastic quantization of the $\Phi^4$ model on the plane whose domain is restricted to the space of smooth cylindrical functions. We prove that its closure is maximally dissipative in $L^p$. Our result implies the so-called $L^p$-uniqueness of the Dirichlet operator for every $1\le p<\infty$. In particular, essential self-adjointness and Markov uniqueness follow.
The key ingredient is an estimate for the linearized equation for the dynamical $\Phi^4_2$ model.
The proof is based on techniques from singular SPDE theory and the stopping time argument introduced by Kunick and Tsatsoulis \cite{KT25}.

\end{abstract}

\maketitle
\tableofcontents

\allowdisplaybreaks

\section{Introduction}\label{sec:intro}

The $\Phi^4_2$-measure is a probability measure on the space of tempered distributions $\mathcal{S}'(\mathbb{R}^2)$, formally given by
\begin{equation}
    \mu^{(\lambda)}(d\phi)
    = \frac{1}{Z}\exp\left(-\frac{\lambda}{4}
      \int_{\mathbb{R}^2}\Wick{\phi^4}(x)\,dx\right)\mu_0(d\phi),
\end{equation}
which arises in constructive quantum field theory; see, for instance, \cite{S74}.
Here, $\mu_0$ is the massive Gaussian free field measure on $\mathbb{R}^2$ with covariance $(1-\Delta)^{-1}$, $\lambda>0$ is the coupling constant, and $\Wick{\phi^k}$ denotes the $k$th Wick power.
The expression above is only formal in infinite volume; the measure used in this paper is specified below.
We are concerned with the associated parabolic stochastic quantization equation
\begin{equation}\label{eq:sq}
    \partial_t\Phi_t
    = \frac{1}{2}(\Delta-1)\Phi_t
      -\frac{\lambda}{2}\Wick{\Phi_t^3}+\partial_tW_t,
\end{equation}
where $W$ is a cylindrical Brownian motion on $L^2(\mathbb{R}^2)$, so that $\partial_tW$ is space-time white noise.
More specifically, we study the associated formal Kolmogorov operator
\begin{equation}\label{eq:op}
    \mathcal{L}F(\phi)
    \coloneqq \frac{1}{2}\operatorname{Tr}_{L^2(\mathbb{R}^2)}D^2F(\phi)
    -\frac{1}{2}DF(\phi)\big[(1-\Delta)\phi+\lambda\Wick{\phi^3}\big].
\end{equation}

We first recall how the measure $\mu^{(\lambda)}$, equation \eqref{eq:sq}, and the operator $\mathcal{L}$ are related through the theory of Dirichlet forms.
Let $\mathfrak{F}C_b^\infty$ denote the space of real-valued smooth cylindrical functions of the form
\begin{equation}\label{eq:cylind}
    F(\phi)=f\big(\langle\phi,l_1\rangle,\ldots,\langle\phi,l_m\rangle\big),
\end{equation}
where $m\in\mathbb{N}$, $l_1,\ldots,l_m\in C_c^\infty(\mathbb{R}^2)$, and $f\in C_b^\infty(\mathbb{R}^m)$. Here, $C_c^\infty(\mathbb{R}^2)$ and $C_b^\infty(\mathbb{R}^m)$ denote the spaces of compactly supported smooth functions and bounded smooth functions, respectively. 
We write $DF(\phi)$ for the $L^2(\mathbb{R}^2)$-gradient, characterized by
\begin{equation}
    \langle DF(\phi),l\rangle_{L^2(\mathbb{R}^2)}
    =DF(\phi)[l]
    \coloneqq\lim_{h\to0}\frac{F(\phi+hl)-F(\phi)}{h},
    \qquad l\in\mathcal{S}(\mathbb{R}^2),
\end{equation}
and define the second derivative $D^2F$ similarly.
Consider the bilinear form
\begin{equation}
    \mathcal{E}(F,G)
    \coloneqq\frac{1}{2}\int
    \langle DF,DG\rangle_{L^2(\mathbb{R}^2)}\,d\mu^{(\lambda)},
    \qquad F,G\in\mathfrak{F}C_b^\infty.
\end{equation}
Following the setting in \cite[Section~4]{RZZ17}, we take $\mu^{(\lambda)}$ to be the weak limit of the cutoff fields
\begin{equation}
    \mu^{(\lambda)}_\Lambda(d\phi)
    =\frac{1}{Z_\Lambda}\exp\left(-\frac{\lambda}{4}
      \int_\Lambda\Wick{\phi^4}(x)\,dx\right)\mu_0(d\phi)
\end{equation}
as bounded regions $\Lambda$ increase to $\mathbb{R}^2$.
Here, $Z_\Lambda$ is the normalizing constant, and the Wick power is defined with respect to $\mu_0$.
The construction and integration by parts formula established in \cite[Section~4]{RZZ17} yield
\begin{equation}\label{eq:ibp}
    \mathcal{E}(F,G)
    =-\int F\mathcal{L}G\,d\mu^{(\lambda)},
    \qquad F,G\in\mathfrak{F}C_b^\infty.
\end{equation}
In particular, the formal expression \eqref{eq:op} has a precise meaning on $\mathfrak{F}C_b^\infty$:
Writing $\mathbf{l}(\phi)=(\langle\phi,l_1\rangle,\ldots,\langle\phi,l_m\rangle)$ for a function $F$ as in \eqref{eq:cylind}, we have
\begin{align}\label{eq:op-cyl}
    \mathcal{L}F(\phi)
    &\coloneqq\frac{1}{2}\sum_{i,j=1}^m
      \partial_i\partial_jf\big(\mathbf{l}(\phi)\big)
      \langle l_i,l_j\rangle_{L^2(\mathbb{R}^2)}-\frac{1}{2}\sum_{i=1}^m
      \partial_if\big(\mathbf{l}(\phi)\big)
      \big\langle(1-\Delta)\phi+\lambda\Wick{\phi^3},l_i\big\rangle,
\end{align}
where the last pairing is the duality between $\mathcal{S}'(\mathbb{R}^2)$ and $\mathcal{S}(\mathbb{R}^2)$, consistent with the $L^2(\mathbb{R}^2)$-inner product.
The Wick power $\Wick{\phi^3}$ is well-defined under $\mu^{(\lambda)}$ and satisfies suitable weighted regularity and moment estimates; see \cite[Lemma~4.1]{RZZ17} and Section~\ref{sec:wick-phi42} below.
In particular, $\mathcal{L}F\in L^p(\mu^{(\lambda)})$ for every $F\in\mathfrak{F}C_b^\infty$ and $1\le p<\infty$.

By \eqref{eq:ibp} and the general Dirichlet form approach to SPDEs developed in \cite{AR91}, the form $(\mathcal{E},\mathfrak{F}C_b^\infty)$ is closable in $L^2(X;\mu^{(\lambda)})$, where $X=H^{-n}(w_\kappa)$ is a weighted Sobolev space with sufficiently large $n,\kappa>0$; our notation for these spaces is recalled in Section~\ref{sec:pre}.
Its closure is a quasi-regular Dirichlet form; see \cite[Section~4]{RZZ17}.
The associated diffusion admits $\mu^{(\lambda)}$ as an invariant measure and gives a probabilistically weak solution to \eqref{eq:sq}.
This realizes the stochastic quantization program of \cite{PW81} for the $\Phi^4_2$ model.
Probabilistically strong solutions were constructed on the torus in \cite{DPD03}, and global well-posedness on the plane was established in \cite{MW17}.
The identification of the solutions obtained by the Dirichlet form and SPDE approaches, in the appropriate setting, was proved in \cite{RZZ17}.
We refer to \cite{AMR15} for an overview of the Dirichlet form approach to stochastic quantization. 

Let $(\bar{\mathcal{E}},\mathcal{D}(\bar{\mathcal{E}}))$ denote the closure of $(\mathcal{E},\mathfrak{F}C_b^\infty)$ in $L^2(\mu^{(\lambda)})$.
This closed form determines a unique nonpositive self-adjoint operator $(L,\mathcal{D}(L))$, called its Dirichlet operator, and a symmetric Markov semigroup $(P_t)_{t\ge0}$; see \cite{FOT11,MR92}. 
Moreover, \eqref{eq:ibp} gives the extension
\begin{equation}\label{eq:ext}
    (\mathcal{L},\mathfrak{F}C_b^\infty)\subset(L,\mathcal{D}(L));
\end{equation}
$L$ is called the Friedrichs extension of $\mathcal{L}$.
The uniqueness of the operator associated with this \emph{fixed closed form} does not, however, imply uniqueness of the extension in \eqref{eq:ext}.
Indeed, $\mathfrak{F}C_b^\infty$ is a core (or an essential domain) of the form $\bar{\mathcal{E}}$ by construction, but need not be an core of the operator $L$: density in the form norm $\|F\|_{L^2(\mu^{(\lambda)})}+\mathcal{E}(F,F)^\frac{1}2$ is weaker than density in the graph norm $\|F\|_{L^2(\mu^{(\lambda)})}+\|LF\|_{L^2(\mu^{(\lambda)})}$.
Determining whether the action of $\mathcal{L}$ on cylindrical functions uniquely specifies the semigroup is the uniqueness problem considered here; see \cite[Chapter~1]{Ebe99} for the general framework.

For $1\le p<\infty$, the semigroup $(P_t)$ has a consistent realization as a contraction $C_0$-semigroup on $L^p(\mu^{(\lambda)})$, whose generator $L_p$ extends $(\mathcal{L},\mathfrak{F}C_b^\infty)$. More precisely, the restriction $P_t \big|_{L^2\cap L^p}$ is uniquely extended to a contraction $C_0$-semigroup on $L^p(\mu^{(\lambda)})$ and its generator $L_p$ extends $(\mathcal{L},\mathfrak{F}C_b^\infty)$. 
In particular, this pre-Dirichlet operator is densely defined, dissipative, and closable on $L^p(\mu^{(\lambda)})$; write $\bar{\mathcal{L}}_p$ for its closure.
The operator is said to be \emph{$L^p(\mu^{(\lambda)})$-unique} if there is only one $C_0$-semigroup on this space whose generator extends $(\mathcal{L},\mathfrak{F}C_b^\infty)$.
In the present setting, this is equivalent to $\bar{\mathcal{L}}_p=L_p$, or, equivalently, to $\mathfrak{F}C_b^\infty$ being a core of $L_p$; see \cite[Chapter~1]{Ebe99}.

The connection with maximal dissipativity is provided by the Lumer--Phillips theorem.
Since dissipativity is already available, $\bar{\mathcal{L}}_p$ generates a contraction $C_0$-semigroup precisely when it is maximal dissipative, that is, when $a-\bar{\mathcal{L}}_p$ is surjective for some $a>0$.
Equivalently, it suffices to establish the dense-range condition
\begin{equation}\label{eq:cond-dense}
    \overline{(a-\mathcal{L})(\mathfrak{F}C_b^\infty)}^{\,L^p(\mu^{(\lambda)})}
    =L^p(\mu^{(\lambda)})
\end{equation}
for some $a>0$; the range of $a-\bar{\mathcal{L}}_p$ is closed by dissipativity and closedness.
Once this holds, $\bar{\mathcal{L}}_p=L_p$.
Moreover, a $C_0$-semigroup generator cannot have a proper extension that is itself a $C_0$-semigroup generator, so uniqueness holds among \emph{all} $C_0$-semigroups, not merely contraction semigroups.
Thus, proving maximal dissipativity of the closure, also called \emph{essential maximal dissipativity} of $(\mathcal{L},\mathfrak{F}C_b^\infty)$, resolves the $L^p$-uniqueness problem; see \cite[Chapter~1]{Ebe99}.

For $p=2$, symmetry implies that $L^2$-uniqueness is equivalent to essential self-adjointness.
This is stronger than \emph{Markov uniqueness}, which asks for uniqueness only among self-adjoint extensions generating Markov $C_0$-semigroups.
$L^1$-uniqueness also implies Markov uniqueness, since such semigroups extend consistently to contraction $C_0$-semigroups on $L^1(\mu^{(\lambda)})$; see \cite[Chapter~1]{Ebe99}.

Uniqueness problems for Dirichlet operators have been studied extensively in stochastic quantization.
For the infinite-volume $P(\Phi)_1$-models, essential self-adjointness and $L^p$-uniqueness were proved in \cite{KR07} and \cite{AKR12}, respectively.
In two dimensions, earlier results concerned \emph{modified} stochastic quantization equations with a smoothing parameter; see \cite{LR98,DT00}.
For the quartic interaction, the corresponding modification takes the form
\begin{equation}\label{eq:msq}
    \partial_t\Phi_t
    =-\frac{1}{2}(1-\Delta)^{1-\varepsilon}\Phi_t
      -\frac{\lambda}{2}(1-\Delta)^{-\varepsilon}\Wick{\Phi_t^3}
      +(1-\Delta)^{-\varepsilon/2}\partial_tW_t,
\end{equation}
with $0<\varepsilon<1$, whereas $\varepsilon=0$ gives \eqref{eq:sq}.
Related uniqueness results for modified exponential and trigonometric models, together with a detailed account of the literature, can be found in \cite{AKMR23}.
For the unmodified $\Phi^4_2$-model on the torus, maximal dissipativity in $L^p$, $1\le p<\infty$, was established in \cite[Theorem~4.1 and Remark~4.2]{DD20}.
The present paper extends this result to the plane.
In the infinite-volume setting, \cite{RZZ17} established a restricted form of Markov uniqueness within a class of quasi-regular Dirichlet forms.
Our result removes that restriction and yields the stronger conclusion of $L^p$-uniqueness for every $1\le p<\infty$. 

In three dimensions, significant progress has been made in the study of the stochastic quantization of $\Phi^4_3$ model in recent years thanks to advances in the field of singular SPDE made by \cite{Hai14,GIP15}; see \cite{HM18,CC18,MW17-phi43,AK20,AK22,GH19,GH21,DHYZ25}.
However, there have been few studies that use the theory of Dirichlet forms; the $\Phi^4_3$ model does not fit within the framework of \cite{AR91} as pointed out in \cite{ALZ06}. Nevertheless, some progress was made in \cite{ZZ18} and the authors constructed the Dirichlet form associated to $\Phi^4_3$ model on the torus using techniques from singular SPDE theory.

A key ingredient in \cite{DD20} is a gradient estimate for the Markov semigroup of the dynamical $\Phi^4_2$-model on the torus, of the form
\begin{equation}\label{eq:gradient}
    |DP_tF(\phi)[h]|
    \lesssim t^{-\beta}\|h\|_{\mathcal{C}^{-s}}\|F\|_\infty
      \big(1+\|\phi\|_{\mathcal{C}^{-\alpha}}\big)^\gamma,
    \qquad 0<t\le1,
\end{equation}
for suitable positive exponents $\alpha,\beta,\gamma,s$.
Such estimates are closely related to control of the linearized equation
\begin{equation}\label{eq:intro-J}
    \begin{cases}
        \displaystyle\partial_tJ_t
        =\frac{1}{2}(\Delta-1)J_t
         -\frac{3\lambda}{2}\Wick{\Phi_t^2}J_t,\\[2pt]
        J_0=h,
    \end{cases}
\end{equation}
which formally describes the derivative of the solution with respect to its initial condition in the direction $h$.
The renormalized coefficient $\Wick{\Phi_t^2}$ is not nonnegative, and a direct Gronwall estimate leads to exponential factors whose integrability is not provided by the available moment bounds.
In \cite{DD20}, this difficulty is handled using an auxiliary Feynman--Kac semigroup, a Duhamel formula, and a Bismut--Elworthy--Li type formula.
Here, instead of relying on an infinite-volume analogue of the strong Feller-type estimate \eqref{eq:gradient}, we adapt the stopping time argument of \cite{KT25} to obtain weighted estimates for the linearized equation.
Combined with suitable approximations of the resolvent, these estimates yield the dense-range condition~\eqref{eq:cond-dense} required for maximal dissipativity.

We can now state our main result.
Recall that $\bar{\mathcal{L}}_p$ denotes the closure in $L^p(\mu^{(\lambda)})$ of the operator defined by \eqref{eq:op-cyl} on $\mathfrak{F}C_b^\infty$.

\begin{thm}\label{thm:main}
    For every $\lambda>0$ and $1\le p<\infty$, the operator
    $(\bar{\mathcal{L}}_p,\mathcal{D}(\bar{\mathcal{L}}_p))$ is maximal dissipative in $L^p(\mu^{(\lambda)})$.
    Equivalently, it is dissipative and $a-\bar{\mathcal{L}}_p$ is surjective for some, and hence every, $a>0$.
\end{thm}

\begin{rem}\label{rem:thm11}
    In our proof, we first establish surjectivity for sufficiently large $a>0$. Then this already implies maximal dissipativity and hence surjectivity for every $a>0$, without any smallness assumption on $\lambda$.
\end{rem}

\begin{rem}\label{rem:thm11-pq}
 For $1\le p\le q <\infty$, the maximal dissipatity of $(\bar{\mathcal{L}}_q,\mathcal{D}(\bar{\mathcal{L}}_q))$ implies that of $(\bar{\mathcal{L}}_p,\mathcal{D}(\bar{\mathcal{L}}_p))$ since the embedding $L^q\hookrightarrow L^p$ is continuous and dense.

\end{rem}

The Lumer--Phillips theorem and the uniqueness criterion recalled above give the following consequence.

\begin{cor}[$L^p(\mu^{(\lambda)})$-uniqueness]
    Let $1\le p<\infty$.
    The operator $(\bar{\mathcal{L}}_p,\mathcal{D}(\bar{\mathcal{L}}_p))$ generates a contraction $C_0$-semigroup $(P_t^{(p)})_{t\ge0}$ on $L^p(\mu^{(\lambda)})$.
    This is the only $C_0$-semigroup on $L^p(\mu^{(\lambda)})$ whose generator extends $(\mathcal{L},\mathfrak{F}C_b^\infty)$, and it coincides with the $L^p$-realization of the semigroup associated with $(\bar{\mathcal{E}},\mathcal{D}(\bar{\mathcal{E}}))$, which is mentioned above.
\end{cor}

In particular, the $L^2$-case identifies the closure of the cylindrical operator with the Dirichlet operator $L$ of the closed form $(\bar{\mathcal{E}},\mathcal{D}(\bar{\mathcal{E}}))$.

\begin{cor}[Essential self-adjointness and Markov uniqueness]
    Let $L$ and $(P_t)_{t\ge0}$ be the Dirichlet operator and the Markov semigroup associated with $(\bar{\mathcal{E}},\mathcal{D}(\bar{\mathcal{E}}))$, respectively.
    Then $L=\bar{\mathcal{L}}_2$.
    Consequently, $(\mathcal{L},\mathfrak{F}C_b^\infty)$ is essentially self-adjoint in $L^2(\mu^{(\lambda)})$, and $\mathfrak{F}C_b^\infty$ is an operator core for $L$.
    Moreover, $(P_t)$ is the only Markov $C_0$-semigroup on $L^2(\mu^{(\lambda)})$ whose generator extends $(\mathcal{L},\mathfrak{F}C_b^\infty)$.
\end{cor}

The paper is organized as follows.
Section~\ref{sec:pre} collects the notation and preliminary results on weighted function spaces and the heat semigroup.
Theorem~\ref{thm:main} is proved in Section~\ref{sec:proof}, using technical results established in Sections~\ref{sec:jacobi} and ~\ref{sec:resolvents}.
Section~\ref{sec:jacobi} contains the key estimates for the linearized equation, while Section~\ref{sec:resolvents} treats the resolvent approximations.




\section{Preliminaries}\label{sec:pre}

\subsection{Weighted function spaces}
We use the unitary Fourier transform
\[
 \widehat f(\xi)=(2\pi)^{-d/2}\int_{\mathbb R^d}e^{-ix\cdot\xi}f(x)\,dx,
 \qquad
 \mathcal F^{-1}f(x)=(2\pi)^{-d/2}\int_{\mathbb R^d}e^{ix\cdot\xi}f(\xi)\,d\xi.
\]
Let $\chi, \rho \in C_c ^\infty (\mathbb{R}^d; [0,\infty))$ be such that
\[ \mathrm{supp}(\chi) \subset{B(4)},  \ \mathrm{supp}(\rho) \subset{B(4)\backslash B(1)}\ \mbox{and}\ \chi (x) + \sum_{i=0}^\infty \rho (2^{-i}x) = 1 \ \ \mbox{for any}\ x \in \mathbb{R}^d,            \]
and write
\[ \rho_{-1} = \chi, \ \rho_{j} = \rho(2^{-j}\cdot) \ \ \ \mathrm{for} \ j \ge 0. \]
The associated Littlewood-Paley blocks are then defined by
\[ \Delta_j f(x) \coloneqq \hat{\rho}_j*f \ \ \mathrm{for} \ f \in \mathcal{S}'(\mathbb{R}^d).\]
For a positive weight $w$, we write
\[
 \|f\|_{L^p(w)}=\|wf\|_{L^p},\qquad
 \|f\|_{B^s_{p,q}(w)}=
 \big\|\big(2^{js}\|w\Delta_jf\|_{L^p}\big)_{j\ge-1}\big\|_{\ell^q}.
\]
Corresponding Besov spaces are then defined by the completion of $C_c^\infty$ in these norms. 
We set $\mathcal C^s(w)=B^s_{\infty,\infty}(w)$ and $H^s(w)=B^s_{2,2}(w)$.
When $w \equiv 1$, we omit the weight and write $B_{p,q}^s$.

The weights used in this paper are
\begin{equation}\label{eq:weights}
 w_\kappa(x)=\langle x\rangle^{-\kappa},\qquad
 e_\delta(x)=\exp(-\delta^{1/2}\langle x\rangle^{\theta}),
 \qquad \langle x\rangle=(1+|x|^2)^{1/2},
\end{equation}
where $\kappa\in\mathbb R$, $\delta>0$ and $0<\theta<1$.
These weights satisfy the following admissible conditions
\begin{align}
 \frac{w_\kappa(x)}{w_\kappa(y)}
 &\lesssim_\kappa\langle x-y\rangle^{|\kappa|},\label{eq:adm-w}\\
 \frac{e_\delta(x)}{e_\delta(y)}
 &\le \exp\big(\delta^{1/2}|x-y|^{\theta}\big).
 \label{eq:adm-e}
\end{align}
When we work with the subexponential weights $e_\delta$ as above, it is often useful to choose various cutoffs in a compactly supported Gevrey class of order $\sigma$ with $1 < \sigma < 1/\theta$ so that their inverse Fourier transforms decay faster than the rates in \eqref{eq:adm-e}, see \cite{MW17}. Indeed, we assume that our Littlewood-Paley cutoffs $\chi,\rho$ are in this class.

\begin{lem}\label{lem:weight-equiv}
For $s,\kappa\in\mathbb R$ and $1\le p,q\le\infty$,
\[
 \|f\|_{B^s_{p,q}(w_\kappa)}\simeq
 \|w_\kappa f\|_{B^s_{p,q}}.
\]
\end{lem}
\begin{proof}
See \cite[Theorem~6.5]{Tri06}.
\end{proof}

For $1<p<\infty$, we also use the Sobolev norm
\begin{equation}\label{eq:sobolev}
 \|f\|_{W^{s,p}(w_\kappa)}:=\|w_\kappa f\|_{W^{s,p}}.
\end{equation}
The usual Besov--Sobolev embeddings and Lemma~\ref{lem:weight-equiv}
give
\begin{equation}\label{eq:embed-sobolev}
 \|f\|_{B^s_{p,\infty}(w_\kappa)}
 \lesssim\|f\|_{W^{s,p}(w_\kappa)}
 \lesssim\|f\|_{B^s_{p,1}(w_\kappa)};
\end{equation}
see \cite{Tri83,BCD11}. 

\begin{lem}[Duality]\label{lem:duality}
Let $s,\kappa\in\mathbb R$, $1/p_1+1/p_2=1$, and $1/q_1+1/q_2=1$.
Then
\[
 |\langle f,g\rangle|
 \lesssim\|f\|_{B^s_{p_1,q_1}(w_\kappa)}
             \|g\|_{B^{-s}_{p_2,q_2}(w_{-\kappa})}.
\]
\end{lem}
\begin{proof}
It follows by applying unweighted Besov duality to $w_\kappa f$ and $w_{-\kappa}g$,
and using Lemma~\ref{lem:weight-equiv}; see \cite[Proposition~2.76]{BCD11}.
\end{proof}

\begin{lem}[Product estimates]\label{lem:product}
Let $a<0<b$, $a+b>0$, $1/p=1/p_1+1/p_2$, and
$1/q=1/q_1+1/q_2$. Then
\[
 \|fg\|_{B^{2a}_{p,q}(w_{\kappa_1+\kappa_2})}
 \lesssim\|f\|_{B^a_{p_1,q_1}(w_{\kappa_1})}
             \|g\|_{B^b_{p_2,q_2}(w_{\kappa_2})}.
\]
We shall also use the sharper special case: if $0<a<b$, then
\begin{equation}\label{eq:rough-product}
 \|fg\|_{B^{-a}_{p,\infty}(w_{\kappa_1+\kappa_2})}
 \lesssim\|f\|_{\mathcal C^{-a}(w_{\kappa_1})}
             \|g\|_{B^b_{p,\infty}(w_{\kappa_2})}.
\end{equation}
For $0<b<1$ and an integer $k\ge1$, we have the fractional Leibniz bound
\begin{equation}\label{eq:power-product}
 \|v^k\|_{B^b_{p,\infty}(w_{\omega+(k-1)\tau})}
 \lesssim\|v\|_{B^b_{q,\infty}(w_\omega)}
             \|v\|_{L^r(w_\tau)}^{k-1}
\end{equation}
when $\omega\ge\tau$ and $1/p=1/q+ (k-1)/r$. 
\end{lem}
\begin{proof}
The first two assertions follow from Bony's decomposition, the usual
paraproduct estimates, and Lemma~\ref{lem:weight-equiv}; see
\cite[Chapter~2]{BCD11}. 
For \eqref{eq:power-product}, the difference characterization of
$B^b_{p,\infty}$ and Lemma~\ref{lem:weight-equiv} give
\begin{align}
 & \|v^k\|_{B^b_{p,\infty}(w_{\omega+(k-1)\tau})} \simeq \|v^k w_{\omega+(k-1)\tau}\|_{L^p} \notag\\
 &\qquad +\sup_{t>0}  t^{-b}  \sup_{|h|\le t} \left\|   (v^k w_{\omega+(k-1)\tau}) (\cdot+ h) - \left( v^k w_{\omega+(k-1)\tau}\right) (\cdot)    \right\|_{L^p}  . 
\end{align}
The first term can be treated easily using H\"older's inequality. For the second term, one has for $|h|\le 1$ that
\begin{align*}
    |w_{\omega+(k-1)\tau} (x+h) | \simeq |w_{\omega+(k-1)\tau}(x) |,
\end{align*} 
\begin{align*}
    |w_{\omega+(k-1)\tau}(x+h) - w_{\omega+(k-1)\tau}(x) | \lesssim |h| |w_{\omega+(k-1)\tau} (x)|,
\end{align*}
\begin{align*}  
&(v^k w_{\omega+(k-1)\tau}) (x+ h) - \left( v^k w_{\omega+(k-1)\tau}\right) (x) \notag\\
&= \left\{ v(x+h)^k-v(x)^k \right\}w_{\omega+(k-1)\tau}(x+h) + v(x)^k \left\{ w_{\omega+(k-1)\tau}(x+h) - w_{\omega+(k-1)\tau}(x)  \right\} ,
\end{align*}
and
$$v(x+h)^k-v(x)^k=(v(x+h)-v(x))\sum_{j=0}^{k-1}v(x+h)^jv(x)^{k-1-j}.$$
Therefore, H\"older's inequality and the difference characterization of
$B^b_{p,\infty}$ again yield the conclusion. 
\end{proof}

The following product estimate will be used only in the proof of Lemma~\ref{lem:J-Lp}.
\begin{lem}\label{lem:product-special}
    Let $a>0, \kappa\in \mathbb{R}, \delta>0$ and $1\le p\le \infty$. Then
    \begin{align}
    \|fg\|_{B^{-2a}_{p,p}(w_{\kappa}e_\delta)}
 \lesssim\|f\|_{\mathcal{C}^{-a}(w_{\kappa})}
             \|g\|_{B^{2a}_{p,p}(e_\delta)}.    
    \end{align}
\end{lem}
\begin{proof}
    The claim follows from \cite[Lemma~3.14]{DHYZ25}.
\end{proof}

\begin{lem}[Interpolation]\label{lem:interpolation}
Suppose $0\le\vartheta\le1$,
\[
 \frac1{p_0}=\frac{\vartheta}{p_1}+\frac{1-\vartheta}{p_2},\quad
 \frac1{q_0}=\frac{\vartheta}{q_1}+\frac{1-\vartheta}{q_2},\quad
 s_0=\vartheta s_1+(1-\vartheta)s_2,\quad
 \kappa_0=\vartheta\kappa_1+(1-\vartheta)\kappa_2.
\]
Then
\begin{equation}\label{eq:interpolation-Besov}
 \|f\|_{B^{s_0}_{p_0,q_0}(w_{\kappa_0})}
 \lesssim\|f\|_{B^{s_1}_{p_1,q_1}(w_{\kappa_1})}^{\vartheta}
             \|f\|_{B^{s_2}_{p_2,q_2}(w_{\kappa_2})}^{1-\vartheta}.
\end{equation}
In particular, $B^0_{r,\infty}$ on the right may be bounded by $L^r$.
Also, for $0<a<b<1$,
\begin{equation}\label{eq:L1-interpolation}
 \|h\|_{B^a_{1,1}}
 \lesssim\|h\|_{L^1}^{1-b}\|h\|_{W^{1,1}}^b.
\end{equation}
\end{lem}
\begin{proof}
H\"older's inequality gives
\begin{align}
    2^{js_0} \| w_{\kappa_0} \Delta_j f\|_{L^{p_0}} &= 2^{js_0}  \| \left| w_{\kappa_1} \Delta_j f \right|^\vartheta \left| w_{\kappa_2} \Delta_j f \right|^{1-\vartheta}    \|_{L^{p_0}} \le \left( 2^{js_1} \|  \Delta_j f\|_{L^{p_1}(w_{\kappa_1})}  \right)^\vartheta \left( 2^{js_2} \|  \Delta_j f\|_{L^{p_2}(w_{\kappa_2})}  \right)^{1- \vartheta}.
\end{align}
Applying H\"older's inequality to the dyadic sequence, we get \eqref{eq:interpolation-Besov}. For \eqref{eq:L1-interpolation}, apply \eqref{eq:interpolation-Besov} and \eqref{eq:embed-sobolev}.
See also \cite[Chapter~2]{BCD11}.
\end{proof}

We use the following lemma in Section~\ref{sec:CDFI}.

\begin{lem}[Embeddings and localization]\label{lem:weighted-embeddings}
If $1\le r\le p\le\infty$, then
\begin{equation}\label{eq:weighted-sobolev-embedding}
 B^s_{r,q}(w_\kappa)\hookrightarrow
 B^{s-d(1/r-1/p)}_{p,q}(w_\kappa).
\end{equation}
If $1\le p\le r\le\infty$ and
$\kappa-\kappa_0>d(1/p-1/r)$, then
\begin{equation}\label{eq:weighted-holder-embedding}
 B^s_{r,q}(w_{\kappa_0})\hookrightarrow B^s_{p,q}(w_\kappa).
\end{equation}
The corresponding weighted Lebesgue embedding also holds.
Finally, if $\chi_z(x)=\chi_*(x-z)$ for a fixed $\chi_* \in C_c^\infty(\mathbb{R}^d)$ and $z\in\mathbb Z^d$, then, for $\zeta\ge0$,
\begin{equation}\label{eq:unit-localization}
 \|\chi_z f\|_{\mathcal C^s}
 \lesssim\langle z\rangle^\zeta\|f\|_{\mathcal C^s(w_\zeta)}.
\end{equation}
\end{lem}
\begin{proof}
The first assertion follows from the unweighted Besov embedding and
Lemma~\ref{lem:weight-equiv}. For the second, H\"older's inequality gives
\begin{align}
    \| \Delta_j f \|_{L^p(w_\kappa)} = \| \Delta_j f w_{\kappa_0} w_{\kappa-\kappa_0}\|_{L^p} \le \| \Delta_j f \|_{L^r(w_{\kappa_0})}\| w_{\kappa-\kappa_0}\|_{L^{\frac{1}{1/p-1/r}}} \lesssim \| \Delta_j f \|_{L^r(w_{\kappa_0})}.
\end{align}
This implies \eqref{eq:weighted-holder-embedding}.
For the last assertion,
Lemma~\ref{lem:weight-equiv} gives 
\begin{align}
    \| \chi_z f \|_{\mathcal{C}^s} = \| (\chi_z w_{-\zeta}) (w_\zeta f )  \|_{\mathcal{C}^s} \lesssim \| \chi_z w_{-\zeta}  \|_{\mathcal{C}^s} \| f \|_{\mathcal{C}^s(w_\zeta)} \lesssim \langle z \rangle^\zeta \| f \|_{\mathcal{C}^s(w_\zeta)}
\end{align}
since every derivatives of $\chi_z w_{-\zeta}$ is bounded by $C\langle z \rangle^\zeta$.
\end{proof}

\subsection{Heat semigroup and approximation operator}
Let $e^{\frac t2(\Delta-1)}$ be the heat semigroup associated to $\frac{1}2(\Delta-1)$.
\begin{lem}[Schauder estimate]\label{lem:schauder}
Let $w$ be a weight in \eqref{eq:weights}, $s\in\mathbb R$, $\varepsilon\ge0$,
$1\le p,q\le\infty$, and $T<\infty$. Then
\[
 \|e^{\frac t2(\Delta-1)} f\|_{B^{s+\varepsilon}_{p,q}(w)}
 \le C_Te^{-t/2} t^{-\varepsilon/2}\|f\|_{B^s_{p,q}(w)},
 \qquad 0<t\le T.
\]
\end{lem}
\begin{proof}
See \cite[Section~3]{MW17}.
\end{proof}

Choose an even $g\in C_c^\infty(\mathbb R^d; \mathbb{R})$ with $g=1$
on the unit ball, and set $g_N(\xi)=g(\xi/N)$. We also assume that $g$ is in the same Gevrey class as above. 
Define
\begin{equation}\label{eq:mollify}
 \Pi_N f=\mathcal F^{-1}(g_N\widehat f)=K_N*f,
 \qquad K_N(x)=(2\pi)^{-d/2}N^2\check g(Nx).
\end{equation}
Here $\check g=\mathcal F^{-1}g$. 

\begin{lem}\label{lem:appr-1}
For $N\ge1$, $\kappa,s\in\mathbb R$, and $1\le p,q\le\infty$,
\begin{align}
 \|\Pi_N f\|_{L^p(w_\kappa)}&\lesssim\|f\|_{L^p(w_\kappa)},\notag\\
 \|\Pi_N f\|_{B^s_{p,q}(w_\kappa)}&\lesssim\|f\|_{B^s_{p,q}(w_\kappa)}.
 \label{eq:uniform-mollifier}
\end{align}
For $\varepsilon>0$,
\begin{equation}\label{eq:mollifier-rate}
 \|(I-\Pi_N)f\|_{B^s_{p,q}(w_\kappa)}
 \lesssim N^{-\varepsilon}\|f\|_{B^{s+\varepsilon}_{p,q}(w_\kappa)}.
\end{equation}
The constants do not depend on $N$. The analogous statements hold for the subexponential weight $e_\delta$.
\end{lem}
\begin{proof} 
From the admissible condition \eqref{eq:adm-w} and Young's convolution inequality, we have
\begin{align}
 \| \Pi_N f \|_{L^p(w_\kappa)} &= \left( \int_{\mathbb{R}^2} \left| K_N*f(x)   \right|^p w_\kappa^{p}(x) \, dx   \right)^{\frac{1}p} \notag\\
 &\lesssim \left\| (K_Nw_{|\kappa|})*(fw_\kappa)   \right\|_{L^p} \notag\\
 &\le \left\| K_Nw_{|\kappa|}   \right\|_{L^1} \left\| fw_\kappa   \right\|_{L^p}.
\end{align}
Together with
\[
 \int |K_N(x)|\langle x\rangle^{|\kappa|}\,dx
 \lesssim\int |\check g(y)|\langle y\rangle^{|\kappa|}\,dy<\infty,
\]
this implies the first assertion.
Since $\Pi_N$
commutates with $\Delta_j$, the second assertion also follows.
Moreover, $(I-\Pi_N)\Delta_j=0$ unless $2^j\gtrsim N$. Therefore, \eqref{eq:uniform-mollifier} implies 
\begin{align}
 & 2^{js_1} \| (I - \Pi_N) \Delta_j f \|_{L^p(w_\kappa)}  \lesssim N^{-\varepsilon}  2^{j(s+\varepsilon)} \| (I - \Pi_N) \Delta_j f \|_{L^p(w_\kappa)} \lesssim N^{-\varepsilon}  2^{j(s+\varepsilon)} \|  \Delta_j f \|_{L^p(w_\kappa)},
 \end{align}
from which \eqref{eq:mollifier-rate} follows.

The case for $e_\delta$ can be proven in the same manner.

\end{proof}

\section{Proof of Theorem~\ref{thm:main}}\label{sec:proof}
In this section, we prove Theorem~\ref{thm:main}, deferring many of the technical results to later sections.
We consider the following approximation of the stochastic quantization equation \eqref{eq:sq}:
\begin{align}
    \begin{cases}
        \partial_t \Phi_t^{(N)}= \frac{1}2(\Delta-1)\Phi_t^{(N)}- \frac{\lambda}2  \left\{ \left(  \Phi_t^{(N)}  \right)^3 - 3C_N\Phi_t^{(N)}\right\} + \Pi_N \partial_t W_t\\
        \Phi_0^{(N)}= \phi, \label{eq:SQ}
    \end{cases}
\end{align}
where $W_t$ is a cylindrical Wiener process on $L^2(\mathbb{R}^2)$, $\Pi_N$ is the mollification operator defined in \eqref{eq:mollify}, 
 and $C_N$ is the renormalization constant defined by \eqref{eq:renor-const}.
The associated Markov semigroup is defined by
\begin{equation}\label{eq:defP}
    P_t^{(N)}F(\psi) \coloneqq \mathbb{E} \left[ F\left( \Phi_t^{(N)}(\psi) \right)  \right]
\end{equation}
and, for any given $F \in \mathfrak{F}C_b^\infty$, we consider the time-truncated resolvent function
\begin{equation}\label{eq:def:G}
 G^{(N)}_T(\psi) = \int^T_0 e^{-at}  P_t^{(N)}F(\psi)\, dt,
\end{equation}
where $\mathfrak{F}C_b^\infty$ is the space of all $\mathbb{R}$-valued cylindrical functions of the form
\begin{equation}
    F(\phi) = f\left( \langle \phi, l_1\rangle, \cdots \langle \phi, l_m \rangle   \right)
\end{equation}
for some $m\in \mathbb{N}, l_1,\cdots, l_m \in C_c^\infty(\mathbb{R}^2)$ and $f \in C_b^\infty(\mathbb{R}^m)$.
We then consider the corrected resolvent function
\begin{equation}\label{eq:def:Gtilde}
 \tilde{G}^{(N)}_T(\phi) \coloneqq G^{(N)}_T(\Pi_N \phi) = \int^T_0 e^{-at}  P_t^{(N)}F( \Pi_N \phi)\, dt. 
\end{equation}
Our strategy for proving Theorem~\ref{thm:main} is to approximate the resolvent of $F \in \mathfrak{F}C_b^\infty$ by the sequence $\{\tilde{G}_T^{(N)}\}_{T,N\in \mathbb{N}}$. 
To implement this approach, we use the following result, whose proof is given in Section~\ref{sec:resolvents}.

\begin{lem}\label{lem:resolvent}
Let $1\le p<\infty$, $N\in \mathbb{N}$ and $T,a>0$. Then  
$\tilde{G}^{(N)}_T \in \mathcal{D} \left(  \bar{\mathcal{L}_p} \right)$ and
    \begin{align}
     (a- \bar{\mathcal{L}}_p)\tilde{G}_T^{(N)}(\phi)
       &= F(\Pi_N \phi) - e^{-aT} P_T^{(N)} F(\Pi_N \phi ) \notag\\
       &\quad + \frac{\lambda}2 DG_T^{(N)} (\Pi_N \phi)
          [  \Pi_N \Wick{\phi^3} - \Wick{(\Pi_N\phi)^3}_{C_N}  ] . \label{eq:LG}
    \end{align}
\end{lem}

It is also proved in Section~\ref{sec:resolvents} that for any $h$ smooth enough,
\begin{align}
    DG^{(N)}_T(\phi) [h] = \int^T_0 e^{-at}\mathbb{E} \left[  DF(\Phi_t^{(N)}(\phi)  )\left[ J_t^h(\phi)\right]  \right]\, dt, \label{eq:DG1}
\end{align}
where $J_t^h(\phi) \coloneqq  D\Phi_t^{(N)}(\phi)[h] $ solves the linearized equation
\begin{align}
    \begin{cases}
        \partial_t J_t =  \frac{1}2 (\Delta -1)J_t - \frac{3\lambda}2 \left\{ \left(  \Phi_t^{(N)}  \right)^2 - C_N\right\}  J_t      \\
        J_0=h. \label{eq:def-JJ}
    \end{cases}
\end{align}

\begin{rem}
One might think that the analogous statement to Lemma~\ref{lem:resolvent} in the case $T=\infty$, 
\begin{align}
     (a- \bar{\mathcal{L}}_p)\tilde{G}_\infty^{(N)}(\phi)
       &= F(\Pi_N \phi)  \notag\\
       &\quad + \frac{\lambda}2 DG_\infty^{(N)} (\Pi_N \phi)
          [  \Pi_N \Wick{\phi^3} - \Wick{(\Pi_N\phi)^3}_{C_N}  ], 
    \end{align}
can be proved directly, but the introduction of the time-truncated resolvents as above turns out to be very useful. We will take the limit $T\to \infty$ by combining the identities~\eqref{eq:LG} and \eqref{eq:DG1} with the bound \eqref{eq:est-Jlinear} below.
In particular, we do not need to establish any bounds for the second derivative that are uniform in $N$.

\end{rem}

The key ingredient in the proof of Theorem~\ref{thm:main} is the following estimate on the solution to \eqref{eq:def-JJ}. The estimate is proved in Section~\ref{sec:jacobi}.

\begin{lem}\label{lem:linear}
  Fix any $\lambda_*>0$.  Assume that $0<\lambda \le \lambda_*, \theta, \eta, \delta, \varepsilon, \zeta, \kappa >0$ and $p\in 2\mathbb{N}$ satisfy
  $2p\delta\theta<1$, $3p\varepsilon + \frac{p(1+\varepsilon)}{2(p-1)(1-\varepsilon)}<1$, $\kappa\le \theta/2$, $\zeta>3p\varepsilon$, $0<\theta<1$ and $6p(\varepsilon+\kappa)<1$.
  Then there exist a polynomial $\mathfrak{p}$ and constants $c, \gamma, a, b>0$ independent of $\lambda$ such that 
  \begin{equation}
    \begin{aligned}
      &\mathbb{E} \left[ \left\| J^h_t (\phi) \right\|_{L^{p}(e_{t+\delta})}^{p}  \right] \le (1+ t^{-\zeta}) \| h\|_{\mathcal{C}^{-\varepsilon}(w_\kappa)}^{p}
        \mathfrak{p} (\| \phi \|_{\mathcal{C}^{-a}(w_b)})
        e^{(c\lambda^\gamma + \eta +3p\theta - p/2)t}.
    \end{aligned}\label{eq:est-Jlinear}
  \end{equation}
\end{lem}

\begin{rem}
As already mentioned in Section~\ref{sec:intro}, the difficulty in dealing with the linearized equation \eqref{eq:def-JJ} is mainly due to the involvement of renormalization. 
In the proof of Lemma~\ref{lem:linear}, we apply the method introduced in \cite{KT25}, where a similar estimate was proved on the two-dimensional torus. This estimate was later extended to $\mathbb{R}^3$ in \cite{DHYZ25}.
However, these works treat the linearized equation \eqref{eq:def-JJ} in settings different from ours, and their estimates assume initial data $h \in L^p$, rather than in a Besov space of negative regularity as in Lemma~\ref{lem:linear}. To deal with this difference, we decompose $J_t^h$ by subtracting the contribution from the initial value $h$ and estimate the resulting processes separately, see Lemma~\ref{lem:J-C}.
\end{rem}

\begin{rem}
    Note that the weight function in the left-hand side of \eqref{eq:est-Jlinear} depends on the time variable $t$.
It should be possible to prove a similar bound for some weight that does not depend on $t$, for example by applying a technique from \cite[Theorem 2.3]{DHYZ25}, but the current claim is sufficient for our purpose.  
\end{rem}

We now prove Theorem~\ref{thm:main}, assuming Lemmas~\ref{lem:resolvent} and \ref{lem:linear}.

\begin{proof}[Proof of Theorem~\ref{thm:main}]
We begin with fixing parameters that satisfy the assumptions of Lemma~\ref{lem:linear}. 
We can assume $\frac{p}{2(p-1)}< \frac{2}3$ in view of Remark~\ref{rem:thm11-pq}. 
Then we fix sufficiently small $\varepsilon,\theta,\eta>0$ such that $\eta<4/5$ and
\begin{align*}
    3p\varepsilon<\frac{1}{10}, \qquad\frac{1+\varepsilon}{1-\varepsilon}< \frac{6}5, \qquad 3p\theta< \eta.
\end{align*}
Moreover, we set $\kappa=\theta/2$ and take any $0<\delta<1$, $\frac{1}{10}<\zeta<1$.

Fix $F \in \mathfrak{F}C_b^\infty$ and write $h_N(\phi) \coloneqq \Pi_N\left(\Wick{\phi^3}\right)- \Wick{(\Pi_N \phi)^3}$. Lemma~\ref{lem:resolvent} and \eqref{eq:DG1} give 
\begin{align}
&(a- \bar{\mathcal{L}}_p)\tilde{G}^{(N)}_T(\phi) - F(\Pi_N \phi) - e^{-aT}P_T^{(N)}F(\Pi_N \phi)  \notag\\
&=  \frac{\lambda}2  DG^{(N)}_T(\Pi_N \phi) \left[  h_N(\phi) \right]   \notag\\
&=  \frac{\lambda}2 \int^T_0 e^{-at}\mathbb{E} \left[  DF(\Phi_t^{(N)}(\Pi_N\phi)  )\left[ J_t^{h_N(\phi)}(\Pi_N\phi) \right] \right]\, dt \notag\\
&=  \frac{\lambda}2 \sum_{i=1}^m \int^T_0 e^{-at}\mathbb{E} \left[ (\partial_i f)\left(\{\langle l_j , \Phi_t^{(N)}(\Pi_N\phi) \rangle\}_{1\le j\le m} \right) \langle l_i, J_t^{h_N(\phi)}(\Pi_N\phi) \rangle \right]\, dt.
\end{align}
Therefore, H\"older's inequality, Lemma~\ref{lem:linear}, and $f\in C_b^\infty(\mathbb{R}^m)$ yield 
\begin{align}
   & \left\| (a- \bar{\mathcal{L}}_p)\tilde{G}^{(N)}_T(\cdot) - F(\Pi_N \cdot) - e^{-aT}P_T^{(N)}F(\Pi_N\cdot)  \right\|_{L^p(\mu^{(\lambda)})} \notag\\
    &\lesssim \sum_{i=1}^m \int^T_0 e^{-at} \left\| \mathbb{E} \left[ \left| \langle l_i, J_t^{h_N(\cdot)}(\Pi_N\cdot) \rangle  \right|  \right]   \right\|_{L^p(\mu^{(\lambda)})} \notag\\
    &= \sum_{i=1}^m \int^T_0 e^{-at} \left\{ \int \| l_i  \|_{L^{\frac{p}{p-1}}(e_{t+\delta}^{-1})}^p \mathbb{E} \left[ \left\|  J_t^{h_N(\phi)}(\Pi_N\phi)   \right\|_{L^p(e_{t+\delta})}  \right]^p \, \mu^{(\lambda)}(d\phi)   \right\}^{1/p} \notag\\
    &\lesssim \sum_{i=1}^m \int^T_0 (1+ t^{-\zeta}) e^{(c\lambda^\gamma+2\eta - \frac{1}2-a)t}
       \| l_i  \|_{L^{\frac{p}{p-1}}(e_{t+\delta}^{-1})}\, dt \notag\\
    &\qquad \times \left\{ \int \| h_N(\phi)\|_{\mathcal{C}^{-\varepsilon}(w_\kappa)}^p
       \mathfrak{p}(\|\Pi_N \phi \|_{\mathcal{C}^{-a}(w_b)})
       \, \mu^{(\lambda)}(d\phi) \right\}^{1/p}. \label{eq:proof-1}
\end{align}
Moreover, choosing $K>0$ such that $\mbox{supp}(l_i) \subset B_K \coloneqq \{ x\in \mathbb{R}^2\ ;\ \langle x \rangle \le K   \}$ for all $1\le i\le m$, we have 
\begin{align}
 \| l_i  \|_{L^{\frac{p}{p-1}}(e_{t+\delta}^{-1})}
   &= \left\{ \int_{\mathbb{R}^2} |l_i(x)|^{\frac{p}{p-1}}
       e^{(t+\delta)^{1/2} \langle x \rangle^{\theta} \frac{p}{p-1}} \, dx \right\}^{\frac{p-1}{p}} \notag\\
   &\le \|l_i\|_{L^\infty(\mathbb{R}^2)}
       \left\{ \int_{B_K }e^{(t+\delta)^{1/2} \langle x \rangle^{\theta} \frac{p}{p-1}}\, dx \right\}^{\frac{p-1}{p}} \notag\\
 &\le \|l_i\|_{L^\infty(\mathbb{R}^2)} |B_K|^{\frac{p-1}{p}} e^{(t+\delta)^{1/2}K^{\theta}}.
\end{align}
Consequently, we obtain
\begin{align}
    &\sup_{1\le i\le m} \int^\infty_0 (1+t^{-\zeta}) e^{(c\lambda^\gamma+2\eta - \frac{1}2-a)t}
       \| l_i  \|_{L^{\frac{p}{p-1}}(e_{t+\delta}^{-1})}\, dt \notag\\
    &\quad\lesssim \int^\infty_0 (1+t^{-\zeta})
       e^{(c\lambda^\gamma+2\eta - \frac{1}2 -a )t + (t+\delta)^{1/2} K^{\theta}}
       \, dt < \infty \label{eq:proof-2}
\end{align}
when we choose $a$ such that
\begin{equation}\label{eq:cond-a}
c\lambda^\gamma +2\eta - \frac{1}2 -a <0   
\end{equation}
since we chose $\zeta$ smaller than $1$.
In particular, $a$ can be chosen independently of $K$.
Furthermore, Lemma~\ref{lem:phi42-wick} gives 
\begin{equation}
    \lim_{N\rightarrow \infty }\int \| h_N(\phi)\|_{\mathcal{C}^{-\varepsilon}(w_\kappa)}^p \mathfrak{p}(\|\Pi_N \phi \|_{\mathcal{C}^{-a}(w_b)})  \,  \mu^{(\lambda)}(d\phi) =0. \label{eq:proof-3}
\end{equation}
From \eqref{eq:proof-1}, \eqref{eq:proof-2} and \eqref{eq:proof-3}, we have
\begin{equation}
  \begin{aligned}
    &\left\| (a- \bar{\mathcal{L}_p})\tilde{G}_T^{(N)}(\cdot) - F( \cdot)  \right\|_{L^p(\mu^{(\lambda)})} \\
    &\quad\le \left\| (a- \bar{\mathcal{L}})\tilde{G}_T^{(N)}(\cdot) - F(\Pi_N \cdot) - e^{-aT}P_T^{(N)}F(\Pi_N\cdot) \right\|_{L^p(\mu^{(\lambda)})} \\
    &\qquad + \left\| e^{-aT}P_T^{(N)}F(\Pi_N\cdot)  \right\|_{L^p(\mu^{(\lambda)})}    + \left\| F(\Pi_N \cdot) - F(\cdot) \right\|_{L^p(\mu^{(\lambda)})} \to 0
  \end{aligned}
\end{equation}
as $T, N\to \infty$.
This proves that the closure of the range of $a -\bar{\mathcal{L}_p}$ contains $\mathfrak{F}C_b^\infty$.
Combining this with the density of $\mathfrak{F}C_b^\infty$ in $L^p(\mu^{(\lambda)})$ and the closedness of the range of $a -\bar{\mathcal{L}_p}$, which follows from the closedness and dissipativity of $\bar{\mathcal{L}_p}$, we obtain the desired result.
    
\end{proof}

\section{Estimates on the linearized equation}\label{sec:jacobi}
This section is devoted to the proof of Lemma~\ref{lem:linear}.
Our proof is based on the stopping time argument introduced in \cite{KT25}. Following their strategy, we consider the random operator $J(s,t): h \mapsto J_{s,t}^h$, where $\{ J_{s,t}^h\}_{t\in [s,\infty)}$ solves the restarted equation
\begin{align}
    \begin{cases}
        \partial_t J_{s,t} =  \frac{1}2 (\Delta -1)J_{s,t} - \frac{3\lambda}2 \left\{ \left(  \Phi_t^{(N)}  \right)^2 - C_N\right\} J_{s,t}      \\
        J_{s,s}=h.
    \end{cases}
\end{align}
Note that $J_{0,t}^h = J_t^h$ and $J_{t}^h = J(t_{n},t) J(t_{n-1},t_{n}) \cdots J(t_1, t_{2}) J(0,t_1)h$ for any $0 \le t_1 \le t_2 \le \cdots \le t_{n}\le t$.
We also introduce the process $\{ \Phi_{s,t}^f\}_{t \in [s,\infty)}$ defined by the equation
\begin{align}
    \begin{cases}
        \partial_t \Phi_{s,t}^{(N)}= \frac{1}2(\Delta-1)\Phi_{s,t}^{(N)}- \frac{\lambda}2  \left\{ \left(  \Phi_t^{(N)}  \right)^3 -3 C_N\Phi_t^{(N)}\right\} + \Pi_N\partial_t W_t\\
        \Phi_{s,s}^{(N)}= f. \label{eq:SQ:restart}
    \end{cases}
\end{align}
Then $\Phi_t(\phi) = \Phi_{s,t}^f\big|_{f=\Phi_s(\phi)}$ for any $t\ge s \ge 0$, where $\Phi_t(\phi)= \Phi_t^{(N)}(\phi)$ is the solution to \eqref{eq:SQ}. We decompose $\Phi_{s,t}^f$ as $\Phi_{s,t}^f = X_{s,t} + Y_{s,t}^f$, where $\{X_{s,t}\}_{t\in [s,\infty)}$ and $\{Y_{s,t}^f\}_{t\in [s,\infty)}$ solve the equations
\begin{align}
    \begin{cases}
        \partial_t X_{s,t}= \frac{1}2(\Delta-1)X_{s,t} + \Pi_N\partial_t W_t\\
        X_{s,s}= 0 \label{eq:X:restart}
    \end{cases}
\end{align}
and
\begin{align}
    \begin{cases}
        \partial_t Y_{s,t}= \frac{1}2(\Delta-1)Y_{s,t}- \frac{\lambda}2 \sum_{k=0}^3 \binom{3}k   \Wick{  X_{s,t}  ^{3-k}}_{C_N}  \left(  Y_{s,t}  \right)^k \\
        Y_{s,s}= f, \label{eq:Y:restart}
    \end{cases}
\end{align}
respectively, where we write
\begin{align}
\Wick{  X_{s,t} ^{3}}_{C_N} = X_{s,t}^3 - 3C_NX_{s,t},\ \Wick{  X_{s,t}^{2}}_{C_N}   = X_{s,t}^2 - C_N,\ \Wick{X_{s,t}^{1}}_{C_N} = X_{s,t},\ \Wick{ X_{s,t}^{0}}_{C_N}=1.
\end{align}

\subsection{Deterministic estimates} 
In this subsection, we establish estimates for the restarted processes $J_{s,t}$. The arguments are purely deterministic. The estimates in Lemma~\ref{lem:est-linear}, the main result of this subsection, will be combined with a stopping time argument in the next subsection to complete the proof of Lemma~\ref{lem:linear}. 

For $0\le s<t$, define
\begin{equation}\label{eq:def-X-norm}
\mathscr{X}_{s,t} \coloneqq \sup_{1\le k\le 3}\sup_{s\le u\le t} (u-s)^{\eta'} \| \Wick{  X_{s,u} ^{3-k}}_{C_N} \|_{\mathcal{C}^{-\varepsilon'}(w_{\kappa'})},
\end{equation}
where $\varepsilon',\eta'\kappa'>0$ are taken sufficiently small in what follows; their precise values are immaterial.

\begin{lem}\label{lem:est-linear}
Let $0\le s\le t\le 1, p\ge 1, r>1,  \delta, \kappa, \varepsilon, \zeta, \alpha, \beta>0$ and $0<\theta<1$. 
\begin{enumerate}
\item
If $0<\kappa\le 1$, $r/(r-1)< 1/(\varepsilon+\kappa)$ and $0<\alpha < 1 - \frac{r(1+\varepsilon)}{2}$, there exist polynomial $\mathfrak{p}_1$ such that 
\begin{equation}\label{eq:bound-J-L}
  \begin{aligned}
        &\| J_{s,t}^h \|_{L^{p}(e_{t -s +\delta})}^{p} \le 2^{p(r-1)/r}  \| h\|_{L^{p}(e_\delta)}^{p}\exp\left( -\frac{p}2(t-s) +  \lambda^{r/2} (t-s)^\alpha \mathfrak{p}_1(\mathscr{X}_{s,t})  \right),
      \end{aligned}
    \end{equation}
    where $\mathfrak{p}_1$ does not depend on $ t, s, N, \lambda$.  
    \item 
    Fix any $\lambda^*>0$. If $p\in \mathbb{N}, 0< \lambda \le \lambda^*, p(t+\delta)\theta<1$, $\kappa\le \theta/2 $,
    $\zeta> 3p\varepsilon$ and $0< \beta < 1 - \frac{p(1+\varepsilon)}{(2p-1)(1-\varepsilon)}-\zeta$, there exist $a,b>0$ and polynomials $\mathfrak{p}_2,\mathfrak{p}_3$  such that
\begin{align}
&\left\| J_{0,t}^h \right\|^{2p}_{L^{2p}(e_{t+\delta})} \le t^{-\zeta} \| h\|_{\mathcal{C}^{-\varepsilon}(w_\kappa)}^{2p} \mathfrak{p}_2 \left( \| \phi\|_{\mathcal{C}^{-a}(w_b)} ,\mathscr{X}_{0,t}\right) \exp\left(-pt + \lambda^{\frac{2p}{2p-1}\cdot \frac{1-2\varepsilon}{2-2\varepsilon}} t^\beta\mathfrak{p}_3(\mathscr{X}_{0,t})   \right), \label{eq:bound-J-C}
\end{align} 
where $\mathfrak{p}_2, \mathfrak{p}_3$ do not depend on $t, N, \lambda$.
\end{enumerate}
\end{lem}

\begin{rem}\label{rem:exp-integrability}
Note that the estimate \eqref{eq:bound-J-C} alone is not sufficient to prove Lemma~\ref{lem:linear} since we do not know how to control the exponential moment of $\mathfrak{p}_3(\mathscr{X}_{0,t})$ directly. 
In Section~\ref{sec:stopping-time}, the bounds in Lemma~\ref{lem:est-linear} will be used repeatedly on suitably defined random intervals $[\tau_i,\tau_{i+1}]$, $i=0,1,2,\cdots n$ such that $0=\tau_0< \tau_1< \cdots<\tau_{n-1}<t\le \tau_n$ to bypass this difficulty. More precisely, the estimate \eqref{eq:bound-J-C} will be used only on the first interval $[0,\tau_1]$ and the estimate \eqref{eq:bound-J-L} will be used repeatedly on $[\tau_i,\tau_{i+1}]$, $i=1,2, \cdots,n$; see Lemma~\ref{lem:stopping1}. 
\end{rem}

We first establish the following estimate, which will be used to prove part (i) of Lemma~\ref{lem:est-linear}. The proof is based on a technique from \cite{DHYZ25} that uses suitable space-time weights. 

\begin{lem}[Estimate on $J_{s,t}^h$ with $ h \in L^p$]\label{lem:J-Lp}
Let $s,t\ge 0, p\ge 1, r>1,  \delta,  \kappa, \varepsilon>0$. Assume that $0< t -s \le 1$, $0<\kappa\le 1$ and $r/(r-1)< 1/(\varepsilon+\kappa)$.
Then there exists $c>0$ such that 
    \begin{equation}
      \begin{aligned}
        &\| J_{s,t}^h \|_{L^{p}(e_{t -s +\delta})} \le 2^{\frac{r-1}r} e^{\frac{s-t}2} \| h\|_{L^{p}(e_\delta)} \\
        &\qquad\qquad\qquad\exp\left(c \lambda^r \int_s^t \left\| \Wick{ X_{s,u}^2 }_{C_N} + 2X_{s,u}Y_{s,u}^f \right\|_{\mathcal{C}^{-\varepsilon/2}(w_{\theta\kappa})}^r \, du\right)
          \end{aligned}
    \end{equation}
    with $f = \Phi_s^{(N)}(\phi)$, where $c$ does not depend on $t,s, N,\lambda$.
\end{lem}

\begin{proof}

Let $\tilde{J}_{s,t}^h$ be the solution to 
\begin{align}
    \begin{cases}
        \partial_t \tilde{J}_{s,t} =  \frac{1}2 (\Delta -1)\tilde{J}_{s,t} - 3\lambda \left\{  \Wick{ X_{s,t}^2 }_{C_N} + 2X_{s,t}Y_{s,t} \right\}\tilde{J}_{s,t}      \\
        \tilde{J}_{s,s}=|h|.
    \end{cases}
\end{align}
Since $(\Phi_t^{(N)})^2 - C_N= \Wick{(X_{s,t})^2}_{C_N} + 2X_{s,t}Y_{s,t} + Y_{s,t}^2$ and $\lambda Y_{s,t}^2$ is nonnegative, the comparison argument in \cite[Lemma 3.27]{DHYZ25} gives $|J^h_{s,t}| \le \tilde{J}^h_{s,t}$ pointwise. In particular, we have
\begin{equation}\label{eq:J-comparison}
    \|J^h_{s,t}\|_{L^p(e_{t-s+\delta})} \le \| \tilde{J}^h_{s,t} \|_{L^p(e_{t-s+\delta})}.
\end{equation}
Recall that we introduced the two weight functions 
\begin{equation}
    e_{u}(x) \coloneqq e^{-u^{1/2}  \langle x \rangle^{\theta}},\qquad w_\kappa(x) \coloneqq \langle x \rangle^{-\kappa}
\end{equation}
in Section~\ref{sec:pre}.
In what follows, we suppress the dependence on $C_N$ and simply write $\Wick{ X_{s,t}^2 }_{C_N} = \Wick{ X_{s,t}^2 }$.
For any $\varepsilon \ge 0$ and $\gamma\le 0$, the Schauder estimate (Lemma~\ref{lem:schauder}) gives 
\begin{align*}
    &\left\|\tilde{J}_{s,t}^h  \right\|_{B^{\varepsilon}_{p,p}(e_{t -s +\delta})} = \left\| e^{\frac{t-s}2(\Delta -1)}|h|  - \frac{3\lambda}2 \int^{t}_s e^{\frac{t-u}2(\Delta -1)} \left(\left\{  \Wick{ X_{s,u}^2 } + 2X_{s,u}Y_{s,u} \right\} \tilde{J}_{s,u}^h  \right)\, du  \right\|_{B^{\varepsilon}_{p,p}(e_{t-s+\delta})}\\
    &\lesssim (t-s)^{-\frac{ \varepsilon}2}e^{\frac{s-t}2} \| h \|_{L^{p}(e_{t-s+\delta})} + \lambda \int^t_s (t-u)^{\frac{\gamma -\varepsilon}2}e^{\frac{u-t}2} \left\|  \left\{  \Wick{ X_{s,t}^2 } + 2X_{s,u}Y_{s,u} \right\} \tilde{J}_{s,u}^h\right\|_{B^{\gamma}_{p,p}(e_{t-s+\delta})}\, du.
\end{align*}
Therefore, letting $\gamma= -\varepsilon$ and noting that 
\begin{align}
\begin{aligned}
    e_{t-s+\delta}(x)&= e_{u-s+\delta}(x)  \exp{\left( \left\{(u-s+\delta)^{1/2} - (t-s+\delta)^{1/2}  \right\}\langle x \rangle^\theta  \right) } \\
 &\le e_{u-s+\delta}(x)  \exp{  \left( \frac{1}2 (1+\delta)^{-1/2}(u-t) \langle x \rangle^\theta  \right) }  \\
 &\lesssim e_{u-s+\delta}(x) (t-u)^{-\kappa}\langle x \rangle^{-\theta \kappa}
 \end{aligned}
 \end{align}   
by $0\le t-u\le t-s\le 1$, we have from the product estimate (Lemma~\ref{lem:product-special}) that 
\begin{align*}
\left\|\tilde{J}_{s,t}^h  \right\|_{B^{\varepsilon}_{p,p}(e_{t-s+\delta})} &\lesssim (t-s)^{-\frac{\varepsilon}2}e^{\frac{s-t}2} \| h \|_{L^{p}(e_{\delta})} \\
&\  + \lambda \int^t_s (t-u)^{-(\varepsilon + \kappa)} e^{\frac{u-t}2} \left\| \Wick{ X_{s,u}^2 } + 2X_{s,u}Y_{s,u} \right\|_{\mathcal{C}^{-\varepsilon/2}(w_{\theta\kappa})} \left\| \tilde{J}_{s,u}^h\right\|_{B^{\varepsilon}_{p,p}(e_{u-s+\delta})}   \, du.
\end{align*} 
A similar argument yields 
\begin{align*}
&\left\|\tilde{J}_{s,t}^h  \right\|_{L^{p}(e_{t-s+\delta})}\le e^{\frac{s-t}2} \| h \|_{L^{p}(e_{\delta})} \\
&\qquad + C\lambda  \int^t_s (t-u)^{-(\varepsilon + \kappa)}e^{\frac{u-t}2} \left\| \Wick{X_{s,u}^2 } + 2X_{s,u}Y_{s,u} \right\|_{\mathcal{C}^{-\varepsilon/2}(w_{\theta\kappa})} \left\| \tilde{J}_{s,u}^h\right\|_{B^{\varepsilon}_{p,p}(e_{u-s+\delta})}   \, du.
\end{align*} 
Thus, setting 
\begin{equation} g(s,t) \coloneqq e^{\frac{t-s}2} \left\{\left\|\tilde{J}_{s,t}^h  \right\|_{L^{p}(e_{t-s+\delta})}    \vee c (t-s)^{\frac{\varepsilon }2} \left\|\tilde{J}_{s,t}^h  \right\|_{B^{\varepsilon}_{p,p}(e_{t-s+\delta})} \right\}    
\end{equation}
with a sufficiently small $c>0$
and letting $q,r>1$ be such that $1/q+1/r=1$, $(\varepsilon + \kappa)q <1$, and $-q\varepsilon/2>-1$, we have 
\begin{align*}
    &g(s,t)  \le \| h \|_{L^{p}(e_{\delta})} \\
    &\qquad + C (t-s)^{\frac{\varepsilon }2} e^{\frac{t-s}2}\lambda  \int^t_s (t-u)^{-(\varepsilon + \kappa)}e^{\frac{u-t}2} \left\|   \Wick{X_{s,u}^2 } + 2X_{s,u}Y_{s,u}\right\|_{\mathcal{C}^{-\varepsilon/2}(w_{\theta\kappa})} (u-s)^{-\frac{ \varepsilon}2}e^{\frac{s-u}2}g(s,u)   \, du\\
    &\qquad\le \| h \|_{L^{p}(e_{\delta})} + C\theta_1^{-\kappa} (t-s)^{\frac{\varepsilon}2} \lambda \left(  \int^t_s \left\{ (t-u)^{-(\varepsilon+\kappa)} (u-s)^{-\frac{\varepsilon}2}  \right\}^q\, du   \right)^{1/q} \\
    &\qquad \qquad\qquad\qquad \times\left(  \int^t_s \left\{ \left\|  \Wick{X_{s,u}^2 } + 2X_{s,u}Y_{s,u} \right\|_{\mathcal{C}^{-\varepsilon/2}(w_{\theta\kappa})}g(s,u)  \right\}^r\, du  \right)^{1/r}\\
    &\qquad\le \| h \|_{L^{p}(e_{\delta})} + C \lambda \left(  \int^t_s \left\{ \left\|  \Wick{ X_{s,u}^2 } + 2X_{s,u}Y_{s,u} \right\|_{\mathcal{C}^{-\varepsilon/2}(w_{\theta\kappa})}g(s,u)  \right\}^r\, du  \right)^{1/r}.
\end{align*}
Note that the implicit constant is independent of $s,t$ under the assumption $0\le t-s \le 1$.
Consequently, the elementary inequality $(|a|^r+ |b|^r)\le 2^{r-1}(|a|^r + |b|^r)$ yields
\begin{equation}
    g(s,t)^r \le 2^{r-1}\| h \|_{L^{p}(e_{\delta})}^r + C \lambda^r \int^t_s \left\| \Wick{ X_{s,u}^2 } + 2X_{s,u}Y_{s,u} \right\|_{\mathcal{C}^{-\varepsilon/2}(w_{\theta\kappa})}^r g(s,u)^r  \, du.
\end{equation}
Applying Gronwall's lemma, we obtain
\begin{equation}
    g(s,t)^r\le 2^{r-1}\| h \|_{L^{p}(e_{\delta})}^r \exp\left(c\lambda^r \int^t_s \left\| \Wick{ X_{s,u}^2 } + 2X_{s,u}Y_{s,u} \right\|_{\mathcal{C}^{-\varepsilon/2}(w_{\theta\kappa})}^r \, du\right)
\end{equation}
for some $c>0$.
Together with \eqref{eq:J-comparison}, this implies the desired result.

\end{proof}



The following result will be used to show the bound~\eqref{eq:bound-J-C}.

\begin{lem}[Estimate on $J_{0,t}^h$ with $h \in \mathcal{C}^{-\delta}$]\label{lem:J-C}
Let $t\ge 0, p\in \mathbb{N}, \delta,\kappa, \varepsilon, \zeta>0$ and $0<\theta<1$.
Fix any $\lambda^*>0$ and assume that $0\le t<1,, 0< \lambda \le \lambda^*, p(t+\delta)\theta<1$, $\kappa\le \theta/2 $.
Then there exists a constant $c>0$, independent of $t, N, \lambda$, such that  
\begin{align}
&\left\| J_{0,t}^h \right\|^{2p}_{L^{2p}(e_{t+\delta})} \lesssim t^{-\zeta} \left\{ \sup_{0< s \le t} s^\zeta D(s)  \right\} \notag\\
&\qquad \qquad\qquad \times\exp\left(-pt +c \lambda^{\frac{2p}{2p-1}} \int^t_0 s^{-\zeta}\left\{  \left\|  \Wick{X_{0,s}^2}_{C_N} + 2X_{0,s} Y_{0,s}^\phi  \right\|_{\mathcal{C}^{-\varepsilon}(w_\kappa)}^{\frac{2p}{(2p-1)(1-\varepsilon)}} +1 \right\}    \, ds \right),
\end{align}
where 
\begin{align}
    D(s) &\coloneqq \left\| I_s  \right\|^{2p}_{L^{2p}(e_{s+\delta})} + \lambda \left\|Y_{0,s}^\phi   \right\|_{L^4(e_{s+\delta}^{p/2} )}^2 \left\|I_s   \right\|_{L^{4p}(e_{s+\delta}^{1/2} )}^{2p}    \notag\\
    &\ +\lambda^{2p} \left\|  \Wick{X_{0,s}^2}_{C_N} + 2X_{0,s} Y_{0,s}^\phi  \right\|_{\mathcal{C}^{-\varepsilon}(w_\kappa)}^{2p} \left\{ \left\| I_s    \right\|_{L^{2p}(e_{s+\delta} )}^{2p} + \left\| I_s     \right\|_{L^{2p}(e_{s+\delta} )}^{2p(1-\varepsilon)} \left\|\left( \nabla  I_s \right)    \right\|_{L^{2p} ( e_{s+\delta})}^{2p\varepsilon}   \right\} \label{eq:def-D}
\end{align}
and $I_s \coloneqq e^{\frac{1}2\Delta}h$.

\end{lem}

\begin{proof}

We first define $\tilde J_{0,t} \coloneqq e^{\frac{t}2}J_{0,t}$ and consider the decomposition $\tilde J_{0,t}^h =  I_{t}+ K_{t}$, where
\begin{align}
    \begin{cases}
        \partial_t I_{t} =  \frac{1}2 \Delta I_{t}      \\
        I_{0}=h,
    \end{cases}
\end{align}
\begin{align}
    \begin{cases}
        \partial_t K_{t} =  \frac{1}2 \Delta K_{t} - \frac{3\lambda}2  \left\{ (\Phi_t^{(N)})^2 - C_N\right\}\left( K_{t} + I_{t} \right)  \label{eq:K}    \\
        K_{0}=0.
    \end{cases}
\end{align}
Note that we have
\begin{align}
     \left(\Phi_t^{(N)}\right)^2 - C_N =\left( X_{t}+Y_{t} \right)^2 - C_N = \left( X_{t}^2 - C_N \right) + 2X_{t} Y_{t} +  Y_{t}^2, \label{eq:XY}
\end{align}
where we write $X_t \coloneqq X_{0,t}$ and $Y_t \coloneqq Y_{0,t}^\phi$ for simplicity.
From \eqref{eq:K}, we have
\begin{align}
    \frac{1}{2p} \frac{d}{dt} \left\| K_t \right\|^{2p}_{L^{2p}(e_{t+\delta})}&= \int_{\mathbb{R}^2} (\partial_t K_t)K_t^{2p-1} e_{t+\delta}^{2p}\, dx - \frac{1}{2} (t+\delta)^{-1/2}\int_{\mathbb{R}^2} K_t^{2p} \langle x \rangle^{\theta} e_{t+\delta}^{2p}\, dx \notag \\
    &= \int_{\mathbb{R}^2} \left\{ \frac{1}2 \Delta K_{t} - \frac{3\lambda}2 \left\{ (\Phi_t^{(N)})^2 - C_N\right\}  \left( K_{t} + I_{t} \right)  \right\}K_t^{2p-1} e_{t+\delta}^{2p}\, dx\notag\\
    &\qquad\qquad\qquad\qquad-  \frac{1}{2} (t+\delta)^{-1/2} \int_{\mathbb{R}^2} (K_t)^{2p} \langle x \rangle^{\theta} e_{t+\delta}^{2p}\, dx. \label{eq:K-L2}
\end{align}
For the first term, integration by parts gives
\begin{align}
    \int_{\mathbb{R}^2} (\Delta K_{t})K_t^{2p-1} e_{t+\delta}^{2p}\, dx  &= - (2p-1) \int_{\mathbb{R}^2} |\nabla K_t|^2 K_t^{2p-2} e_{t+\delta}^{2p}\, dx - \int_{\mathbb{R}^2} \left(\nabla K_t \cdot \nabla e_{t+\delta}^{2p} \right) K_t^{2p-1}\, dx. \label{eq:KIbP}
\end{align}
Therefore, we have from \eqref{eq:XY}, \eqref{eq:K-L2} and \eqref{eq:KIbP} that
\begin{align}
    &\frac{1}{2p} \frac{d}{dt}\left\| K_t \right\|^{2p}_{L^{2p}(e_{t+\delta})} + \frac{2p-1}2 \int_{\mathbb{R}^2} |\nabla K_t|^2 K_t^{2p-2} e_{t+\delta}^{2p}\, dx \notag\\
    &\qquad + \frac{1}2(t+\delta)^{-1/2} \int_{\mathbb{R}^2}K_t^{2p} w_{-\theta} e_{t+\delta}^{2p}\, dx  + \frac{3\lambda}2 \int_{\mathbb{R}^2} Y_t^2 K_t^{2p} e_{t+\delta}^{2p}\, dx \notag\\
    &\eqqcolon  \frac{1}{2p}  \frac{d}{dt} \left\| K_t \right\|^{2p}_{L^{2p}(e_{t+\delta})} + A_1 + A_2 + A_3 \notag\\
    &\le \frac{1}2\left| \int_{\mathbb{R}^2} \left(\nabla K_t \cdot \nabla e_{t+\delta}^{2p} \right) K_t^{2p-1}   \, dx  \right| + \frac{3\lambda}2 \left| \int_{\mathbb{R}^2} \left\{ \Wick{X_t^2} + 2X_t Y_t \right\} K_t^{2p} e_{t+\delta}^{2p}    \, dx  \right| \notag \\
    &\qquad +3\lambda \left| \int_{\mathbb{R}^2} \left\{ \Wick{X_t^2} + 2X_t Y_t \right\} K_t^{2p-1} I_t e_{t+\delta}^{2p}    \, dx  \right| +3\lambda \left| \int_{\mathbb{R}^2} Y_t^2 K_t^{2p-1} I_t e_{t+\delta}^{2p}    \, dx  \right|\notag \\
    & \eqqcolon B_1 + B_2 + B_3 + B_4. \label{eq:est-K}
\end{align}
We estimate each term on the right-hand side separately. First, using
\begin{equation}
    \nabla e_{t+\delta}^{2p} = -2p(t+\delta)^{1/2} \theta \langle x \rangle^{\theta -2} e_{t+\delta}^{2p}(x)x
\end{equation}
and Young's inequality, we obtain
\begin{align}
    B_1 &= \left| \int_{\mathbb{R}^2}  p(t+\delta)^{1/2} \theta \left(\nabla K_t \cdot x \right) \langle x \rangle^{\theta -2} K_t^{2p-1} e_{t+\delta}^{2p}(x)  \, dx  \right|\notag \\
    &\le \int_{\mathbb{R}^2} p(t+\delta)^{1/2}\theta |\nabla K_t| |K_t|^{2p-1} \langle x \rangle^{\theta -1} e_{t+\delta}^{2p}(x) \, dx\notag \\
    &\le \frac{p(t+\delta)^{1/2}\theta}2 \int_{\mathbb{R}^2} K_t^{2p} \langle x \rangle^{2\theta -2}e_{t+\delta}^{2p}(x) \, dx + \frac{p(t+\delta)^{1/2}\theta}2 \int_{\mathbb{R}^2} |\nabla K_t|^2 K_t^{2p-2}e_{t+\delta}^{2p}(x) \, dx. \label{eq:est-B1}
\end{align}
Therefore, if $p(t+\delta)\theta <1$, this term can be absorbed into the nonnegative terms on the left-hand side.
For the second term, Lemma~\ref{lem:weight-equiv}, duality (Lemma~\ref{lem:duality}), interpolation (Lemma~\ref{lem:interpolation}), the embedding $W^{\varepsilon,1}(\mathbb{R}^2) \hookrightarrow B_{1,1}^{\varepsilon/2}(\mathbb{R}^2)$, and H\"older's inequality give 
\begin{align}
    B_2 &= \frac{3\lambda}2 \left| \int_{\mathbb{R}^2} \left\{ \Wick{X_t^2} + 2X_t Y_t \right\}w_\kappa w_{-\kappa} K_t^{2p} e_{t+\delta}^{2p}    \, dx  \right|\notag \\
    &\lesssim  \lambda \left\| \left\{ \Wick{X_t^2} + 2X_t Y_t \right\}w_\kappa \right\|_{\mathcal{C}^{-\varepsilon/2}} \left\| K_t^{2p} w_{-\kappa} e_{t+\delta}^{2p}  \right\|_{B^{\varepsilon/2}_{1,1}}\notag \\
    &\simeq   \lambda \left\|  \Wick{X_t^2} + 2X_t Y_t  \right\|_{\mathcal{C}^{-\varepsilon/2}(w_\kappa)} \left\| K_t^{2p} w_{-\kappa} e_{t+\delta}^{2p}  \right\|_{L^{1}}^{1-\varepsilon} \left\| K_t^{2p} w_{-\kappa} e_{t+\delta}^{2p}  \right\|_{W^{1,1}}^\varepsilon \notag \\
    &\lesssim \lambda \left\|  \Wick{X_t^2} + 2X_t Y_t  \right\|_{\mathcal{C}^{-\varepsilon/2}(w_\kappa)} \left\| K_t w_{-\kappa} e_{t+\delta}  \right\|_{L^{2p}}^{1-\varepsilon} \left\| K_t^{2p-1} e_{t+\delta}^{2p-1}  \right\|_{L^{2p/(2p-1)}}^{1-\varepsilon} \notag\\
    &\qquad \times \left\{    \left\| K_t^{2p} w_{-\kappa} e_{t+\delta}^{2p}  \right\|_{L^{1}}^\varepsilon +  \left\| \nabla \left( K_t^{2p} w_{-\kappa} e_{t+\delta}^{2p} \right)  \right\|_{L^{1}}^\varepsilon  \right\} \notag\\
    &\lesssim \lambda \left\|  \Wick{X_t^2} + 2X_t Y_t  \right\|_{\mathcal{C}^{-\varepsilon/2}(w_\kappa)} A_2^{\frac{1-\varepsilon}{2p}} \| K_t \|_{L^{2p}(e_{t+\delta})}^{(2p-1)(1-\varepsilon)} \times \left\{    A_2^\varepsilon +  \left\| \nabla \left( K_t^{2p} w_{-\kappa} e_{t+\delta}^{2p} \right)  \right\|_{L^{1}}^\varepsilon  \right\} \label{eq:est-B2}
\end{align}
provided that $\kappa\le \theta/2$.
Moreover, we have from H\"older's inequality that
\begin{align}
 \left\| \nabla \left( K_t^{2p} w_{-\kappa} e_{t+\delta}^{2p} \right)  \right\|_{L^{1}}&\lesssim \left\| (\nabla K_t ) K_t^{2p-1} w_{-\kappa} e_{t+\delta}^{2p}   \right\|_{L^{1}} + \left\|   K_t^{2p} \nabla \left( w_{-\kappa} e_{t+\delta}^{2p} \right)  \right\|_{L^{1}} \notag\\
 &\lesssim \left\| (\nabla K_t ) K_t^{p-1} e_{t+\delta}^{p}   \right\|_{L^{2}} \left\|  K_t^p w_{-\kappa} e_{t+\delta}^p   \right\|_{L^{2}} \notag\\
 &\qquad + \left\|   K_t  e_{t+\delta}   \right\|_{L^{2p}} \left\|   K_t^{2p-1} e_{t+\delta}^{-1} \nabla \left( w_{-\kappa} e_{t+\delta}^{2p} \right)  \right\|_{L^{2p/(2p-1)}}\notag\\
 &\lesssim  \left\| (\nabla K_t ) K_t^{p-1}  e_{t+\delta}^p   \right\|_{L^{2}}\left\|  K_t^p w_{-\theta/2} e_{t+\delta}^p   \right\|_{L^{2}} \notag\\
 &\qquad + \left\|   K_t  e_{t+\delta}   \right\|_{L^{2p}} \left\|  K_t^p w_{-\theta/2} e_{t+\delta}^p   \right\|_{L^{2}}^{(2p-1)/p}  \notag \\
 &\simeq A_1^{1/2} A_2^{1/2} + \| K_t \|_{L^{2p}(e_{t+\delta})} A_2^{\frac{2p-1}{2p}}\label{eq:est-B2-2} 
\end{align}
provided that $\kappa \le \theta/2 $. 
Therefore, from \eqref{eq:est-B2}, \eqref{eq:est-B2-2} and Young's inequality, we obtain
\begin{align}
B_2 &\lesssim  \lambda \left\|  \Wick{X_t^2} + 2X_t Y_t  \right\|_{\mathcal{C}^{-\varepsilon/2}(w_\kappa)} \notag\\
&\qquad\times\left\{ A_2^{\frac{1-\varepsilon}{2p}+ \varepsilon } \| K_t \|_{L^{2p}(e_{t+\delta})}^{(2p-1)(1-\varepsilon)} +  A_2^{\frac{1-\varepsilon}{2p}+ \frac{\varepsilon}2 } A_1^{\varepsilon/2} \| K_t \|_{L^{2p}(e_{t+\delta})}^{(2p-1)(1-\varepsilon)}  +  A_2^{\frac{1-\varepsilon}{2p}+ \frac{2p-1}{2p}\varepsilon } \| K_t \|_{L^{2p}(e_{t+\delta})}^{(2p-1)(1-\varepsilon)+ \varepsilon}     \right\}        \notag\\
&\le \eta \left( A_1+A_2 \right)  + C_{\varepsilon, \eta} \left\| K_t   \right\|_{L^{2p}(e_{t+\delta} )}^{2p} \left\{  \lambda^{\frac{2p}{(2p-1)(1-\varepsilon)}} \left\|  \Wick{X_t^2} + 2X_t Y_t  \right\|_{\mathcal{C}^{-\varepsilon/2}(w_\kappa)}^{\frac{2p}{(2p-1)(1-\varepsilon)}} \right. \notag\\
&\qquad\qquad\qquad\qquad\qquad\qquad\qquad \qquad \left.+ \lambda^{\frac{2p}{(2p-1)(1-\varepsilon) + \varepsilon}} \left\|  \Wick{X_t^2} + 2X_t Y_t  \right\|_{\mathcal{C}^{-\varepsilon/2}(w_\kappa)}^{\frac{2p}{(2p-1)(1-\varepsilon) + \varepsilon}}  \right\}\notag\\
&\lesssim \eta \left(   A_1+A_2\right)  + C_{\varepsilon, \eta, \lambda_*} \left\| K_t   \right\|_{L^{2p}( e_{t+\delta} )}^{2p}\lambda^{\frac{2p}{2p-1}} \left\{  \left\|  \Wick{X_t^2} + 2X_t Y_t  \right\|_{\mathcal{C}^{-\varepsilon/2}(w_\kappa)}^{\frac{2p}{(2p-1)(1-\varepsilon)}} +1 \right\}  , \label{eq:est-B2-3}
\end{align}
for any $\eta>0$ and $0 <\lambda \le \lambda_*$. Taking $\eta>0$ sufficiently small allows us to absorb the first term into $A_1+ A_2$.
Arguing as in \eqref{eq:est-B2} and \eqref{eq:est-B2-2}, we obtain
\begin{align*}
    B_3 &\lesssim \lambda \left\|  \Wick{X_t^2} + 2X_t Y_t  \right\|_{\mathcal{C}^{-\varepsilon}(w_\kappa)} \left\| I_t  e_{t+\delta}  \right\|_{L^{2p}}^{1-\varepsilon} \left\| K_t^{2p-1}w_{-\kappa} e_{t+\delta}^{p-1}  \right\|_{L^{\frac{2p}{2p-1}}}^{1-\varepsilon} \\
    &\qquad \times\left\{ \left\| K_t^{2p-1} I_t w_{-\kappa} e_{t+\delta}^{2p}  \right\|_{L^{1}}^\varepsilon +  \left\| \nabla \left( K_t^{2p-1} I_t w_{-\kappa} e_{t+\delta}^{2p} \right)  \right\|_{L^{1}}^\varepsilon  \right\}
\end{align*}
and
\begin{align*}
 \left\| \nabla \left( K_t^{2p-1} I_t w_{-\kappa} e_{t+\delta}^{2p} \right)  \right\|_{L^{1}}&\lesssim \left\|\left( \nabla  I_t \right) K_t^{2p-1} w_{-\kappa} e_{t+\delta}^{2p}   \right\|_{L^{1}}+ \left\|  I_t  K_t^{2p-1} \nabla \left(w_{-\kappa} e_{t+\delta}^{2p}\right)   \right\|_{L^{1}} \\
 &\qquad + \left\|\left( \nabla  K_t \right) K_t^{2p-2} I_t w_{-\kappa} e_{t+\delta}^{2p}   \right\|_{L^{1}}  \\
 &\lesssim \left\|\left( \nabla  I_t \right)  e_{t+\delta}   \right\|_{L^{2p}}\left\| K_t^{2p-1} w_{-\kappa} e_{t+\delta}^{p-1}   \right\|_{L^{\frac{2p}{2p-1}}} \\
 &\qquad +\left\|   K_t^{2p-1}  e_{t+\delta} \nabla \left(w_{-\kappa} e_{t+\delta}^{2p}\right)  \right\|_{L^{\frac{2p}{2p-1}}}\left\|  I_t   e_{t+\delta}   \right\|_{L^{2p}}  \\
 &\qquad + \left\|\left( \nabla  K_t \right) K_t^{p-1}  e_{t+\delta}^p   \right\|_{L^{2}} \left\| K_t^{p-1}w_{-\kappa} e_{t+\delta}^{p-1}   \right\|_{L^{\frac{2p}{p-1}}}  \left\| I_t e_{t+\delta}   \right\|_{L^{2p}}  .
\end{align*}
Therefore, Young's inequality gives
\begin{align}
    B_3 &\lesssim \lambda \left\|  \Wick{X_t^2} + 2X_t Y_t  \right\|_{\mathcal{C}^{-\varepsilon/2}(w_\kappa)} \left\{ \left\| I_t  e_{t+\delta}   \right\|_{L^{2p}} \left\| K_t^{p-1}w_{-\kappa} e_{t+\delta}^{p-1}   \right\|_{L^{\frac{2p}{p-1}}} \right. \notag\\
    &\left.+ \left\| I_t  e_{t+\delta}   \right\|_{L^{2p}}^{1-\varepsilon} \left\| (\nabla I_t )  e_{t+\delta}   \right\|_{L^{2p}}^\varepsilon \left\| K_t^{p-1}w_{-\kappa} e_{t+\delta}^{p-1}   \right\|_{L^{\frac{2p}{p-1}}} \right.\notag\\
    &\qquad\left. + \left\| I_t  e_{t+\delta}   \right\|_{L^{2p}} \left\| K_t^{p}w_{-\theta/2} e_{t+\delta}^p   \right\|_{L^{2}}^{\frac{2p-1}p} \right.  \notag \\
    &\left. \qquad +  \left\| I_t  e_{t+\delta}   \right\|_{L^{2p}} \left\| K_t^{p}w_{-\theta/2} e_{t+\delta}^p   \right\|_{L^{2}}^{\frac{2p-1}p (1-\varepsilon)+ \frac{p-1}p \varepsilon} \left\|\left( \nabla  K_t \right) K_t^{p-1}  e_{t+\delta}^p   \right\|_{L^{2}}^\varepsilon       \right\} \notag  \\
    &\le \eta \left( A_1+ A_2  \right) + C_{\varepsilon,\eta} \lambda^{2p} \left\|  \Wick{X_t^2} + 2X_t Y_t  \right\|_{\mathcal{C}^{-\frac{\varepsilon}2}(w_\kappa)}^{2p} \left\{ \left\| I_t  e_{t+\delta}   \right\|_{L^{2p}}^{2p} \right.\notag\\
    &\qquad \qquad \qquad \qquad \qquad \qquad \qquad \left. + \left\| I_t  e_{t+\delta}   \right\|_{L^{2p}}^{2p(1-\varepsilon)} \left\|\left( \nabla  I_t \right)  e_{t+\delta}   \right\|_{L^{2p} }^{2p\varepsilon}   \right\}. \label{eq:est-B3}
\end{align}
For the fourth term, Young's inequality gives
\begin{align}
 B_4 &\le \frac{3(2p-1)}{2p} \lambda \int_{\mathbb{R}^2} Y_t^2 K_t^{2p} e_{t+\delta}^{2p}\, dx + \frac3{2p} \lambda \int_{\mathbb{R}^2} Y_t^2 I_t^{2p} e_{t+\delta}^{2p}\, dx \notag\\
 &\le \frac{2p-1}{2p} A_3 + \frac{3}{2p} \lambda \left\|Y_t^2 e_{t+\delta}^p  \right\|_{L^2} \left\|I_t^{2p} e_{t+\delta}^p  \right\|_{L^2}. \label{eq:est-B4}
\end{align}
Combining \eqref{eq:est-K}, \eqref{eq:est-B1}, \eqref{eq:est-B2-3}, \eqref{eq:est-B3}, and \eqref{eq:est-B4}, we obtain
\begin{align}
 \left\| K_t \right\|^{2p}_{L^{2p}(e_{t+\delta})} &\lesssim  C(t) + \int^t_0 \lambda^{\frac{2p}{2p-1}} \left\{  \left\|  \Wick{X_s^2} + 2X_s Y_s  \right\|_{\mathcal{C}^{-\frac{\varepsilon}2}(w_\kappa)}^{\frac{2p}{(2p-1)(1-\varepsilon)}} +1 \right\}\left\| K_s   \right\|_{L^{2p}( e_{s+\delta} )}^{2p}  \, ds,
\end{align}
where
\begin{align}
\begin{aligned}
    C(t) &\coloneqq \lambda \left\|Y_t^2 e_{t+\delta}^p  \right\|_{L^2} \left\|I_t^{2p} e_{t+\delta}^p  \right\|_{L^2}    \notag\\
    &\ +\lambda^{2p} \left\|  \Wick{X_t^2} + 2X_t Y_t  \right\|_{\mathcal{C}^{-\varepsilon/2}(w_\kappa)}^{2p} \left\{ \left\| I_t  e_{t+\delta}   \right\|_{L^{2p}}^{2p} + \left\| I_t  e_{t+\delta}   \right\|_{L^{2p}}^{2p(1-\varepsilon)} \left\|\left( \nabla  I_t \right)  e_{t+\delta}   \right\|_{L^{2p} }^{2p\varepsilon}   \right\}   .
    \end{aligned}
\end{align}
Thus, setting $G(t) \coloneqq t^\zeta \left\| K_t \right\|^{2p}_{L^{2p}(e_{t+\delta})}$, we get
\begin{align}
    G(t) \lesssim t^\zeta C(t) + \int^t_0 \lambda^{\frac{2p}{2p-1}} \left\{  \left\|  \Wick{X_s^2} + 2X_s Y_s  \right\|_{\mathcal{C}^{-\varepsilon/2}(w_\kappa)}^{\frac{2p}{(2p-1)(1-\varepsilon)}} +1 \right\}s^{-\zeta}G(s)  \, ds.
\end{align}
Gronwall's inequality then gives
\begin{align}
    G(t) \lesssim \left\{ \sup_{0<s \le t} s^\zeta C(s) \right\}\exp\left(c \lambda^{\frac{2p}{2p-1}} \int^t_0 s^{-\zeta} \left\{  \left\|  \Wick{X_s^2} + 2X_s Y_s  \right\|_{\mathcal{C}^{-\varepsilon/2}(w_\kappa)}^{\frac{2p}{(2p-1)(1-\varepsilon)}} +1 \right\}    \, ds \right)
\end{align}
which is equivalent to
\begin{align}
\left\| K_t \right\|^{2p}_{L^{2p}(e_{t+\delta})} \lesssim t^{-\zeta}\left\{ \sup_{0<s \le t} s^\zeta C(s) \right\} \exp\left(c \lambda^{\frac{2p}{2p-1}} \int^t_0 s^{-\zeta} \left\{  \left\|  \Wick{X_s^2} + 2X_s Y_s  \right\|_{\mathcal{C}^{-\varepsilon/2}(w_\kappa)}^{\frac{2p}{(2p-1)(1-\varepsilon)}} +1 \right\}    \, ds \right) .   
\end{align}
Consequently, we obtain
\begin{align}
&\left\| J_{0,t}^h \right\|^{2p}_{L^{2p}(e_{t+\delta})} = e^{-pt} \left\| \tilde J_{0,t}^h\right\|^{2p}_{L^{2p}(e_{t+\delta})} \lesssim e^{-pt} \left(\left\| I_{t} \right\|^{2p}_{L^{2p}(e_{t+\delta})} + \left\| K_{t} \right\|^{2p}_{L^{2p}(e_{t+\delta})}  \right) \notag\\
&\lesssim t^{-\zeta}\left\{ \sup_{0<s \le t} s^\zeta \left(\| I_s\|^{2p}_{L^{2p}(e_{s+\delta})} +C(s)\right) \right\} \exp\left(c \lambda^{\frac{2p}{2p-1}} \int^t_0 s^{-\zeta} \left\{  \left\|  \Wick{X_s^2} + 2X_s Y_s  \right\|_{\mathcal{C}^{-\varepsilon/2}(w_\kappa)}^{\frac{2p}{(2p-1)(1-\varepsilon)}} +1 \right\}    \, ds \right). 
\end{align}

\end{proof}

We also prepare the following bounds for $Y_{s,t}^f$, which follows from the results given in Section~\ref{sec:CDFI}. Recall that $\mathscr{X}_{s,t}$ is defined by \eqref{eq:def-X-norm}.

\begin{lem}\label{lem:apriori}
 Fix $\lambda^*>0$. Let $Y_{s,t}^f$ be the solution to the equation~\eqref{eq:Y:restart} and let $0\le s< t\le 1$. Fix any $0<\eta<1$, $\varepsilon>0$, $\kappa>0$. Then there exist some polynomials $\mathfrak{p}_1, \mathfrak{p}_2$ and small constants $a,b>0$ such that 
    \begin{align}
      \lambda^{1/2}  \| Y_{s,t}^f \|_{\mathcal{C}^\eta(w_\kappa)} \lesssim (t-s)^{-(1+\eta+\varepsilon )/2} \mathfrak{p}_1(\mathscr{X}_{s,t}), \label{eq:est-YY2}
    \end{align}
    \begin{align}
     \| Y_{0,t}^\phi \|_{\mathcal{C}^\eta(w_\kappa)} \lesssim t^{-(\eta+a+\varepsilon )/2} \mathfrak{p}_2(\| \phi \|_{\mathcal{C}^{-a}(w_b)}, \mathscr{X}_{0,t}) ,  \label{eq:est-YY3}
    \end{align}
    where the first estimate is uniform in the initial condition $f$ and $0< \lambda\le \lambda^*$.
\end{lem}

\begin{proof}
The bound \eqref{eq:est-YY2} follows from \eqref{est:YY2}, while \eqref{eq:est-YY3} follows from \eqref{est:Y3}. 
\end{proof}

Combining Lemmas~\ref{lem:J-Lp}, \ref{lem:J-C}, and~\ref{lem:apriori}, we can now finish the proof of Lemma~\ref{lem:est-linear}.

\begin{proof}[Proof of Lemma~\ref{lem:est-linear}]
We first prove part (i). From the product estimate (Lemma~\ref{lem:product}) and the coming-down-from-infinity estimate \eqref{eq:est-YY2}, we have for any $\tilde{\varepsilon}, \tilde{\kappa}>0$ that
\begin{align}
\left\| \Wick{ X_{s,t}^2 } + 2X_{s,t}Y_{s,t}^f \right\|_{\mathcal{C}^{-\tilde\varepsilon}(w_{\tilde\kappa})}  &\lesssim \left\| \Wick{ X_{s,t}^2 }  \right\|_{\mathcal{C}^{-\tilde\varepsilon/2}(w_{\tilde\kappa})}  + \left\| X_{s,t} \right\|_{\mathcal{C}^{-\tilde\varepsilon}(w_{\tilde\kappa/2})}\left\| Y_{s,t} \right\|_{\mathcal{C}^{\tilde\varepsilon}(w_{\tilde\kappa/2})} \notag\\
&\lesssim \lambda^{-1/2} (t-s)^{-(1 +\tilde\varepsilon)/2} \mathfrak{p} ( \mathscr{X}_{s,t} ). \label{eq:5-1}
\end{align}
Taking $\tilde\varepsilon=\varepsilon/2$ and applying Lemma~\ref{lem:J-Lp}, we get (i). We now turn to part (ii). In view of Lemma~\ref{lem:J-C} and \eqref{eq:5-1}, it suffices to estimate each term in \eqref{eq:def-D}. From Lemma~\ref{lem:schauder} and $e_\delta \lesssim w_\kappa$, we have
\begin{align}
    \| I_s\|_{L^{2p}(e_{s+\delta})}^{2p} \lesssim s^{-2p\varepsilon} \| h\|^{2p}_{\mathcal{C}^{-\varepsilon}(w_\kappa)},
\end{align}
\begin{align}
 \| I_s\|_{L^{4p}(e_{s+\delta}^{1/2})}^{2p} \lesssim s^{-2p\varepsilon} \| h\|^{2p}_{\mathcal{C}^{-\varepsilon}(w_\kappa)},   
\end{align}
\begin{align}
    \| \nabla I_s \|_{L^{2p}(e_{s+\delta})}^{2p\varepsilon} \lesssim s^{-2p\varepsilon(1+\varepsilon)} \| h\|^{2p}_{\mathcal{C}^{-\varepsilon}(w_\kappa)}.
\end{align}
Moreover, from the estimate \eqref{eq:est-YY3} and Lemma~\ref{lem:product}, we have for sufficiently small $a,b>0$ that
\begin{align}
    \left\| Y_{0,s}^\phi \right\|^2_{L^4 \left( e_{s+\delta}^{p/2} \right)} \lesssim \left\| Y_{0,s}^\phi \right\|^2_{L^\infty \left(w_\kappa \right)} \lesssim s^{-2a} \mathfrak{p}( \| \phi\|_{\mathcal{C}^{-a}(w_b)}  , \mathscr{X}_{0,s}),
\end{align}
\begin{align}
\left\| \Wick{ X_{0,s}^2 } + 2X_{0,s}Y_{0,s}^\phi \right\|_{\mathcal{C}^{-\varepsilon}(w_{\kappa})}^{2p} &\lesssim  \left\| \Wick{ X_{0,s}^2 }  \right\|_{\mathcal{C}^{-\varepsilon}(w_{\kappa})}^{2p}  + \left\| X_{0,s} \right\|^{2p}_{\mathcal{C}^{-\varepsilon/2}(w_{\kappa/2})}\left\| Y_{0,s} \right\|^{2p}_{\mathcal{C}^{\varepsilon}(w_{\kappa/2})}  \notag\\
&\lesssim s^{-p(\varepsilon + 2a )} \mathfrak{p}( \| \phi\|_{\mathcal{C}^{-a}(w_b)}  , \mathscr{X}_{0,s}).
\end{align}
The desired result now follows by choosing $a>0$ so small that $a<\varepsilon$.

\end{proof}

\subsection{Proof of Lemma~\ref{lem:linear}}\label{sec:stopping-time}
We adapt the argument in \cite{KT25} to bypass the problem of exponential integrability mentioned in Remark~\ref{rem:exp-integrability}. Let $\varepsilon',\eta',\kappa'>0$ and define the sequence $\{ \tau_n \}_{n\in \mathbb{N}_0}$ of stopping times by 
\begin{equation}
  \begin{aligned}
    &\tilde{\tau}_0 = \tau_0 =0, \\
    &\tilde{\tau}_n \coloneqq \inf\left\{t \ge \tau_{n-1}\ ;\ \sup_{k=1,2,3} \sup_{\tau_{n-1} \le u\le t} (u-\tau_{n-1} )^{\eta'} \left\| \Wick{X_{\tau_{n-1},t}^k}_{C_N}  \right\|_{\mathcal{C}^{-\varepsilon'}(w_{\kappa'})} \ge K   \right\}, \\
   & \tau_n \coloneqq  \tau_{n-1} +  (\tilde{\tau}_n - \tau_{n-1})\wedge \nu  \qquad \mbox{for}\ n\ge 1,
  \end{aligned}
\end{equation}
where $0<\nu<1$ and $K>0$ are chosen so that
\begin{equation}
    \mathbb{P}( \tilde{\tau}_1  \le \nu) \le \frac{1}4 .
\end{equation}
We define $N(t) \coloneqq \inf \{ n\ge 1\ ;\ \tau_n \ge t\}$. Applying the estimates in Lemma~\ref{lem:est-linear} repeatedly on the intervals $(\tau_i,\tau_{i-1} ]$, we get the following result.

\begin{lem}\label{lem:stopping1}
Fix $\lambda^*>0$.
Let $p \in 2\mathbb{N}, r>1, \zeta,\varepsilon,\kappa,\alpha>0$, $0<\lambda<\lambda*$, and $0<\theta<1$ be such that $0<\kappa\le \theta/2$, $\frac{r}{r-1}< \frac{1}{\varepsilon+\kappa}$, $0< \alpha < 1 - r(1+\varepsilon)/2$ and $3p\varepsilon < \zeta< 1- \frac{p(1+\varepsilon)}{(2p-1)(1-\varepsilon)}$. Moreover, assume that $p(\nu+ \delta)\theta< 1$.
  Then there exist $\varepsilon',\eta',\kappa'>0$ and $a,b$ such that, for $N(t)$ defined as above, we have
   \begin{equation}
     \begin{aligned}
        &\left\| J^h_t \right\|_{L^{p}(e_{t+\delta})}^{p}  \le (\tau_1 \wedge t)^{-\zeta} \| h\|_{\mathcal{C}^{-\varepsilon}(w_\kappa)}^{p} \mathfrak{p}_1 (\| \phi \|_{\mathcal{C}^{-a}(w_b)}) \\
        &\qquad \qquad \qquad \times\exp\left(  \frac{p(r-1)}{r} N(t) \log 2 + \lambda^{r/2}\mathfrak{p}_2(K)\nu^\alpha N(t)  -   \frac{p}2  t \right)
        \end{aligned}
   \end{equation} 
   for polynomials $\mathfrak{p}_1, \mathfrak{p}_2$ independent of $t, \lambda, N$.
\end{lem}

\begin{proof}
We first consider the case $\tau_1\le t$.
Writing $J_{s,t}^h = J(s,t)[h]$ and applying \eqref{eq:bound-J-L}, we have for $\tau_{n-1}\le t \le \tau_n$ with $n\ge 2$ that
\begin{align}
   & \left\| J^h_t \right\|_{L^{p}(e_{t+\delta})}^{p} = \left\| J(\tau_{n-1},t)[J_{\tau_{n-1}}^h] \right\|_{L^{p}(e_{t+\delta})}^{p}\notag\\
    &\le \exp\left( \frac{p(r-1)}r \log 2 - \frac{p}2(t-\tau_{n-1}) +  \lambda^{r/2} \mathfrak{p}_1(K)(t-\tau_{n-1})^\alpha   \right) \left\| J^h_{\tau_{n-1}} \right\|_{L^{p}(e_{\tau_{n-1}+\delta})}^{p} \notag\\
    &\le \exp\left( \frac{p(r-1)}r \log 2 - \frac{p}2(t-\tau_{n-1}) + \lambda^{r/2} \mathfrak{p}_1(K) \nu^\alpha   \right) \left\| J^h_{\tau_{n-1}} \right\|_{L^{p}(e_{\tau_{n-1}+\delta})}^{p}. 
    \end{align}
Iterating this bound, we obtain
\begin{align}
\left\| J^h_t \right\|_{L^{p}(e_{t+\delta})}^{p}    &\le  \exp\left( \frac{p(r-1)}r (n-1) \log 2  + \lambda^{r/2} \mathfrak{p}_1(K) \nu^\alpha (n-1) - \frac{p}2(t-\tau_{1})  \right) \left\| J^h_{\tau_{1}} \right\|_{L^{p}(e_{\tau_{1}+\delta})}^{p}.   \label{eq:J-st}
\end{align}
Thus, we have from \eqref{eq:bound-J-C} that 
\begin{align}
    \left\| J^h_{t} \right\|_{L^{p}(e_{t+\delta})}^{p}   &\lesssim \tau_1^{-\zeta} \| h\|_{\mathcal{C}^{-\varepsilon}(w_\kappa)}^{p} \mathfrak{p}_2 (\| \phi \|_{\mathcal{C}^{-a}(w_b)}, K) \notag\\
    &\ \times\exp\left(  \frac{p(r-1)}{r} N(t) \log 2 + \lambda^{r/2}\mathfrak{p}_1(K)\nu^\alpha N(t)  -   \frac{p}2  t  + \mathfrak{p}_3(K)\tau_1 \right) ,
\end{align}
which implies the desired result since $\tau_1\le \nu< 1$.
The case $\tau_1>t$ follows similarly from Lemma~\ref{lem:J-C}.

\end{proof}

The following lemma is a crucial step in our argument. The difference from \cite[Proposition~3.3]{KT25} is that we have to deal with the inverse moment of $\tau_1$.

\begin{lem}\label{lem:stopping2}
    For any $\gamma_1, \gamma_2$ such that $3\gamma_2 \le \log 2$, we have
    \begin{equation}
        \mathbb{E} \left[ (\tau_1 \wedge t )^{-\gamma_1} e^{\gamma_2N(t)} \right] \lesssim (t^{-\gamma_1}+1 )e^{\frac{2\log2}{\nu} t}.
    \end{equation}   
\end{lem}

\begin{proof}
By the Cauchy--Schwarz inequality, we have
\begin{align}
 \mathbb{E} \left[ (\tau_1 \wedge t)^{-\gamma_1} e^{\gamma_2N(t)} \right]  &\le \mathbb{E} \left[ (\tau_1 \wedge t )^{-2\gamma_1}  \right]^{1/2} \mathbb{E} \left[ e^{2\gamma_2N(t)} \right]^{1/2} \notag\\
 &\le \left\{ t^{-\gamma_1 } + \mathbb{E} \left[ \tau_1^{-2\gamma_1}  \right]^{1/2} \right\} \mathbb{E} \left[ e^{2\gamma_2N(t)} \right]^{1/2} .
\end{align}
    Markov's inequality and Corollary~\ref{cor:wick} give
    \begin{align}
    \mathbb{E} \left[ \tau_1^{-2\gamma_1} \right] &= \mathbb{E} \left[ \tau_1^{-2\gamma_1}\ ;\ \tau_1\ge 1 \right] + \sum_{\ell=1}^\infty \mathbb{E} \left[ \tau_1^{-2\gamma_1}\ ;\ 2^{-\ell} \le \tau_1 < 2^{-\ell+1} \right]\notag \\
    &\le 1 + \sum_{\ell=1}^\infty 2^{2\ell \gamma_1} \mathbb{P}( \tau_1 < 2^{-\ell+1}) = 1 + \sum_{\ell=1}^\infty 2^{2\ell \gamma_1} \mathbb{P}\left( \sup_{0\le t < 2^{-\ell +1}}\sup_{k=1,2,3} t^{\eta'} \left\| \Wick{X_{0,t}^k}_{C_N}  \right\|_{\mathcal{C}^{-\varepsilon'}(w_{\kappa'})}  \ge K   \right)\notag \\
    &\le 1+ \sum_{\ell=1}^\infty 2^{2\ell \gamma_1}\frac{1}{K^P} \mathbb{E} \left[ \sup_{0\le t < 2^{-\ell +1}} \sup_{k=1,2,3} t^{P\eta'} \left\| \Wick{X_{0,t}^k}_{C_N}  \right\|_{\mathcal{C}^{-\varepsilon'}(w_{\kappa'})}^P  \right] \lesssim 1+ \sum_{\ell=1}^\infty 2^{2\ell \gamma_1}\left(2^{-\ell +1}\right)^{P\eta'/2}\notag\\
    &< \infty \notag
    \end{align}
    provided that $P\ge1$ is chosen sufficiently large.

We now estimate $\mathbb{E}\left[ e^{2\gamma_2N(t)} \right]$ following the argument in \cite[Section 3.2]{KT25}.
Let $\mathcal{F}_t$ be the filtration generated by $\{W_t\}_{t\ge 0}$.
Then the process $\{X_{t,t+ \cdot}\}$ is independent of $\mathcal{F}_{t}$. The same is true for any stopping time by the strong Markov property and hence $\{X_{\tau_n,\tau_n+ \cdot}\}$ is independent of $\mathcal{F}_{\tau_n}$. This implies that $\tau_{n+1} - \tau_n$ is independent of $\mathcal{F}_{\tau_n}$ and $\tau_{k+1}-\tau_k$, $k=0,1,2,\cdots$ are identically distributed by their definition.
Therefore, using $\{N(t) \ge n\} = \{ \tau_n\le t  \} = \{ \sum_{k=1}^n(\tau_k - \tau_{k-1}) \le t \}$, we obtain
\begin{align}
 \mathbb{E}\left[ e^{2\gamma_2N(t)} \right] &\lesssim 1+  \sum_{n=1}^\infty e^{2\gamma_2n }\mathbb{P}\left( N(t) \ge n \right) = 1+  \sum_{n=1}^\infty e^{2\gamma_2n }\mathbb{P}\left( \sum_{k=1}^n(\tau_k - \tau_{k-1}) \le t \right) \notag\\
 &= 1+  \sum_{n=1}^\infty e^{2\gamma_2n }\mathbb{P}\left( e^{-c\sum_{k=1}^n(\tau_k - \tau_{k-1})} \ge e^{-ct} \right) \notag\\
 &\le 1+  \sum_{n=1}^\infty e^{2\gamma_2n +ct} \mathbb{E} \left[e^{-c\sum_{k=1}^n(\tau_k - \tau_{k-1})}   \right] = 1+  e^{ct}\sum_{n=1}^\infty e^{2\gamma_2n} \mathbb{E} \left[e^{-c(\tau_1 - \tau_0)}   \right]^{n} 
\end{align}
for any $c>0$.
We also have
\begin{align}
\mathbb{E} \left[e^{-c(\tau_1- \tau_0)}   \right] \le e^{-c\nu} + \mathbb{P}(\tilde{\tau}_1 - \tau_0\le \nu) \le e^{-c\nu} + \frac{1}4,    
\end{align}
Therefore, taking $c = \frac{2\log 2}{\nu}$, we obtain
\begin{align}
\mathbb{E}\left[ e^{2\gamma_2N(t)} \right] &\lesssim 1+  e^{\frac{2\log 2}{\nu} t}\sum_{n=1}^\infty   \left(  e^{2\gamma_2 } \cdot \frac{1}2\right)^n \lesssim  e^{\frac{2\log 2}{\nu} t}
\end{align}
since $3\gamma_2 \le \log 2$. \qedhere
\end{proof}

We can finally give the proof of Lemma~\ref{lem:linear} and finish this section.

\begin{proof}[Proof of Lemma~\ref{lem:linear}]
Fix $\eta>0$ and choose $K>0$ so large that 
\begin{align} 
    \mathbb{P} \left( \tilde{\tau}_1  \le \frac{2\log 2}\eta  \right) \le \frac{1}4
\end{align}
and let  
\begin{align}
    \nu \coloneqq \frac{2\log 2}\eta \wedge \left( \frac{\log 2}{6\lambda^r \mathfrak{p}_2(K)}   \right)^{1/\alpha} \wedge \frac{1}{3p\theta}.
\end{align}
Then, we have
$$  \mathbb{P} \left( \tilde{\tau}_1  \le \nu \right) \le \frac{1}4 $$
and
\begin{align}
    3 \left\{ \frac{p(r-1)}r \log 2 + \lambda^{r/2} \mathfrak{p}_2(K)\nu^\alpha    \right\} \le \log 2
\end{align}
provided that $r$ is so small that $\frac{p(r-1)}r \le \frac{1}6$. The conclusion then follows from Lemmas~\ref{lem:stopping1} and \ref{lem:stopping2}, together with the bound 
\begin{align*}
    \frac{2\log 2}\nu = 2\log 2 \left\{\frac{\eta}{2\log 2} \vee \left( \frac{6\lambda^{r/2} \mathfrak{p}_2(K)}{\log 2}   \right)^{1/\alpha} \vee 3p\theta  \right\} \le \eta + 2\log 2 \left( \frac{6\lambda^{r/2} \mathfrak{p}_2(K)}{\log 2}   \right)^{1/\alpha}+3p\theta .
\end{align*}

\end{proof}

\section{Analysis of the time-truncated resolvents}\label{sec:resolvents}
In this section, we prove Lemma~\ref{lem:resolvent}. The proof uses a standard approximation argument. Similar arguments can be found in the proofs of \cite[Lemma 3.4]{DGT22} and \cite[Proposition 3.5]{DGT22}.

\subsection{Differentiability of the stochastic flow}
The following result will be used to study the regularity of the time-truncated resolvents $G_T^{(N)}$ in the next subsection.
 
\begin{lem}\label{lem:derivative-flow}
Let $\kappa>1$, $T>0$ and $V\coloneqq H^2(w_\kappa)$.
Let $\Phi^{(N)}(\psi)$ be the solution to \eqref{eq:SQ} with initial condition $\Phi^{(N)}_0= \psi$.
    Then the map
    \begin{align}
       V\ni \psi \mapsto \Phi^{(N)}(\psi) \in C([0,T]; L^2(w_{4\kappa}))  \eqqcolon C_TL^2(w_{4\kappa})
    \end{align}
is twice continuously Fr\'echet differentiable. Its derivatives
\begin{align}
    J_t^h(\psi) \coloneqq D\Phi_t^{(N)}(\psi)[h],\ \ \ J_t^{h,k}(\psi) \coloneqq D^2 \Phi_t^{(N)}(\psi)[h,k]
\end{align}
solve 
\begin{align}
    \begin{cases}
        \partial_t J_t^h = \frac{1}2(\Delta-1)J_t^h -\frac{3}2\lambda\Wick{(\Phi_t^{(N)})^2}_{C_N}J_t^h \\
        J_0^h=h, \label{eq:JN1}
    \end{cases}
\end{align}
\begin{align}
    \begin{cases}
        \partial_t J_t^{h,k} = \frac{1}2(\Delta-1)J_t^{h,k} -\frac{3}2\lambda\Wick{(\Phi_t^{(N)})^2}_{C_N}J_t^{h,k}  -3\lambda \Phi^{(N)}_t J_t^h J_t^k \\
        J_0^{h,k}=0, \label{eq:JN2}
    \end{cases}
\end{align}
respectively.
Moreover, there exists $D_N>0$ such that
\begin{align}
 \|  J^h(\psi)\|_{C_TL^r(w_\kappa)} \le e^{D_N T}  \| h \|_{L^r(w_\kappa)} \label{eq:bound-JN}
\end{align}
for every $1\le r <\infty$ and
\begin{align}
 \|  J^{h,k}(\psi)\|_{C_TL^2(w_{2\kappa})} \le e^{D_N T} \| h \|_{L^4(w_\kappa)} \| k \|_{L^4(w_{\kappa})} . \label{eq:bounds-J2}
\end{align}
These bounds hold uniformly in $\psi \in V$.
\end{lem}

\begin{proof}
The proof follows a standard argument; see, for example, \cite{DGT22,Cer01} for similar arguments.
To establish \eqref{eq:bound-JN} and \eqref{eq:bounds-J2}, we consider equations of the form
\begin{align}
    \partial_t v = \frac{1}2\Delta v + E_N v - a_0 \lambda ( \Phi^{(N)} )^2 v + \Phi^{(N)}f_1 + f_0
\end{align}
for some $E_N, a_0>0$ and $f_0,f_1$. Multiplying this equation by $w_\alpha^r |v|^{r-2}v$ ($\alpha>0$) and integrating by parts, we obtain
\begin{align}
    &\frac{1}r \frac{d}{dt} \| v\|^r_{L^r(w_\alpha)} + \frac{r-1}2 \int w_\alpha^r |v|^{r-2} |\nabla v|^2 + a_0\lambda \int w_\alpha^r ( \Phi^{(N)}  )^2 |v|^r \notag\\
    &= -\frac{1}2 \int \left\{ \nabla v \cdot \nabla(w_\alpha^r)  \right\} |v|^{r-2}v + E_N \| v \|^r_{L^r(w_\alpha)} + \int \left\{ \Phi^{(N)} f_1 +  f_0  \right\} w_\alpha^r |v|^{r-2}v.
\end{align}
The first term on the right-hand side is bounded by
\begin{align}
  \frac{r}2 \int w_\alpha^r |v|^{r-1}|\nabla v| |\nabla \log w_\alpha | \le \frac{r-1}4 \int w_\alpha^r |v|^{r-2} |\nabla v|^2 + \frac{r^2}{4(r-1)} \| \nabla \log w_\alpha \|_{L^\infty}^2 \| v \|^r_{L^r(w_\alpha)}.  
\end{align}
The third term is bounded by
\begin{align}
    &\int \left| \Phi^{(N)} |v|^{r/2} w_\alpha^{r/2}  \right| \cdot  \left| f_1 |v|^{r/2 -1} w_\alpha^{r/2}  \right| +  \| f_0 \|_{L^r(w_\alpha)}\| v\|_{L^r(w_\alpha)}^{r-1}\notag\\
    &\le \frac{a_0 \lambda}2 \int w_\alpha^r (\Phi^{(N)}  )^2 |v|^r  + C_{a_0,\lambda} \| f_1 \|_{L^r(w_\alpha)}^2\| v\|_{L^r(w_\alpha)}^{r-2} + \| f_0 \|_{L^r(w_\alpha)}\| v\|_{L^r(w_\alpha)}^{r-1}.
\end{align}
Combining these estimates, we obtain
\begin{align}
\frac{1}r \frac{d}{dt} \| v\|^r_{L^r(w_\alpha)}  \le C_{N,\alpha} \| v\|_{L^r(w_\alpha)}^r +C_{a_0,\lambda} \| f_1 \|_{L^r(w_\alpha)}^2\| v\|_{L^r(w_\alpha)}^{r-2} + \| f_0 \|_{L^r(w_\alpha)}\| v\|_{L^r(w_\alpha)}^{r-1}.  \label{eq:bound-v-linear}
\end{align}
Dividing it by $\| v\|_{L^r(w_\alpha)}^{r-2}$ and applying Gronwall's inequality, we obtain
\begin{align}
 \| v\|_{C_TL^r(w_\alpha)}  \le C_{N,T,r,\alpha,a_0,\lambda} \left( \| v_0\|_{L^r(w_\alpha)} +  \| f_1 \|_{L^2_TL^r(w_\alpha)} +  \| f_0 \|_{L_T^1 L^r(w_\alpha)}  \right)  .\label{eq:bound-v-linear-2}
\end{align}
Applying \eqref{eq:bound-v-linear-2} with $f_0=f_1=0$ to \eqref{eq:JN1}, we obtain \eqref{eq:bound-JN}. Moreover, applying \eqref{eq:bound-v-linear-2} to \eqref{eq:JN2} with $r=2$ gives
\begin{align}
\|  J^{h,k}(\psi)\|_{C_TL^2(w_{2\kappa})}\lesssim \|J^h J^k   \|_{L_T^2L^2(w_{2\kappa})}  \lesssim \| J^h \|_{C_TL^4(w_\kappa)} \| J^k \|_{C_TL^4(w_{\kappa})} ,  
\end{align}
which implies \eqref{eq:bounds-J2} in combination with \eqref{eq:bound-JN}.

We now turn to the proof of differentiability. Let $d^{h,\varepsilon} \coloneqq \Phi^{(N)}(\psi + \varepsilon h)- \Phi^{(N)}(\psi)$, $h\in V, \varepsilon\in \mathbb{R}$. This difference solves 
\begin{align}
    \begin{cases}
        \begin{aligned}[t]
        \partial_t d_t^{h,\varepsilon}
        &= \frac{1}2(\Delta-1) d_t^{h,\varepsilon} \\
        &\quad + \frac{\lambda}2 \left(3C_N - \left\{
          \Phi^{(N)}(\psi + \varepsilon h)^2
          + \Phi^{(N)}(\psi + \varepsilon h)\Phi^{(N)}(\psi )
          + \Phi^{(N)}(\psi)^2
          \right\}\right) d_t^{h,\varepsilon}
        \end{aligned} \\[2pt]
        d_t^{h,\varepsilon} = \varepsilon h. \label{eq:d1}
    \end{cases}
\end{align}
Since $\Phi^{(N)}(\psi + \varepsilon h)^2 + \Phi^{(N)}(\psi + \varepsilon h)\Phi^{(N)}(\psi ) + \Phi^{(N)}(\psi )^2 \ge 0$, the same energy computation as above, together with the Sobolev embedding, yields
\begin{align}
  \|  d^{h,\varepsilon} \|_{C_TL^r(w_\kappa)} \le  e^{D_N T} \|\varepsilon h \|_{L^r(w_\kappa)} \le C_r  e^{D_N T} \|\varepsilon h \|_{V} \label{eq:ac1}
\end{align}
for some $D_N>0$.
Moreover, setting $R \coloneqq d^{h,\varepsilon} - J^{\varepsilon h}$, we have
\begin{align}
    \begin{cases}
        \begin{aligned}[t]
        \partial_t R
        &= \frac{1}2(\Delta-1) R
          + \frac{3\lambda}2 \left(C_N - (\Phi_t^{(N)})^2\right) R  - \frac{3\lambda}2 \Phi^{(N)} (d^{h,\varepsilon})^2
          - \frac{\lambda}2 (d^{h,\varepsilon})^3
        \end{aligned} \\[2pt]
        R_0 = 0.
    \end{cases}
\end{align}
Therefore, applying the bound \eqref{eq:bound-v-linear-2}, we obtain
\begin{align}
    \begin{aligned}
    \| R \|_{C_TL^r(w_{3\kappa})}
    &\lesssim \| (d^{h,\varepsilon})^2 \|_{L^2_TL^r(w_{3\kappa})}
      + \| (d^{h,\varepsilon})^3 \|_{L_T^1L^r(w_{3\kappa})} \\
    &\lesssim \| d^{h,\varepsilon} \|_{C_TL^{2r}(w_\kappa)}^2
      + \| d^{h,\varepsilon} \|_{C_TL^{3r}(w_\kappa)}^3 \\
    &\lesssim  \|\varepsilon h \|_V^2
      +  \|\varepsilon h \|_V^3.
    \end{aligned}
    \label{eq:bound-R}
\end{align}
This proves Fr\'echet differentiability.
Next, $d^{h,k,\varepsilon} \coloneqq J^h(\psi + \varepsilon k) - J^h(\psi)$ solves 
\begin{align}
    \begin{cases}
        \begin{aligned}[t]
        \partial_t d_t^{h,k,\varepsilon}
        &= \frac{1}2(\Delta-1) d_t^{h,k,\varepsilon}
          + \frac{3\lambda}2 \left(C_N - (\Phi_t^{(N)})^2\right)
            d_t^{h,k,\varepsilon} \\
        &\quad - 3\lambda \Phi^{(N)} d^{k,\varepsilon}
            J^h(\psi + \varepsilon k)
          - \frac{3\lambda}2 (d^{k,\varepsilon})^2
            J^h(\psi + \varepsilon k)
        \end{aligned} \\[2pt]
        d_t^{h,k,\varepsilon} = 0. \label{eq:d2}
    \end{cases}
\end{align}
Therefore, a computation similar to the one above, together with the bounds \eqref{eq:bound-JN} and \eqref{eq:ac1}, yields
\begin{align}
  \|  d^{h,k,\varepsilon} \|_{C_TL^2(w_{3\kappa})} &\lesssim \| d^{k,\varepsilon} J^h(\psi+ \varepsilon k) \|_{L_T^2 L^r(w_{3\kappa})} + \| (d^{k,\varepsilon})^2 J^h(\psi+ \varepsilon k) \|_{L_T^1 L^r(w_{3\kappa})} \notag\\
  &\lesssim  \| d^{k,\varepsilon} \|_{C_TL^{2r}(w_\kappa)} \| J^h(\psi+ \varepsilon k) \|_{C_TL^{2r}(w_\kappa)} +  \| d^{k,\varepsilon} \|_{C_TL^{3r}(w_\kappa)}^2 \| J^h(\psi+ \varepsilon k) \|_{C_TL^{3r}(w_\kappa)}    \notag\\
  &\lesssim \left(  \| \varepsilon k \|_V +  \| \varepsilon k \|_V^2  \right) \| \varepsilon h \|_V . \label{eq:bound-dhk}
\end{align}
Moreover, letting $R' \coloneqq d^{h,k,\varepsilon} - J^{\varepsilon k, h}$, we have
\begin{align}
    \begin{cases}
        \begin{aligned}[t]
        \partial_t R'
        &= \frac{1}2(\Delta-1) R'
          + \frac{3\lambda}2 \left(C_N - (\Phi_t^{(N)})^2\right) R' \\
        &\quad - 3\lambda \Phi^{(N)}
          \left(R_k J^h(\psi + \varepsilon k)
            + J^{\varepsilon k} d^{h,k,\varepsilon}\right)  - \frac{3\lambda}2 (d^{k,\varepsilon})^2
          J^h(\psi + \varepsilon k)
        \end{aligned} \\[2pt]
        R_0' = 0,
    \end{cases}
\end{align}
where $R_k \coloneqq d^{k,\varepsilon} - J^{\varepsilon k}$.
Thus, \eqref{eq:bound-v-linear-2}, \eqref{eq:bound-R} and \eqref{eq:bound-dhk} yield
\begin{align}
    \| R'\|_{C_T L^r(w_{4\kappa})} \lesssim \left( \| \varepsilon k \|_V^2 +  \| \varepsilon k \|_V^3   \right) \| \varepsilon h \|_V.
\end{align}
This proves that the solution map is twice Fr\'echet differentiable. Continuity of the derivatives follows from a similar computation.\qedhere

\end{proof}

\subsection{Proof of Lemma~\ref{lem:resolvent}}

Fix $\beta_0>1$ and $\kappa>1$. Let $E=H^{-\beta_0}(w_\kappa)$ and let $\mathfrak{F}C^2_{\mathrm{pol}} (E)$ be the space of all functions $F$ of the form \eqref{eq:cylind} for which $f\in C^2(\mathbb{R}^m)$ and $f$, $\nabla f, \nabla^2f$ have polynomial growth.
We first prove that functions in this class belong to the domain of $\bar{\mathcal{L}}_p$.
\begin{lem}\label{lem:domain-1}
Let $F\in \mathfrak{F}C^2_{\mathrm{pol}} (E)$ be of the form \eqref{eq:cylind} for $f$ as above. Then $F \in \mbox{Dom}(\bar{\mathcal{L}}_p )$ and $\bar{\mathcal{L}}_pF$ is given by \eqref{eq:op-cyl}.
    
\end{lem}

\begin{proof}
Write $F(\phi)=f(\mathbf l(\phi))$, where
$\mathbf l(\phi)=(\langle\phi,l_i\rangle)_{i=1}^m$.
Choose $\chi\in C_c^\infty(\mathbb R^m)$ such that $\chi=1$ near the origin,
and put $f_R(y)=\chi(y/R)f(y)$, $R\ge1$.
For $|\alpha|\le2$, we have
\begin{equation}\label{eq:bound-f-derivative}
 \partial^\alpha f_R(y)\to\partial^\alpha f(y),\qquad
 |\partial^\alpha f_R(y)| \lesssim (1+|y|)^M,
\end{equation}
where the bound on $\partial^\alpha f_R$ are uniform in $R$. 
Moreover, letting $\{\varphi_\epsilon\}_{\epsilon>0}$ be a standard mollifier on $\mathbb{R}^m$, we have for $f_{R,\epsilon}=\varphi_\epsilon*f_R \in C_c^\infty(\mathbb{R}^m)$ that
\begin{equation}
\| \partial^\alpha f_{R,\epsilon}-\partial^\alpha f_R\|_{L^\infty(\mathbb{R}^m)}\to 0,\qquad |\alpha|\le 2 
\end{equation}
as $\epsilon\searrow0$.
Therefore, we can choose a sequence
$f_n\in C_c^\infty(\mathbb R^m)$ such that 
\begin{equation}
 \partial^\alpha f_n(y)\to\partial^\alpha f(y),\qquad
 |\partial^\alpha f_n(y)| \lesssim (1+|y|)^M,
\end{equation}
for any $|\alpha|\le 2$, where the bound on the derivatives are uniform in $n$.

Set $F_n=f_n\circ\mathbf l\in\mathfrak F C_b^\infty$. It follows from
\eqref{eq:op-cyl} that
\begin{align*}
 |F_n(\phi)|+|\mathcal LF_n(\phi)|
 &\le C(1+|\mathbf l(\phi)|)^M
       \left(1+\sum_{i=1}^m|\langle (1-\Delta)\phi + \lambda \Wick{\phi}^3  ,l_i\rangle|\right)\\
 &\le C'(1+\|\phi\|_E)^M(1+\| (1-\Delta)\phi + \lambda \Wick{\phi}^3  \|_{H^{-\beta_0-2}}(w_\kappa)).
\end{align*}
The last expression belongs to $L^p(\mu^{(\lambda)})$ for every $1\le p<\infty$ by Lemma~\ref{lem:phi42-wick}. Thus the dominated convergence theorem gives convergence of
$F_n$ and $\mathcal LF_n$ in $L^p$ to $F$ and the expression
\eqref{eq:op-cyl}, respectively. The closedness of $\bar{\mathcal
{L}}_p$ proves the desired assertion.
\end{proof}

The following lemma will be used in the approximation argument below.

\begin{lem}\label{lem:res-spatial-approximation}
There exist real finite-rank operators $Q_n$ with ranges contained in
$C_c^\infty(\mathbb R^2)$ such that
\begin{equation}\label{eq:res-Q-properties}
 \sup_n\big(\|Q_n\|_{\mathcal L(E)}+
              \|Q_n\|_{\mathcal L(H^{-\beta_0-2}(w_\kappa))}\big)<\infty,
 \qquad Q_nu\to u\quad\hbox{in }E\hbox{ and }H^{-\beta_0-2}(w_\kappa) ,
\end{equation}
where $\mathcal{L}(E)$ denotes the space of bounded operators on $E$.
They can be chosen self-adjoint contractions on $L^2(\mathbb{R}^2)$, converging strongly
to the identity on $L^2(\mathbb{R}^2)$. In particular, 
\begin{equation}\label{eq:res-HS-convergence}
 \|Q_n\iota-\iota\|_{L_{(2)}( L^2(\mathbb{R}^2) ,E)}\to 0
\end{equation}
because the inclusion $\iota:L^2(\mathbb{R}^2) \hookrightarrow E = H^{-\beta_0}(w_\kappa)$ is Hilbert--Schmidt when $\beta_0,\kappa>1$.
Here $L_{(2)}$ and $L_{(1)}$ denote the Hilbert--Schmidt and
trace-class ideals, respectively.
\end{lem}
\begin{proof}
Choose even cutoffs $f,\vartheta\in C_c^\infty(\mathbb R^2)$,
with values in $[0,1]$, equal to one near the origin, and with
$\operatorname{supp}f\subset(-\pi,\pi)^2$.
On the torus $\mathbb T_L^2=\mathbb R^2/(2\pi L\mathbb Z)^2$, we consider the normalized Fourier basis
$e_{L,\ell}(x)=(2\pi L)^{-1}e^{i\ell\cdot x/L}$.
For $R,L\ge1$, define
\begin{equation}\label{eq:res-Q-definition}
 Q_{R,L}u=f_L\sum_{\ell\in\mathbb Z^2}
   \vartheta\!\left(\frac{\ell}{RL}\right)
   \langle f_Lu,\overline{e_{L,\ell}}\rangle e_{L,\ell},
 \qquad f_L(x)=f(x/L).
\end{equation}
Note that the sum is finite and this is a self-adjoint contraction on $L^2(\mathbb{R}^2)$. Let $K_{R,L}$ be the kernel associated with the Fourier multiplier appearing in the definition of $Q_{R,L}$. For any $\alpha\in \mathbb{R}$ and $\gamma\in \mathbb{N}_0^2$, it follows from $w_\alpha(x) \lesssim w_{|\alpha|}(x-y)w_\alpha(y)$, the rapid decay of $K_{R,L}$, the Leibniz rule, and Young's convolution inequality that
\begin{align}
    &\|  \partial^\gamma Q_{R,L} u   \|_{L^2(w_\alpha)}^2 \lesssim \sum_{\gamma_1+\gamma_2= \gamma} \int |\partial^{\gamma_1}f_L |^2 \left(  K_{R,L}*\partial^{\gamma_2} u    \right)^2(x) w_{2\alpha}(x) dx \notag\\* 
    &\lesssim \sup_{|\gamma_2|\le |\gamma|} \int \int \int  K_{R,L}(x-y_1) w_{|\alpha|}(x-y_1)  K_{R,L}(x-y_2) w_{|\alpha|}(x-y_2) \notag\\
    &\qquad \qquad \qquad \qquad \qquad\qquad\times\partial^{\gamma_2} u(y_1) w_\alpha(y_1)\partial^{\gamma_2} u(y_2) w_\alpha(y_2) dxdy_1 dy_2 \notag\\*
    &\lesssim \sup_{|\gamma_2|\le |\gamma|}  \int \left\{ \left( w_{|\alpha|} K_{R,L}  \right) * \left( w_\alpha \partial^{\gamma_2} u   \right)   (x)    \right\}^2  dx        \notag \\*
    &\lesssim \sup_{|\gamma_2|\le |\gamma|} \|\partial^{\gamma_2} u \|_{L^2(w_\alpha)}^2
\end{align}
uniformly in $R,L>0$. This implies $\| Q_{R,L} v\|_{H^m(w_\alpha)} \lesssim \| v\|_{H^m(w_\alpha)}$ uniformly in $R,L>0$ for any $m\in \mathbb{N}_0$. As a consequence,
\begin{align}
    \left|\langle  Q_{R,L}u, v \rangle  \right| = \left|\langle  u,  Q_{R,L} v \rangle  \right|  &\lesssim \| u \|_{H^{-m}(w_\alpha)}   \| Q_{R,L} v \|_{H^{m}(w_{-\alpha})} \notag\\
    &\lesssim \| u \|_{H^{-m}(w_\alpha)}   \|  v \|_{H^{m}(w_{-\alpha})},
\end{align}
which implies the bound \eqref{eq:res-Q-properties} by duality and interpolation.\qedhere

\end{proof}

\begin{lem}\label{lem:domain-extension}
Let $G:E \rightarrow \mathbb{R}$ be twice continuously Fr\'echet differentiable. Assume $DG \in C(E;H^{\beta_0+2}(w_{-\kappa}) )$ and, for some $M>0$,
\begin{align}
    |G(\phi)| + \| DG(\phi) \|_{H^{\beta_0+2}(w_{-\kappa})} + \| D^2 G(\phi) \|_{\mathcal{L}^{(2)}(E;\mathbb{R})} \lesssim \left( 1+ \| \phi \|_{E}  \right)^M, \label{eq:poly-bound-G}
\end{align}
where $\mathcal{L}^{(2)}(E;\mathbb{R})$ denotes the space of bounded bilinear operators from $E$ to $\mathbb{R}$.
Then $G \in \mathcal{D}(\bar{\mathcal{L}}_p )$ for any $1\le p<\infty$, and
\begin{align}
    \bar{\mathcal{L} }_p G(\phi) = \frac{1}2\operatorname{Tr}_{L^2(\mathbb{R}^2)} D^2G(\phi)   - \frac{1}2DG(\phi)\left[ (1-\Delta)\phi+\lambda\Wick{\phi^3}   \right].
\end{align}
    
\end{lem}

\begin{proof}
Using the operators from Lemma~\ref{lem:res-spatial-approximation}, consider the approximations
\begin{align}
    G_{n} (\phi) \coloneqq G \left( Q_{n}  \phi  \right).
\end{align}
Then $G_n\in \mathfrak{F}C^2_{\mathrm{pol}} (E)$, so Lemma~\ref{lem:domain-1} applies. 
By the closedness of $\bar{\mathcal{L}_p}$, the conclusion follows once we prove
\begin{align}
    \begin{aligned}
    \lim_{n\to\infty} G_n &= G, \\
    \lim_{n\to\infty} \bar{\mathcal{L}}_p G_n
    &= \frac{1}2\operatorname{Tr}_{L^2(\mathbb{R}^2)}
       \left(D^2G(\phi)\right)  - \frac{1}2DG(\phi)
       \left[(1-\Delta)\phi + \lambda\Wick{\phi^3}\right],
    \end{aligned}
\end{align}
where both limits are taken in $L^p(\mu^{(\lambda)})$. The first convergence follows immediately from Lemma~\ref{lem:res-spatial-approximation} and the bound \eqref{eq:poly-bound-G}. We therefore focus on the second one.
The chain rule yields
\begin{align}
    DG_n(\phi)[h]=DG(Q_n\phi)[Q_nh],\qquad
 D^2G_n(\phi)[h,k]=D^2G(Q_n\phi)[Q_nh,Q_nk].
\end{align}

For the drift term, we use the decomposition
\begin{align}\label{eq:first-derivative}
 &DG(Q_n\phi)[Q_nb(\phi)]-DG(\phi)[b(\phi)]\notag\\
 &\quad=\big(DG(Q_n\phi)-DG(\phi)\big)[Q_nb(\phi)]
       +DG(\phi)[(Q_n-I)b(\phi)],
\end{align}
where 
\begin{align}
    b(\phi) \coloneqq (1-\Delta)\phi+\lambda\Wick{\phi^3}.
\end{align}
The two terms converge pointwise to zero by continuity of $DG$ and
strong convergence of $Q_n$ on $H^{-\beta_0-2}(w_\kappa)$ (Lemma~\ref{lem:res-spatial-approximation}). They are dominated by 
$(1+\|\phi\|_E)^M\|b(\phi)\|_{H^{-\beta_0-2}(w_\kappa)}$ by \eqref{eq:poly-bound-G}, and this bound belongs to $L^p(\mu^{(\lambda)})$ by Lemma~\ref{lem:phi42-wick}.

We now turn to the trace term. For the inclusion map $\iota:L^2(\mathbb{R}^2) \hookrightarrow E = H^{-\beta_0}(w_\kappa)$, it suffices to prove the convergence
\begin{align}
    \left\| \tilde{Q}_n^* \left( D^2G \circ Q_n \right) \tilde{Q}_n -  \iota^* D^2G \iota    \right\|_{L_{(1)}(L^2(\mathbb{R}^2) )} \to 0
\end{align}
where $\tilde{Q}_n \coloneqq Q_n \iota$. Note that the bound
\begin{align}
    \| \iota^* D^2G \iota \|_{L_{(1)}(L^2(\mathbb{R}^2) )} \le \| D^2G \|_{\mathcal{L}(E)} \| \iota \|_{L_{(2)}(L^2(\mathbb{R}^2):E )}^2
\end{align}
holds, where we identified $\mathcal{L}^{(2)}(E;\mathbb{R})$ with $\mathcal{L}(E)$ in the usual way.
From the identity
\begin{align}
    \begin{aligned}
    &\tilde{Q}_n^* \left(D^2G \circ Q_n\right)\tilde{Q}_n
       - \iota^* D^2G \iota  = \tilde{Q}_n^* \left(D^2G \circ Q_n - D^2G\right)\tilde{Q}_n \\
    &\qquad \qquad\qquad\qquad\qquad\qquad\qquad + (\tilde{Q}_n - \iota)^* D^2G\tilde{Q}_n
       + \iota^* D^2G(\tilde{Q}_n - \iota)
    \end{aligned}
\end{align}
and Lemma~\ref{lem:res-spatial-approximation},
we obtain the pointwise convergence
\begin{align}
&\left\| \tilde{Q}_n^* \left( D^2G \circ Q_n \right) \tilde{Q}_n -  \iota^* D^2G \iota    \right\|_{L_{(1)}(L^2(\mathbb{R}^2) )}\notag\\*
&\le  \| D^2G \circ Q_n - D^2 G    \|_{\mathcal{L}(E)} \| \tilde{Q}_n \|^2_{L_{(2)}(L^2(\mathbb{R}^2):E )} \notag\\*
&\ + \| D^2 G\|_{\mathcal{L}(E)} \left( \| \tilde{Q}_n \|_{L_{(2)}(L^2(\mathbb{R}^2):E )} +  \| \iota \|_{L_{(2)}(L^2(\mathbb{R}^2):E )}    \right) \| \tilde{Q}_n - \iota  \|_{L_{(2)}(L^2(\mathbb{R}^2):E )} \to 0.
\end{align}
Since these terms are bounded by $(1+\|\phi\|_E)^M$ by \eqref{eq:poly-bound-G}, convergence also holds in $L^p(\mu^{(\lambda)})$ by the dominated convergence theorem. This completes the proof.\qedhere

\end{proof}

\begin{lem}\label{lem:domain-Gtilde}
Fix $T>0$ and $N\in \mathbb{N}$. Recall the following definition from Section~\ref{sec:proof}:
\begin{align}
    G_T^{(N)}(\psi)=\int_0^T e^{-at}P_t^{(N)}F(\psi)\,dt,
 \qquad \tilde G_T^{(N)}(\phi)=G_T^{(N)}(\Pi_N\phi),
 \qquad F\in\mathfrak F C_b^\infty.
\end{align}
Then $\tilde{G}^{(N)}_T$ satisfies the assumptions of Lemma~\ref{lem:domain-extension}. In particular, we have $\tilde{G}^{(N)}_T \in \mathcal{D}(\bar{\mathcal{L}}_p ) $ for any $1\le p<\infty$.
\end{lem}

\begin{proof}
By Lemma~\ref{lem:derivative-flow}, we have $G_T^{(N)} \in  C^2(H^2(w_\kappa ) )$ and its derivatives are given by
\begin{align}
    DG^{(N)}_T(\phi) [h] = \int^T_0 e^{-at} \mathbb{E} \left[ DF\left( \Phi_t^{(N)}(\phi)  \right) \left[ J_t^h(\phi)  \right]   \right]\, dt
\end{align}
and
\begin{align}
    D^2G^{(N)}_T(\phi) [h,k]
    &= \int^T_0 e^{-at} \mathbb{E} \left[ D^2F\left( \Phi_t^{(N)}(\phi)  \right) \left[ J_t^h(\phi), J^k_t(\phi)  \right]   \right]\, dt \notag\\
    &\ \ \ \ \ \ \ \ \ + \int^T_0 e^{-at} \mathbb{E} \left[ DF\left( \Phi_t^{(N)}(\phi)  \right) \left[ J_t^{h,k}(\phi)  \right]   \right]\, dt,
\end{align}
where the bounds given in the proof of Lemma~\ref{lem:derivative-flow}  justify differentiation under both the expectation and the time integral.
It also follows from Lemma~\ref{lem:derivative-flow} that
\begin{align}
    \sup_{\psi\in H^2(w_\kappa)} \left(  \| D G_T^{(N)}(\psi)   \|_{H^{-2}(w_{-\kappa})} + \| D^2 G_T^{(N)}(\psi)  \|_{\mathcal{L}^{(2)}(H^2(w_\kappa);\mathbb{R})}    \right) <\infty . \label{eq:bound-G-derivatives}
\end{align}
The chain rule and $\Pi_N \in \mathcal{L} \left( H^{-\beta_0-2}(w_\kappa), H^2(w_\kappa)  \right)$ yield $\tilde{G}_T^{(N)}\in C^2(E;\mathbb{R})$, 
\begin{align}
    \begin{aligned}
    D\tilde{G}^{(N)}_T(\phi)[h]
    &= DG^{(N)}_T(\Pi_N\phi)[\Pi_N h], \\
    D^2\tilde{G}^{(N)}_T(\phi)[h,k]
    &= D^2G^{(N)}_T(\Pi_N\phi)[\Pi_N h,\Pi_N k],
    \end{aligned}
\end{align}
and $D\tilde{G}_T^{(N)} \in C(E;H^{\beta_0+2}(w_{-\kappa}))$.
For the first derivative, \eqref{eq:bound-G-derivatives} and Lemma~\ref{lem:appr-1} yield 
\begin{align}
    \begin{aligned}
    \sup_{\phi\in E}
       \|D\tilde{G}_T^{(N)}(\phi)\|_{H^{\beta_0+2}(w_{-\kappa})}  &= \sup_{\phi\in E}
       \sup_{\|h\|_{H^{-\beta_0-2}(w_\kappa)}\le1}
       |D\tilde{G}_T^{(N)}(\Pi_N\phi)[\Pi_N h]| \\
    &\quad \lesssim \sup_{\psi\in H^2(w_\kappa)}
       \|DG_T^{(N)}(\psi)\|_{H^{-2}(w_{-\kappa})} < \infty.
    \end{aligned}
\end{align}
For the second derivative, we similarly obtain
\begin{align}
    \begin{aligned}
    \sup_{\phi\in E}
       \|D^2\tilde{G}_T^{(N)}(\phi)\|_{\mathcal{L}^{(2)}(E;\mathbb{R})} & = \sup_{\phi\in E}\sup_{\|h\|_E,\|k\|_E\le1}
       |D^2\tilde{G}_T^{(N)}(\Pi_N\phi)[\Pi_N h,\Pi_N k]| \\
    & \lesssim \sup_{\psi\in H^2(w_\kappa)}
       \|D^2G_T^{(N)}(\psi)\|_{\mathcal{L}^{(2)}(H^2(w_\kappa);\mathbb{R})}
       < \infty.
    \end{aligned}
\end{align}
Together with $ |\tilde{G}_T^{(N)}| \le \| F\|_\infty/a$, these bounds prove the assertion.

\end{proof}

We can now prove Lemma~\ref{lem:resolvent}.

\begin{proof}[Proof of Lemma~\ref{lem:resolvent}]
The semigroup property of $P_t^{(N)}$ yields, for any $\psi \in H^2(w_\kappa)$,
\begin{align}
  P_t^{(N)} G^{(N)}_T  = \int^T_0 e^{-as} P_{t+s}^{(N)} F(\psi)\, ds = e^{at}\int^{T+t}_t e^{-as} P_s^{(N)}F(\psi)\, ds. 
\end{align}
Therefore, we obtain
\begin{align}
    \begin{aligned}
    \frac{P_t^{(N)}G_T^{(N)}-G_T^{(N)}}t
    &= \frac{e^{at}-1}t \int_t^{T+t} e^{-as}P_s^{(N)}F\,ds \\
    &\quad - \frac{1}t\int_0^t e^{-as}P_s^{(N)}F\,ds
       - \frac{1}t\int_t^{T+t} e^{-as}P_s^{(N)}F\,ds \\
    &\to aG_T^{(N)}-F-e^{-aT}P_T^{(N)}F
    \end{aligned}
    \label{eq:kolmogorov1}
\end{align}
as $t\searrow0$.
On the other hand, the infinite-dimensional It\^o formula (cf. \cite[Theorem 4.32]{DPZ14}), together with $\Phi^{(N)}(\cdot) \in C([0,T]; H^2(w_\kappa))$ and the regularity property of $G_T^{(N)}$, yields
\begin{align}
    \begin{aligned}
    \frac{P_t^{(N)}G_T^{(N)}(\psi)-G_T^{(N)}(\psi)}t     &\quad = \mathbb{E}\left[\frac{1}t\int_0^t
       \left\{\frac{1}2\operatorname{Tr}_{L^2(\mathbb{R}^2)}
       \left(\Pi_ND^2G_T^{(N)}(\Phi_s^{(N)}(\psi))\Pi_N\right)
       \right.\right. \\
    &\qquad\ \  \left.\left.
       - \frac{1}2 DG_T^{(N)}(\Phi_s^{(N)}(\psi))
       \left[(1-\Delta)\Phi_s^{(N)}(\psi)
          + \lambda\Wick{\Phi_s^{(N)}(\psi)^3}_{C_N}\right]
       \right\}\,ds\right] \\
    &\quad \to \frac{1}2\operatorname{Tr}_{L^2(\mathbb{R}^2)}
       \left(\Pi_ND^2G_T^{(N)}(\psi)\Pi_N\right)  - \frac{1}2DG_T^{(N)}(\psi)
       \left[(1-\Delta)\psi+\lambda\Wick{\psi^3}_{C_N}\right]
    \end{aligned}
    \label{eq:kolmorogov2}
\end{align}
as $t\searrow0$.
From \eqref{eq:kolmogorov1} and \eqref{eq:kolmorogov2}, we get
\begin{align}
&\frac{1}2\operatorname{Tr}_{L^2(\mathbb{R}^2)}\left( \Pi_N D^2G^{(N)}_T(\psi) \Pi_N  \right)   -\frac{1}2DG^{(N)}_T(\psi)\left[  (1-\Delta)\psi + \lambda\Wick{\psi^3}_{C_N}  \right]\notag\\
&=aG^{(N)}_T(\psi)-F(\psi) - e^{-aT}P_T^{(N)} F(\psi)
\end{align}
for any $\psi \in H^2(w_\kappa)$.
Together with the chain rule, this yields
\begin{align}
    \begin{aligned}
    \frac{1}2\operatorname{Tr}_{L^2(\mathbb{R}^2)}
       \left(D^2\tilde{G}_T^{(N)}(\phi)\right) &= \frac{1}2\operatorname{Tr}_{L^2(\mathbb{R}^2)}
       \left(\Pi_ND^2G_T^{(N)}(\Pi_N\phi)\Pi_N\right) \\
    &= \frac{1}2DG_T^{(N)}(\Pi_N\phi)
       \left[(1-\Delta)\Pi_N\phi
         + \lambda\Wick{(\Pi_N\phi)^3}_{C_N}\right] \\
    &\quad + aG_T^{(N)}(\Pi_N\phi)-F(\Pi_N\phi)
       - e^{-aT}P_T^{(N)}F(\Pi_N\phi) \\
    &= \frac{1}2D\tilde{G}_T^{(N)}(\phi)
       \left[(1-\Delta)\phi+\lambda\Wick{\phi^3}\right] \\
    &\quad + a\tilde{G}_T^{(N)}(\phi)-F(\Pi_N\phi)
       - e^{-aT}P_T^{(N)}F(\Pi_N\phi) \\
    &\quad + \frac{\lambda}2DG_T^{(N)}(\Pi_N\phi)
       \left[\Wick{(\Pi_N\phi)^3}_{C_N}-\Pi_N\Wick{\phi^3}\right].
    \end{aligned}
\end{align}
Combining this with Lemmas~\ref{lem:domain-extension} and \ref{lem:domain-Gtilde}, we obtain the desired result.\qedhere
    
\end{proof}

\appendix


\section{A priori estimates}\label{sec:CDFI}

We consider the deterministic equation
\begin{equation}\label{eq:YY}
\begin{cases}
 \partial_t\mathcal Y_t=\frac12(\Delta-1)\mathcal Y_t
 -\nu\big(\mathcal Y_t^3+\mathcal X_t^{(1)}\mathcal Y_t^2
                   +\mathcal X_t^{(2)}\mathcal Y_t+\mathcal X_t^{(3)}\big),\\
 \mathcal Y_0=f,
 \end{cases}
\end{equation}
where $\nu>0$. The coefficients $\mathcal X^{(j)}$ are prescribed, measurable
distribution-valued functions of time. One can recover the equation \eqref{eq:Y:restart} by taking
$\nu=\lambda/2$ and
$$(\mathcal X^{(1)},\mathcal X^{(2)},\mathcal X^{(3)})
=(3X_{s,\cdot},3\Wick{X_{s,\cdot}^2}_{C_N},\Wick{X_{s,\cdot}^3}_{C_N}).$$

The aim of this section is to show the following estimates by applying a weighted energy argument. Similar arguments can be found in \cite{TW18,OTWZ26}. 

\begin{lem}
\label{lem:CDFI-Y}
Fix $0<\eta<1$, $\varepsilon>0$, $\nu^*>0$, and
$\kappa>0$.
There exist $a,\alpha,\kappa_X,\kappa_f>0$ and polynomials
$\mathfrak p_1,\mathfrak p_2,\mathfrak p_3$ with the following property.
Put
\begin{equation}\label{eq:def-mathcalX}
 \mathfrak X_t=\max_{1\le j\le3}\sup_{0<u\le t}
 u^\alpha\|\mathcal X_u^{(j)}\|_{\mathcal C^{-a}(w_{\kappa_X})}, 
 \ (0<t\le1), \qquad F:=\|f\|_{\mathcal C^{-a}(w_{\kappa_f})}.
\end{equation}
Then one has
\begin{align}
\nu^{1/2}  \|\mathcal Y_t\|_{\mathcal{C}^\eta(w_\kappa)}
 &\lesssim t^{-(1+\eta+\varepsilon)/2}
                   \mathfrak p_2(\mathfrak X_t),
 \label{est:YY2}
\end{align}
and
\begin{equation}\label{est:Y3}
 \|\mathcal Y_t\|_{\mathcal{C}^\eta(w_\kappa)}
 \lesssim t^{-(\eta+a+\varepsilon)/2}
                   \mathfrak p_3(F,\mathfrak X_t).
\end{equation}
The first estimate is independent of $f$ and uniform in $0<\nu< \nu^*$. 
\end{lem}

\begin{rem}
    Note that it suffices to consider the case $\nu=1$. Indeed, multiplying the equation~\eqref{eq:YY} by $\nu^{1/2}$, it follows that $\mathcal{Y}_t' \coloneqq \nu^{1/2} \mathcal{Y}_t$ solves
    \begin{equation}
    \begin{cases}
 \partial_t  \mathcal Y_t'=\frac12(\Delta-1) \mathcal Y_t'
 -\big(\mathcal Y_t'^3+\nu^{1/2} \mathcal X_t^{(1)}\mathcal Y_t'^2
                   +\nu \mathcal X_t^{(2)}\mathcal Y_t'+ \nu^{3/2}\mathcal X_t^{(3)}\big),\\
 \mathcal Y_0'=\nu^{1/2} f,
 \end{cases}
\end{equation}
so the proof is reduced to the case $\nu=1$.
    
\end{rem}

We first record the elementary comparison principle used below.
\begin{lem}[\textup{cf.\ \cite[Lemma~3.8]{TW18}}]\label{lem:comparison}
Let $U$ be nonnegative and differentiable, and suppose
$U'+c_1U^m\le c_2$ on $(s,t)$, where $m>1$ and
$c_1,c_2>0$. Then
\[
 U(t)\le
 \left(\frac{c_1(m-1)}2(t-s)\right)^{-1/(m-1)}
 \vee\left(\frac{2c_2}{c_1}\right)^{1/m}.
\]
\end{lem}

Lemma~\ref{lem:comparison} will be combined with the following result. 

\begin{lem}\label{lem:weighted-energy}
Let $r\ge2$, $\tau>2/r$, $0<a<b<1$, and $\zeta>0$ satisfy
\begin{equation}\label{eq:energy-weight-condition}
 r\tau-\zeta \frac{r+2}{1-b}>2.
\end{equation}
Set
\begin{align*}
 Z_j(t)=\|\mathcal X_t^{(j)}\|_{\mathcal C^{-a}(w_\zeta)},\qquad
 q_{r,k}=\frac{r+2}{3-k-b}\quad(0\le k\le2).
\end{align*}
There are $c,C>0$, independent of the solution and coefficients, such that
\begin{equation}\label{est:Y-proof}
 \frac{d}{dt}\int |\mathcal Y_t|^r w_\tau^r\,dx +c\int |\nabla\mathcal Y_t|^2|\mathcal Y_t|^{r-2}w_\tau^r\,dx +c\int |\mathcal Y_t|^{r+2}w_\tau^r\,dx
 \le C\left(1+\sum_{k=0}^2 Z_{3-k}(t)^{q_{r,k}}\right).
\end{equation}
In particular, the left-hand side may be replaced by
$$\frac{d}{dt} \int |\mathcal Y_t|^r w_\tau^r\,dx +c' \left( \int |\mathcal Y_t|^r w_\tau^r\,dx \right)^{1+2/r},$$ 
for a suitable $c'>0$.
\end{lem}

\begin{proof}
Testing \eqref{eq:YY} against
$|\mathcal Y|^{r-2}\mathcal Yw_\tau^r$ and integrating by parts give
\begin{align}\label{eq:energy-identity}
 \frac1r \frac{d}{dt} \int |\mathcal Y_t|^r w_\tau^r\,dx +\frac{r-1}{2}&\int |\nabla\mathcal Y_t|^2|\mathcal Y_t|^{r-2}w_\tau^r\,dx +\frac12\int |\mathcal Y_t|^r w_\tau^r\,dx +\nu \int |\mathcal Y_t|^{r+2}w_\tau^r\,dx
 \notag\\
 ={}&-\frac12\int\nabla\mathcal Y\cdot\nabla(w_\tau^r)
                         |\mathcal Y|^{r-2}\mathcal Y\,dx -\nu\sum_{k=0}^2
 \big\langle\mathcal X^{(3-k)},\mathcal Y^k
                         |\mathcal Y|^{r-2}\mathcal Yw_\tau^r\big\rangle,
\end{align}
where we used $\nabla( w_\tau^r  |\mathcal{Y}|^{r-2}\mathcal{Y} ) = \nabla(w_\tau^r)  |\mathcal Y|^{r-2}\mathcal Y + w_\tau^r (r-1)|\mathcal Y|^{r-2}\nabla \mathcal Y$.
Since $|\nabla w_\tau^r|\lesssim w_\tau^r$ and $\int w_\tau^r<\infty$ by $\tau r>2$,
the first term on the right is at most
$$\epsilon \int |\nabla\mathcal Y_t|^2|\mathcal Y_t|^{r-2}w_\tau^r\,dx +\epsilon \int |\mathcal Y_t|^{r+2}w_\tau^r\,dx+C_\epsilon$$ 
by Young's inequality.

To estimate the second term, choose a smooth nonnegative partition
of unity $(\psi_z)_{z\in\mathbb Z^2}$ with uniformly bounded derivatives
and supports in the ball $B_z (2) = \{ x\in \mathbb{R}^2; |x-z|\le2\}$.
 Let $\chi_z=1$ on
$B_z(2)$ and $\mbox{supp}(\chi_z) \subset B_z(3)$.
Note that the weight $w_\tau$ is comparable to $w_\tau(z)$ on $B_z(3)$.
For
$$h_z=\psi_z|\mathcal Y|^{r-2}\mathcal Y^{k+1}
              w_\tau^r/w_\tau(z)^r,$$
H\"older's inequality on $B_z(3)$ yields
\begin{align*}
 \|h_z\|_{L^1}&\lesssim \left( \int_{B_z(3) }|\mathcal Y|^{r+2}\,dx  \right)^{\frac{r+k-1}{r+2} },\\
 \|h_z\|_{W^{1,1}}&\simeq \|h_z\|_{L^1} + \|\nabla h_z\|_{L^1}  \\
 &\lesssim \left( \int_{B_z(3)}|\mathcal Y|^{r+2}\,dx \right)^{\frac{r+k-1}{r+2} } + \int_{B_z(3) } |\nabla \mathcal{Y} | |\mathcal{Y}|^{r-2+k} \frac{\psi_z w_\tau^r}{w_\tau(z)^r} dx + \int_{B_z(3) }  |\mathcal{Y}|^{r-1+k} \frac{\nabla (\psi_z w_\tau^r)}{w_\tau(z)^r} dx        \\
 &\lesssim \left( \int_{B_z(3) }|\mathcal Y|^{r+2}\,dx \right)^{\frac{r+k-1}{r+2} }
       +\left( \int_{B_z(3) }|\nabla\mathcal Y|^2|\mathcal Y|^{r-2}\,dx \right)^{1/2} \left( \int_{B_z(3) }|\mathcal Y|^{r+2}\,dx \right)^{(r/2+k-1)/(r+2)}.
\end{align*}
Therefore, duality (Lemma~\ref{lem:duality}), interpolation (Lemma~\ref{lem:interpolation}), and
\eqref{eq:unit-localization} yield 
\begin{align}
\left| \big\langle \psi_z \mathcal X^{(3-k)},\mathcal Y^k |\mathcal Y|^{r-2}\mathcal Yw_\tau^r\big\rangle   \right| &\lesssim w_\tau(z)^r \| \chi_z \mathcal{X}^{(3-k)}  \|_{\mathcal{C}^{-a}} \| h_z \|_{B_{1,1}^a} \notag \\
&\lesssim w_\tau(z)^r \langle z \rangle^\zeta \|  \mathcal{X}^{(3-k)}  \|_{\mathcal{C}^{-a}(w_\zeta)} \| h_z \|_{L^1}^{1-b} \| h_z \|_{W^{1,1}}^b \notag  \\ 
 &\lesssim w_\tau(z)^r\langle z\rangle^\zeta Z_{3-k}
 \left\{\left( \int_{B_z(3) }|\mathcal Y|^{r+2}\,dx \right)^{(r+k-1)/(r+2)} \right. \notag\\
&\qquad \qquad\qquad \left. +\left( \int_{B_z(3) }|\nabla\mathcal Y|^2|\mathcal Y|^{r-2}\,dx  \right)^{b/2}\left( \int_{B_z(3) }|\mathcal Y|^{r+2}\,dx  \right)^{B_{r,k}}\right\},\label{eq:localized-pairing}
\end{align}
where 
\begin{equation}
     B_{r,k}=\frac{(r+k-1)(1-b)+(r/2+k-1)b}{r+2}.
\end{equation}
Noting that
\[
 \frac b2+B_{r,k}=\frac{r+k-1+b}{r+2}<1,
 \qquad
 \left(1-\frac b2-B_{r,k}\right)^{-1}= \frac{r+2}{3-k-b} =q_{r,k} \left(\le \frac{r+2}{1-b}\right),
\]
Young's inequality therefore bounds \eqref{eq:localized-pairing} by
\[
 \epsilon w_\tau(z)^r\left( \int_{B_z(3) }|\nabla\mathcal Y|^2|\mathcal Y|^{r-2}\,dx +  \int_{B_z(3) }|\mathcal Y|^{r+2}\,dx\right)
 +C_\epsilon w_\tau(z)^r
       \big(1+\langle z\rangle^{\zeta q_{r,k}}
                         Z_{3-k}^{q_{r,k}}\big).
\]
Since $r\tau-\zeta q_{r,k}>2$ by \eqref{eq:energy-weight-condition}, we get \eqref{est:Y-proof} by taking sum over $z\in \mathbb{Z}^2$ and choosing $\varepsilon$ sufficiently small.
Finally, the last assertion follows from
\[
 \int |\mathcal Y_t|^r w_\tau^r\,dx \le \left(\int |\mathcal Y_t|^{r+2}w_\tau^r\,dx\right)^{r/(r+2)}
          \left(\int w_\tau^r\,dx\right)^{2/(r+2)}
\]
and $\int w_\tau^r\,dx<\infty$.
\end{proof}

If $Z_j(u)\le u^{-\alpha}M_t$ for $0<u\le t$, application of Lemma~\ref{lem:comparison} on $[t/2,t]$
gives
\begin{equation}\label{eq:energy-cdfi-explicit}
 \|\mathcal Y_t\|_{L^r(w_\tau)}
 \lesssim t^{-1/2}+1+
 \sum_{k=0}^2t^{-\alpha/(3-k-b)}M_t^{1/(3-k-b)}.
\end{equation}
In particular, this is bounded by $t^{-1/2}\mathfrak p(M_t)$ when
$\alpha\le(1-b)/2$. If instead $\mathcal Y_0=0$ and
$Z_j(u)\le u^{-\beta}M_t$ with $\beta (r+2)/(1-b)<1$, integtating 
\eqref{est:Y-proof} over $[0,t]$, one obtains
\begin{equation}\label{eq:zero-initial-energy}
 \sup_{0<u\le t} \|\mathcal Y_t\|_{L^r(w_\tau)}^r 
 \le C\mathfrak p(M_t),\qquad t\le1.
\end{equation}

We also prepare the next estimate, which will be used to upgrade the bounds~\eqref{eq:energy-cdfi-explicit} and~\eqref{eq:zero-initial-energy} on the Lebesgue spaces to the ones on Besov spaces with positive regularity.
\begin{lem}\label{lem:smoothing-heat-weight}
Let $2\le P\le \infty$ $\kappa >2/P$, $0<a<b<s<1$, $\vartheta=b/s$, and $k\in\{1,2\}$. Set $r_k=P(k-\vartheta)/(1-\vartheta)$
If $\tau_k<\kappa$ and
\begin{equation}\label{eq:smoothing-weight-condition}
 \zeta+\vartheta \kappa +(k-\vartheta)\tau_k\le \kappa,
\end{equation}
then, for $0<v\le1$,
\begin{equation}\label{eq:smoothing-product}
 \|e^{\frac{v}2(\Delta-1)} (XY^k)\|_{\mathcal{C}^s(w_\kappa)}
 \lesssim v^{-(s+a)/2}\|X\|_{\mathcal C^{-a}(w_\zeta)}
 \|Y\|_{L^{r_k}(w_{\tau_k})}^{k-\vartheta}
 \|Y\|_{\mathcal{C}^s(w_\kappa)}^\vartheta.
\end{equation}
Moreover, if $\kappa-\zeta>2/P$,
\begin{align}
 \|e^{\frac{v}2(\Delta-1)}  X\|_{B^s_{P,\infty}(w_\kappa)}
 &\lesssim v^{-(s+a)/2}\|X\|_{\mathcal C^{-a}(w_\zeta)},
 \label{eq:smoothing-forcing}\\
 \|e^{\frac{v}2(\Delta-1)}  (Y^3)\|_{B^s_{P,\infty}(w_\kappa)}
 &\lesssim v^{-s/2}\|Y\|_{L^{3P}(w_{\kappa/3})}^3.
 \label{eq:smoothing-cubic}
\end{align}
\end{lem}
\begin{proof}
Choose $q$ by $1/q=\vartheta/P+(1-\vartheta)/r_k$ and put $\omega=\vartheta \kappa+(1-\vartheta)\tau_k$.
Then $1/P=1/q+(k-1)/\tau_k$ and interpolation (Lemma~\ref{lem:interpolation}) gives
\[
 \|Y\|_{B_{q,\infty}^b(w_\omega)}
 \lesssim\|Y\|_{B^s_{P,\infty}(w_\kappa)}^\vartheta
             \|Y\|_{L^{r_k}(w_{\tau_k})}^{1-\vartheta}.
\]
Therefore, applying the product estimates \eqref{eq:power-product}, \eqref{eq:rough-product}, and the Schauder
estimate (Lemma~\ref{lem:schauder}) and using \eqref{eq:smoothing-weight-condition}, we obtain
\begin{align}
\|e^{\frac{v}2(\Delta-1)} (XY^k)\|_{B^s_{P,\infty}(w_\kappa)} &\lesssim v^{-(s+a)/2}\|X Y^k\|_{B^{-a}_{P,\infty}(w_\kappa)}\notag\\
&\lesssim v^{-(s+a)/2}\|X\|_{\mathcal C^{-a}(w_\zeta)} 
 \|Y^k\|_{B^b_{P,\infty}(w_{\omega+(k-1)\tau_k})} \notag\\
 &\lesssim v^{-(s+a)/2}\|X\|_{\mathcal C^{-a}(w_\zeta)} 
 \|Y\|_{B^b_{q,\infty}(w_{\omega})} \|Y\|_{L^{r_k}(w_{\tau_k})}^{k-1} \notag\\
&\lesssim v^{-(s+a)/2}\|X\|_{\mathcal C^{-a}(w_\zeta)} 
 \|Y\|_{L^{r_k}(w_{\tau_\kappa})}^{k-\vartheta}
 \|Y\|_{B^s_{P,\infty}(w_\kappa)}^\vartheta.    
\end{align}
For \eqref{eq:smoothing-forcing}, use the Schauder estimate (Lemma~\ref{lem:schauder}) and $\mathcal{C}^s(w_\zeta) \hookrightarrow B^s_{P,\infty}(w_\kappa)$.
For \eqref{eq:smoothing-cubic}, use
$\|Y^3\|_{L^P(w_\kappa)}=\|Y\|_{L^{3P}(w_{\kappa/3})}^3$ and
$L^P\hookrightarrow B^0_{P,\infty}$.
\end{proof}

\begin{proof}[Proof of Lemma~\ref{lem:CDFI-Y}]
We first prove the estimate \eqref{est:YY2}, which is independent of the initial condition.
Fix any $2\le P< \infty$ and $\kappa>2/P$.
Let $s\in(\eta,\min\{1,\eta+\varepsilon\})$, $0<2a<s$, $\vartheta=2a/s$, $\zeta<(1-\vartheta)\kappa/2$, and, for $k=1,2$, set
$$  r_k = \frac{P(k-\vartheta)}{1-\theta}    
 $$
 and choose $\tau_k>0$ sufficiently small such that $\vartheta\kappa + (k-\vartheta)\tau_k < \kappa $. Fix $t\le1$, let $t_0=t/2$, and set
\[
 H=\sup_{t_0<u\le t}(u-t_0)^{s/2}
                  \|\mathcal Y_u\|_{B_{P,\infty}^s(w_\kappa)}.
\]
The mild equation on $[t_0,t]$ and Lemma~\ref{lem:smoothing-heat-weight} give for some $\zeta>0$ that
\begin{align}
\|\mathcal Y_u\|_{B^s_{P,\infty}(w_\kappa)} &\lesssim \left\| e^{\frac{u- t_0}2(\Delta-1)} \mathcal{Y}_{t_0}  \right\|_{B^s_{P,\infty}(w_\kappa)} + \sum_{k=1}^2 \int^u_{t_0} \left\| e^{\frac{u-v}2(\Delta-1)} \left( \mathcal{X}_v^{(3-k)} \mathcal{Y}_v^k   \right)   \right\|_{B^s_{P,\infty}(w_\kappa)} dv  \notag\\
&\ \ \ \ +\int^u_{t_0} \left\| e^{\frac{u-v}2(\Delta-1)}  \mathcal{X}_v^{(3)}    \right\|_{B^s_{P,\infty}(w_\kappa)} dv  + \int^u_{t_0} \left\| e^{\frac{u-v}2(\Delta-1)} \left(  \mathcal{Y}_v^3   \right)   \right\|_{B^s_{P,\infty}(w_\kappa)} dv \notag\\
&\lesssim (u-t_0)^{-\frac{s}2} \| \mathcal{Y}_{t_0} \|_{L^P(w_\kappa)} + \sum_{k=1}^2 \int^u_{t_0} (u-v)^{-(s+a)/2}\| \mathcal{X}_v^{(k)} \|_{\mathcal C^{-a}(w_\zeta)}
 \| \mathcal{Y}_v \|_{L^{r_k}(w_{\tau_k})}^{k-\vartheta}
 \| \mathcal{Y}_v \|_{B^s_{P,\infty}(w_\kappa)}^\vartheta dv  \notag\\
&\ \ \ \ +\int^u_{t_0} (u-v)^{-(s+a)/2} \left\|  \mathcal{X}_v^{(3)}    \right\|_{\mathcal{C}^{-a}(w_\zeta)} dv  + \int^u_{t_0}  (u-v)^{-s/2} \left\|   \mathcal{Y}_v \right\|_{L^{3P}(w_{\kappa/3})}^3 dv.
\end{align}
Therefore, \eqref{eq:energy-cdfi-explicit}  yields for $\kappa_X\le \zeta$
\begin{align}\label{eq:A-H-inequality}
 H&\lesssim t^{-1/2}\mathfrak p(\mathfrak X_t) +  \sum_{k=1}^2 \sup_{0<h\le t-t_0} h^{s/2} \int^h_{0} (h-v)^{-(s+a)/2} (v+t_0)^{-\alpha -(k-\vartheta)2} v^{-a}    dv\, p(\mathfrak X_t) H^\vartheta   \notag\\
 &\ + \sup_{0<h\le t-t_0} h^{s/2} \int^h_{0} \left\{ (h-v)^{-(s+a)/2} (v+t_0)^{-\alpha } \mathfrak X_t  + (h-v)^{-s/2} (v+t_0)^{-3/2} \mathfrak p(\mathfrak X_t)  \right\} dv     \notag\\
 &\lesssim{}t^{-1/2}\mathfrak p(\mathfrak X_t)
       +t^{1-a/2-\alpha}\mathfrak X_t  +\sum_{k=1}^2
 t^{1-3a/2-\alpha-(k-\vartheta)/2}
             \mathfrak p(\mathfrak X_t)H^\vartheta,
\end{align}
where we used $\vartheta s=2a$ and
\[
 h^{s/2}\int_0^h(h-v)^{-(s+a)/2}v^{-a}\,dv
   =C h^{1-3a/2},\qquad 0<h\le t/2.
\]
We choose $\alpha>0$ so small that
$\alpha+3a/2<1/2$. Then the exponent in the third line of
\eqref{eq:A-H-inequality} is at least $-(1-\vartheta)/2$;
$$ -\frac{1-\vartheta}2 < 1 - \frac{3}2a - \alpha - \frac{k-\vartheta}2.  $$
Young's inequality thus absorbs $H^\vartheta$ and gives
\[
 H\lesssim t^{-1/2}\mathfrak p(\mathfrak X_t),\qquad
 \|\mathcal Y_t\|_{B^s_{P,\infty}(w_\kappa)}
       \lesssim t^{-(1+s)/2}\mathfrak p(\mathfrak X_t).
\]
Together with $s<\eta+\varepsilon$ and the Besov embedding,
this proves \eqref{est:YY2} by choosing sufficiently large $P$.

We now turn to the proof of \eqref{est:Y3}.
Write $I_t= e^{\frac{t}2(\Delta-1)} f$ and $V_t=\mathcal Y_t-I_t$.
Then $V_0=0$ and $V$ solves \eqref{eq:YY} with coefficients
\begin{align}\label{eq:shifted-coefficients}
 \widetilde{\mathcal X}^{(1)}&=\mathcal X^{(1)}+3I,\notag\\
 \widetilde{\mathcal X}^{(2)}&=\mathcal X^{(2)}
                              +2I\mathcal X^{(1)}+3I^2,\\
 \widetilde{\mathcal X}^{(3)}&=\mathcal X^{(3)}
                         +I\mathcal X^{(2)}+I^2\mathcal X^{(1)}+I^3.
 \notag
\end{align}
Let all parameters are as above.
By the Schauder and product estimates (Lemmas~\ref{lem:schauder} and~\ref{lem:product}),
\[
 \|I_u\|_{\mathcal C^{2a}(w_{\kappa_f})}\lesssim u^{-3a/2}F,
 \qquad
 \max_j\|\widetilde{\mathcal X}^{(j)}_u\|_{
                      \mathcal C^{-a}(w_{\zeta})}
 \lesssim u^{-\beta}\mathfrak p(F,\mathfrak X_t),\quad 0<u\le t.
\]
for $\beta=\alpha + 9a/2$ and sufficiently small $\kappa_X, \kappa_f$ such that
$\kappa_X+3\kappa_f<\zeta$.
Then \eqref{eq:zero-initial-energy} yields that for $(r, \tau)= (r_k,\tau_k), (3P,\kappa/3) $, 
\begin{align}
    \| V_t \|_{L^r(w_\tau)} \lesssim t^{-\beta} \mathfrak{p}(F,\mathfrak X_t )
\end{align}
for some polynomial $\mathfrak{p}$, if $\beta, \kappa_f,\kappa_X$ are sufficiently small.
Applying the mild equation from time zero,
\begin{align}
\|V_u\|_{B^s_{P,\infty}(w_\kappa)} &\lesssim \sum_{k=1}^2 \int^u_{0} \left\| e^{\frac{u-v}2(\Delta-1)} \left( \widetilde{\mathcal{X}}_v^{(3-k)} V_v^k   \right)   \right\|_{B^s_{P,\infty}(w_\kappa)} dv  \notag\\
&\ \ \ \ +\int^u_{0} \left\| e^{\frac{u-v}2(\Delta-1)}  \widetilde{\mathcal{X}}_v^{(3)}    \right\|_{B^s_{P,\infty}(w_\kappa)} dv  + \int^u_{0} \left\| e^{\frac{u-v}2(\Delta-1)} \left(  V_v^3   \right)   \right\|_{B^s_{P,\infty}(w_\kappa)} dv .
\end{align}
Therefore, similarly to above, Lemma~\ref{lem:smoothing-heat-weight} yields for 
$H_0=\sup_{0<u\le t}\|V_u\|_{B^s_{P,\infty}(w_\kappa)}$ that
\begin{align}
 H_0 &\lesssim  \sup_{0<u \le t}  \int^u_{0} (u-v)^{-(s+a)/2}  v^{-\beta}    dv\, p(F, \mathfrak X_t) H_0^\vartheta   \notag\\
 &\ + \sup_{0<u\le t}  \int^u_{0} \left\{ (u-v)^{-(s+a)/2} v^{-\beta}  + (u-v)^{-s/2}  \right\} dv \, p(F, \mathfrak X_t)      \notag \\
 &\lesssim \mathfrak p(F,\mathfrak X_t)(1+H_0^\vartheta),
\end{align}
if $\beta$ is so small that
\[
 \int_0^u(u-v)^{-(s+a)/2}v^{-\beta} \,dv
       =C u^{1-(s+a)/2-\beta}\le C.
\]
  Absorption by Young's inequality again gives a polynomial bound for
$H_0$. Finally, if $\kappa>\kappa_f$, the Schauder estimate implies
\[
 \|I_t\|_{B^s_{P,\infty}(w_\kappa)}\lesssim t^{-(s+a)/2}F.
\]
Combining these estimates with Besov embedding, we get \eqref{est:Y3} by choosing sufficiently large $P$.

\end{proof}

\section{Stochastic estimates}
In this section, we collect stochastic estimates used in the paper. 

\subsection{Estimates on the Wick powers and renormalization constants}
Recall that $W$ is a cylindrical Brownian motion on $L^2(\mathbb{R}^2)$ and $\Pi_N$ is the mollification operator introduced in Section~\ref{sec:pre}.
We define
\begin{align}
    Z_t = \int^t_0 e^{\frac{t-s}2(\Delta -1)}\, dW_s,\ \ \  Z_{t,\infty} = \int^t_{-\infty} e^{\frac{t-s}2(\Delta -1)}\, dW_s \label{eq:def-X}
\end{align}
and
\begin{align}
    C_{N,t} \coloneqq \mathbb{E} \left[ |\Pi_N Z_t(x)|^2  \right], \ \ \ C_{N} \coloneqq \mathbb{E} \left[ |\Pi_N Z_{t,\infty}(x)|^2  \right].
\end{align}
We will use the following identity repeatedly.
\begin{lem}\label{lem:covariance}
    We have
    \begin{align}
        \mathbb{E} \left[ \Pi_N Z_t(x) \Pi_N Z_s(y)  \right] = \frac{1}{(2\pi)^2} \int_{\mathbb{R}^2}  \left( e^{-\frac{|t-s|}{2} (|\xi|^2 +1)}  - e^{-\frac{t+s}{2} (|\xi|^2 +1)} \right)  \frac{1}{|\xi|^2+1}  |g_N(\xi)|^2 e^{i(x-y)\cdot \xi} \, d\xi.
    \end{align}
\end{lem}

\begin{proof}
The definition \eqref{eq:def-X} of $X$ and the It\^o isometry yield
\begin{align}
    &\mathbb{E} \left[ \Pi_N Z_t(x) \Pi_N Z_s(y)  \right] = \mathbb{E} \left[ \langle Z_t , \Pi_N \delta_x \rangle \langle  Z_s, \Pi_N \delta_y \rangle  \right] = \left\langle  \int^{t\wedge s}_0 e^{ (\frac{t+s}{2} - u )(\Delta -1)} \Pi_N \delta_x \, du, \Pi_N \delta_y     \right\rangle.
\end{align}
Noting that 
\begin{align}
    \Pi_N\delta_x(z) = \frac{1}{2\pi} \int_{\mathbb{R}^2} g_N(\xi) e^{i(z-x)\cdot \xi}\, d\xi,
\end{align}
we also obtain, by Plancherel's theorem, 
\begin{align}
 &\left\langle  \int^{t\wedge s}_0 e^{ (\frac{t+s}{2} - u )(\Delta -1)} \Pi_N \delta_x \, du, \Pi_N \delta_y     \right\rangle \notag\\
 &= \frac{1}{(2\pi)^2} \int_{\mathbb{R}^2}  \left(\int^{t\wedge s}_0 e^{- (\frac{t+s}{2} - u )(|\xi|^2 +1)} \, du \right) |g_N(\xi)|^2 e^{i(x-y)\cdot \xi} \, d\xi   \notag\\
 &=\frac{1}{(2\pi)^2} \int_{\mathbb{R}^2}  \left( e^{-\frac{|t-s|}{2} (|\xi|^2 +1)}  - e^{-\frac{t+s}{2} (|\xi|^2 +1)} \right)  \frac{1}{|\xi|^2+1}  |g_N(\xi)|^2 e^{i(x-y)\cdot \xi} \, d\xi. 
\end{align}

\end{proof}

Lemma~\ref{lem:covariance} yields
\begin{align}
    C_{N,t} &= \mathbb{E} \left[ \langle X_t, \Pi_N \delta_x \rangle^2  \right]  
    = \frac{1}{(2\pi)^2}\int_{\mathbb{R}^2}  \left|g\left(\frac{\xi}N  \right) \right|^2 \frac{1}{|\xi|^2+1} \left\{ 1-  e^{-t(|\xi|^2+1)}   \right\}    \, d\xi,
\end{align}
\begin{align}
    C_{N} &=  \frac{1}{(2\pi)^2} \int_{\mathbb{R}^2}  \left|g\left(\frac{\xi}N  \right) \right|^2 \frac{1}{|\xi|^2+1}   \, d\xi  = O(\log N ), \label{eq:renor-const}
\end{align}
and thus
\begin{align}
  0\le  C_{N} - C_{N,t} &=\frac{1}{(2\pi)^2}   \int_{\mathbb{R}^2}  \left|g\left(\frac{\xi}N  \right) \right|^2 \frac{1}{|\xi|^2+1}e^{-t(|\xi|^2+1)}   \, d\xi \lesssim_\varepsilon t^{-\varepsilon} \label{eq:const-differ}
\end{align}
uniformly in $N\in \mathbb{N}$ for any $\varepsilon>0$.

We also have
\begin{align}
    \Wick{(\Pi_N Z_t)^2}_{C_N} = (\Pi_N Z_t )^2 - C_N = \Wick{(\Pi_N Z_t)^2}_{C_{N,t}} + (C_{N,t} - C_N) \label{eq:const-differ-2}
\end{align}
and 
\begin{align}
    \Wick{(\Pi_N Z_t)^3}_{C_N} = (\Pi_N Z_t )^3 - 3C_N \Pi_N Z_t = \Wick{(\Pi_N Z_t)^3}_{C_{N,t}} + 3(C_{N,t} - C_N)\Pi_N Z_t. \label{eq:const-differ-3}
\end{align}

We record the following elementary integrability criterion.
\begin{lem}\label{lem:sum-conv}
    If $0< \alpha_0,\alpha_1, \cdots, \alpha_k \le d$ and $\alpha_0+\alpha_1+\cdots+\alpha_k>dk$, we have
    \begin{align}
        \int_{\mathbb{R}^{dk}} d\xi_1\cdots d\xi_k \frac{1}{|\xi_1+ \cdots + \xi_k |^{\alpha_0}+1} \prod_{m=1}^k \frac{1}{|\xi_m|^{\alpha_m}+1}< \infty.
    \end{align}
\end{lem}
\begin{proof}
    Writing $\rho_{\alpha}(\xi) \coloneqq \frac{1}{|\xi|^\alpha+1}$, we have
    \begin{align}
     &\int_{\mathbb{R}^{dk}} d\xi_1\cdots d\xi_k\, \frac{1}{|\xi_1+ \cdots + \xi_k |^{\alpha_0}+1} \prod_{m=1}^k \frac{1}{|\xi_m|^{\alpha_m}+1}\notag\\
     &= \int_{\mathbb{R}^{d(k-1)}} d\xi_2\cdots d\xi_k\, \prod_{m=2}^k \rho_{\alpha_m}(\xi_m) \int_{\mathbb{R}^d} d\xi_1\, \rho_{\alpha_0}(\xi_1+\cdots+\xi_k) \rho_{\alpha_1}(\xi_1)\notag\\
    &=\int_{\mathbb{R}^{d(k-1)}} d\xi_2\cdots d\xi_k\, \prod_{m=2}^k \rho_{\alpha_m}(\xi_m) \int_{\mathbb{R}^d} d\xi_1\, \rho_{\alpha_0}(-\xi_1+\cdots+\xi_k) \rho_{\alpha_1}(\xi_1) \notag\\
    &=\int_{\mathbb{R}^{d(k-1)}} d\xi_2\cdots d\xi_k\, (\rho_{\alpha_0}*\rho_{\alpha_1}) (\xi_2+ \cdots + \xi_k) \prod_{m=2}^k \rho_{\alpha_m}(\xi_m) .
    \end{align}
    By iteration, we obtain
    \begin{align}
     &\int_{\mathbb{R}^{dk}} d\xi_1\cdots d\xi_k \,\frac{1}{|\xi_1+ \cdots + \xi_k |^{\alpha_0}+1} \prod_{m=1}^k \frac{1}{|\xi_m|^{\alpha_m}+1}\notag\\
     &= \int_{\mathbb{R}^d} d\xi_k\, (\rho_{\alpha_0}* \rho_{\alpha_1}*\cdots*\rho_{\alpha_{k-1}})(\xi_k) \rho_{\alpha_k}(\xi_k) .  
    \end{align}
Applying H\"older's inequality and Young's convolution inequality, we get
\begin{align}
 &\int_{\mathbb{R}^{dk}} d\xi_1\cdots d\xi_k\,
   \frac{1}{|\xi_1+\cdots+\xi_k|^{\alpha_0}+1}
   \prod_{m=1}^k\frac{1}{|\xi_m|^{\alpha_m}+1}\notag\\*
 &\qquad\le \|\rho_{\alpha_0}\|_{L^{\beta_0}(\mathbb{R}^d)}
   \times\cdots\times\|\rho_{\alpha_k}\|_{L^{\beta_k}(\mathbb{R}^d)}
\end{align}
for any $1\le \beta_0,\cdots, \beta_k \le \infty$ such that $k = 1/\beta_0+ \cdots + 1/\beta_k$.
Choosing $\beta_m = (\alpha_0+ \cdots + \alpha_k)/(k\alpha_m) $, we obtain the desired result.

\end{proof}

\begin{lem}\label{lem:est-wick}
    Let $k=1,2,3$, $\alpha<0$, $1\le p<\infty$ and $\kappa,\delta>0$ satisfy $p\kappa>1$, $\delta k<1 \wedge(-\alpha) $ and $p\delta>1$. Then
    \begin{align}
      \sup_{N\ge 1}  \mathbb{E} \left[ \sup_{t \in[0,T]} \left\|  \Wick{(\Pi_N Z_t)^k}_{C_{N,t}} \right\|^{2p}_{B_{2p,2p}^\alpha(w_\kappa)}    \right] \lesssim T^{p\delta} .
    \end{align}
\end{lem}

\begin{proof}
In this proof, we suppress the dependence on $C_{N,t}$ and simply write $\Wick{(\Pi_N Z_t)^k}= \Wick{(\Pi_N Z_t)^k}_{C_{N,t}}$.
By the hypercontractivity of Gaussian polynomials (cf. \cite{Jan97}), we have
\begin{align}
&\mathbb{E} \left[  \left\|  \Wick{(\Pi_N Z_t)^k} - \Wick{(\Pi_N Z_s)^k}  \right\|^{2p}_{B_{2p,2p}^\alpha(w_\kappa)}    \right]\notag\\
&= \mathbb{E} \left[ \sum_{j\ge -1} 2^{2pj\alpha} \left\| \Delta_j\left\{ \Wick{(\Pi_N Z_t)^k} - \Wick{(\Pi_N Z_s)^k} \right\} \right\|^{2p}_{L^{2p}(w_\kappa)}    \right]  \notag\\
&= \sum_{j\ge -1} 2^{2pj\alpha} \int_{\mathbb{R}^2}  \mathbb{E} \left[ \left| \Delta_j\left\{ \Wick{(\Pi_N Z_t)^k} - \Wick{(\Pi_N Z_s)^k} \right\}(x) w_\kappa(x)    \right|^{2p}   \right]       \, dx \notag\\
&\lesssim \sum_{j\ge -1} 2^{2pj\alpha} \int_{\mathbb{R}^2} w_\kappa^{2p}(x) \mathbb{E} \left[ \left| \Delta_j\left\{ \Wick{(\Pi_N Z_t)^k} - \Wick{(\Pi_N Z_s)^k} \right\}(x)     \right|^{2}   \right]^p      \, dx. \label{eq:wick-0}
\end{align}
Furthermore, we have
\begin{align}
 &\mathbb{E} \left[ \left| \Delta_j\left\{ \Wick{(\Pi_N Z_t)^k} - \Wick{(\Pi_N Z_s)^k} \right\}(x)     \right|^{2}   \right]\notag\\
 &= \Delta_{j,x} \Delta_{j,y} \mathbb{E} \left[  \left\{ \Wick{(\Pi_N Z_t)^k} - \Wick{(\Pi_N Z_s)^k} \right\}(x) \left\{ \Wick{(\Pi_N Z_t)^k} - \Wick{(\Pi_N Z_s)^k} \right\}(y)     \right]  \Big|_{y=x} \label{eq:wick-1}
\end{align}
and the standard formula for Wick products (cf. \cite[Theorem 3.9]{Jan97}) yields
\begin{align}
     \mathbb{E} \left[ \Wick{(\Pi_N Z_t)^k}(x) \Wick{(\Pi_N Z_s)^k}(y)  \right] &= k!\mathbb{E} \left[ \Pi_N Z_t(x) \Pi_N Z_s(y)  \right]^k. \label{eq:wick-2} 
\end{align}
Lemma~\ref{lem:covariance} also gives
\begin{align}
 &\Delta_{j,x}\Delta_{j,y}
   \mathbb{E}\left[\Pi_N Z_t(x)\Pi_N Z_s(y)\right]^k\notag\\  
 &=\frac{1}{(2\pi)^{2k}}\int_{\mathbb{R}^{2 k}}d\xi_1\cdots d\xi_k
   |\rho_j(\xi_1+\cdots+\xi_k)|^2\notag\\
 &\qquad\prod_{m=1}^k\left(
       e^{-\frac{|t-s|}{2}(|\xi_m|^2+1)}
       -e^{-\frac{t+s}{2}(|\xi_m|^2+1)}\right)
   \frac{|g_N(\xi_m)|^2}{|\xi_m|^2+1}
   e^{i(x-y)\cdot\xi_m}.\label{eq:wick-3}
\end{align}
From \eqref{eq:wick-1}, \eqref{eq:wick-2} and \eqref{eq:wick-3}, we obtain
\begin{align}
 &\mathbb{E} \left[ \left| \Delta_j\left\{ \Wick{(\Pi_N Z_t)^k} - \Wick{(\Pi_N Z_s)^k} \right\}(x)     \right|^{2}   \right]\notag\\
 &= \frac{k!}{(2\pi)^{2k}} \int_{\mathbb{R}^{2 k}} d\xi_1\cdots d\xi_k |\rho_j(\xi_1 + \cdots + \xi_k)|^2 A(t,s;\xi_1,\cdots,\xi_k) \prod_{m=1}^k   \frac{1}{|\xi_m|^2+1}  |g_N(\xi_m)|^2, \label{eq:wick-4}  
\end{align}
where
\begin{align}
 A(t,s;\xi_1,\cdots,\xi_k)
 &\coloneqq\prod_{m=1}^k(1-e^{-t(|\xi_m|^2+1)})
       +\prod_{m=1}^k(1-e^{-s(|\xi_m|^2+1)})\notag\\*
 &\qquad-2\prod_{m=1}^k\left(
       e^{-\frac{|t-s|}{2}(|\xi_m|^2+1)}
       -e^{-\frac{t+s}{2}(|\xi_m|^2+1)}\right).
\end{align}
Since 
\begin{align}
 A(t,s;\xi_1,\cdots,\xi_k) \lesssim |t-s|^\delta (|\xi_1|^{2\delta}+ \cdots +|\xi_k|^{2\delta} + 1 ) 
\end{align}
and
\begin{align}
|\rho_j(\xi_1 + \cdots + \xi_k)|^2 \lesssim 2^{2j\varepsilon} \frac{1}{|\xi_1 + \cdots \xi_k  |^{2\varepsilon}+1} ,  
\end{align}
it follows from Lemma~\ref{lem:sum-conv}, \eqref{eq:wick-0}, and \eqref{eq:wick-4} that 
\begin{align}
 &\mathbb{E}\left[\left\|
     \Wick{(\Pi_N Z_t)^k}-\Wick{(\Pi_N Z_s)^k}
     \right\|^{2p}_{B_{2p,2p}^\alpha(w_\kappa)}\right]\notag\\
 &\lesssim |t-s|^{p\delta}\sum_{j\ge-1}2^{2pj(\alpha+\varepsilon)} \left(\int_{\mathbb{R}^{2 k}}d\xi_1\cdots d\xi_k
     \frac{(|\xi_1|^{2\delta}+\cdots+|\xi_k|^{2\delta}+1)}
          {|\xi_1+\cdots+\xi_k|^{2\varepsilon}+1}
     \prod_{m=1}^k\frac{1}{|\xi_m|^2+1}\right)^p
     \int_{\mathbb{R}^2}w_\kappa^{2p}(x)\,dx\notag\\
 &\lesssim |t-s|^{p\delta}
\end{align}
if $\alpha+\varepsilon<0, \delta k <1 \wedge \varepsilon$, and $2p\kappa>2$.
The Kolmogorov continuity theorem (cf. \cite[Theorem 4.3.2]{Str11}) therefore yields the bound
\begin{align}
\mathbb{E} \left[ \sup_{0\le s < t \le T} \frac{\left\|  \Wick{(\Pi_N Z_t)^k} - \Wick{(\Pi_N Z_s)^k}  \right\|^{2p}_{B_{2p,2p}^\alpha(w_\kappa)} }{(t-s)^{2p\beta}}   \right] \lesssim T^{p\delta- 2p\beta}    
\end{align}
for any $\beta, p$ such that $0< \beta < (p\delta-1)/2p\le1$. Combining this bound with $\Wick{(\Pi_N Z_0)^k}=0$, we obtain
\begin{equation}
 \mathbb{E} \left[ \sup_{0<t\le T} \left\|  \Wick{(\Pi_N Z_t)^k}   \right\|^{2p}_{B_{2p,2p}^\alpha(w_\kappa)}    \right]   \lesssim T^{p\delta}.\qedhere
\end{equation}

\end{proof}

\begin{cor}\label{cor:wick}
    Let $k=1,2,3$, $\gamma,\eta,\kappa>0$ and $1\le p<\infty$. The sequence $\Wick{(\Pi_N Z)^k}_{C_N}$ converges to a limit $\Wick{Z^k}$ in the norm
    \begin{align}
        \|  F\|_{\gamma,T;\eta,\kappa} \coloneqq \sup_{0<t\le T}t^\gamma \| F_t \|_{\mathcal{C}^{-\eta}(w_\kappa) } \label{eq:cor-wick-norm}
    \end{align}
in $L^p(\Omega)$. Moreover, we have the bound
\begin{align}
 \sup_{N\ge 1} \mathbb{E}\left[ \| \Wick{(\Pi_N Z)^k}_{C_N} \|_{\gamma,T;\eta,\kappa}^p  \right]   \lesssim T^{p\sigma} \label{eq:wick-cor}
\end{align}
for any $0< \sigma < \gamma$ and $0<T\le 1$.
\end{cor}

\begin{proof}
In view of the embedding $B_{2p,2p}^{-\eta_0}(w_\kappa) \hookrightarrow \mathcal{C}^{-\eta}(w_\kappa) $ for $2/(2p)< \eta - \eta_0$, it follows by taking $p\ge1$ sufficiently large in Lemma~\ref{lem:est-wick} that
\begin{align}
    \sup_{N\ge 1} \mathbb{E} \left[ \| \Wick{(\Pi_N Z)^k}_{C_{N,t}}\|_{\mathcal{C}^{-\eta}(w_\kappa) }^p  \right] \lesssim 1
\end{align}
for any $0< T\le 1$. Together with \eqref{eq:const-differ}--\eqref{eq:const-differ-3}, this implies \eqref{eq:wick-cor}. The convergence assertion follows from an analogous estimate on the difference $\Wick{(\Pi_N Z)^k}_{C_N} - \Wick{(\Pi_M Z)^k}_{C_M}$, which can be obtained by a straightforward modification of the proof of Lemma~\ref{lem:est-wick}.
    
\end{proof}

\subsection{Wick powers under the $\Phi^4_2$-measure}\label{sec:wick-phi42}
We prove that the Wick powers converge in a suitable topology under the $\Phi^4_2$-measure $\mu^{(\lambda)}$. Our proof uses stochastic quantization and the bounds established in Section~\ref{sec:CDFI}. Corollary~\ref{cor:wick} states that $\Wick{(\Pi_N Z_t)^k}_{C_N}$ converges to a limit $\Wick{Z_t^k}$ in the norm defined in \eqref{eq:cor-wick-norm}. Let $Y_t=Y_t^\phi$ be the solution to
\begin{align}
    \begin{cases}
        \partial_t Y_{t}= \frac{1}2(\Delta-1)Y_{t}-\frac{\lambda}2 \sum_{k=0}^2      \Wick{Z_{t}^{3-k} }  Y_{t}^k - \lambda Y_t^3 \\
        Y_{0}= \phi.
    \end{cases}
\end{align}
The solution to the stochastic quantization equation \eqref{eq:sq} is then given by $\Phi_t(\phi) \coloneqq Z_t + Y_t^\phi$, and the $\Phi^4_2$-measure $\mu^{(\lambda)}$ is invariant under $\Phi_t$; see, for example, \cite{RZZ17}. Combining this invariance with the results of Section~\ref{sec:CDFI}, we obtain the following.

\begin{lem}\label{lem:phi42-wick}
    Let $1\le p<\infty$ and $\eta,\kappa>0$. We have
    \begin{align}
        \int \left\|  \Wick{(\Pi_N \phi )^k}_{C_N} - \Wick{(\Pi_M \phi)^k}_{C_M}    \right\|^p_{\mathcal{C}^{-\eta}(w_\kappa)} \mu^{(\lambda)}(d\phi) \to 0
    \end{align}
    as $M,N \rightarrow \infty$. In particular, there exists $\Wick{\phi^k}\in L^p(\mu^{(\lambda)};\mathcal{C}^{-\eta}(w_\kappa) )$ such that 
    \begin{align}
        \int \left\|  \Wick{(\Pi_N \phi )^k}_{C_N} - \Wick{\phi^k}    \right\|^p_{\mathcal{C}^{-\eta}(w_\kappa)} \mu^{(\lambda)}(d\phi) \to 0.
    \end{align}
\end{lem}

\begin{proof}
As noted above, the $\Phi^4_2$-measure $\mu^{(\lambda)}$ is invariant under the dynamical $\Phi^4_2$-model, so $Z_1 + Y_1^\phi$ is distributed according to $\mu^{(\lambda)}$ when $\phi$ is distributed according to $\mu^{(\lambda)}$ and independent of the noise. For simplicity, we write $Z= Z_1$ and $Y=Y_1$ below.
By the coming-down-from-infinity estimate in Lemma~\ref{lem:CDFI-Y} and the stochastic estimate in Corollary~\ref{cor:wick}, we have, for any $0<\rho<1, \xi>0, 1\le q<\infty$,
\begin{align}
  \int \mathbb{E} \left[ \| Y^\phi \|^q_{\mathcal{C}^\rho(w_\xi)}    \right]     \mu^{(\lambda)}(d\phi) < \infty.  \label{eq:bound-ex-Y}
\end{align}
On the other hand, noting that
\begin{align}
    \Wick{\left\{ \Pi_N(Z + Y^\phi)_{C_N}   \right\}^k} = \sum_{l=0}^k \binom{k}l (\Pi_N Y^\phi )^l \Wick{(\Pi_N Z )^{k-l}}_{C_N},
\end{align}
we have
\begin{align}
 &\int\left\|\Wick{(\Pi_N\phi)^k}_{C_N}
       -\Wick{(\Pi_M\phi)^k}_{C_M}\right\|^p_{\mathcal{C}^{-\eta}(w_\kappa)}
       \mu^{(\lambda)}(d\phi)\notag\\*
 &=\int\mathbb{E}\left[\left\|\sum_{l=0}^k\binom{k}l
       \left\{\begin{aligned}
        &(\Pi_NY^\phi)^l\Wick{(\Pi_NZ)^{k-l}}_{C_N}\\
        &\quad-(\Pi_MY^\phi)^l\Wick{(\Pi_MZ)^{k-l}}_{C_M}
       \end{aligned}\right\}
       \right\|^p_{\mathcal{C}^{-\eta}(w_\kappa)}\right]
       \mu^{(\lambda)}(d\phi)\notag\\*
 &\lesssim\sum_{l=0}^k\int\mathbb{E}\left[\begin{aligned}
       &\left\|\left\{(\Pi_NY^\phi)^l-(\Pi_MY^\phi)^l\right\}
       \Wick{(\Pi_NZ)^{k-l}}_{C_N}\right\|^p_{\mathcal{C}^{-\eta}(w_\kappa)}\\
       &\quad+\left\|(\Pi_MY^\phi)^l\left\{
       \Wick{(\Pi_NZ)^{k-l}}_{C_N}-\Wick{(\Pi_MZ)^{k-l}}_{C_M}
       \right\}\right\|^p_{\mathcal{C}^{-\eta}(w_\kappa)}
       \end{aligned}\right]\,\mu^{(\lambda)}(d\phi).
\end{align}
Moreover, the product estimate (Lemma~\ref{lem:product}) yields
\begin{align}
 &\left\|\left\{(\Pi_NY^\phi)^l-(\Pi_MY^\phi)^l\right\}
       \Wick{(\Pi_NZ)^{k-l}}\right\|_{\mathcal{C}^{-\eta}(w_\kappa)}\notag\\
 &\lesssim\left\|(\Pi_NY^\phi)^l-(\Pi_MY^\phi)^l
       \right\|_{\mathcal{C}^{\eta}(w_{\kappa/2})}
       \left\|\Wick{(\Pi_NZ)^{k-l}}
       \right\|_{\mathcal{C}^{-\eta/2}(w_{\kappa/2})}\notag\\
 &\lesssim\left\|\Pi_NY^\phi-\Pi_MY^\phi
       \right\|_{\mathcal{C}^{\eta}(w_{\kappa/(2l)})}\notag\\*
 &\qquad \times \left\{\left\|\Pi_NY^\phi\right\|^{l-1}_{\mathcal{C}^{\eta}(w_{\kappa/(2l)})}
       +\left\|\Pi_MY^\phi\right\|^{l-1}_{\mathcal{C}^{\eta}(w_{\kappa/(2l)})}\right\}
       \left\|\Wick{(\Pi_NZ)^{k-l}}
       \right\|_{\mathcal{C}^{-\eta/2}(w_{\kappa/2})}
\end{align}
and 
\begin{align}
 &\left\|   (\Pi_M Y^\phi )^l \left\{\Wick{(\Pi_N Z)^{k-l}} -  \Wick{(\Pi_M Z )^{k-l}}   \right\}  \right\|_{\mathcal{C}^{-\eta}(w_\kappa)} \notag\\
 &\lesssim \left\| \Pi_M Y^\phi \right\|^{l}_{\mathcal{C}^{\eta}(w_{\kappa/(2l)})}   \left\|  \Wick{(\Pi_N Z )^{k-l}} -  \Wick{(\Pi_M Z )^{k-l}}     \right\|_{\mathcal{C}^{-\eta/2}(w_{\kappa/2})}  .  
\end{align}
Together with Lemma~\ref{lem:appr-1}, Corollary~\ref{cor:wick} and \eqref{eq:bound-ex-Y}, these estimates imply the desired result.

\end{proof}


\ \\
\noindent {\Large \textbf{Acknowledgement.}}
The author would like to thank Professor Hiroshi Kawabi for helpful comments and constant encouragement.
This work is partially supported by JSPS KAKENHI Grant Number 24KJ1329.

\ \\
\noindent {\Large \textbf{Declaration of Generative AI Use} }
The author used generative AI for the purpose of improving the English phrasing and checking the mathematical arguments in Sections~\ref{sec:resolvents} and~\ref{sec:CDFI} based on the initial version of the manuscript fully written by the author. 
The author revised the output as needed and assume responsibility for the final publication.

\bibliographystyle{myalpha_v3}

\bibliography{phi42_dirichlet}

@article {AR91,
    AUTHOR = {Albeverio, S. and R\"ockner, M.},
     TITLE = {Stochastic differential equations in infinite dimensions:
              solutions via {D}irichlet forms},
   JOURNAL = {Probab. Theory Related Fields},
  FJOURNAL = {Probability Theory and Related Fields},
    VOLUME = {89},
      YEAR = {1991},
    NUMBER = {3},
     PAGES = {347--386},
      ISSN = {0178-8051,1432-2064},
   MRCLASS = {60H10 (34F05 34G99)},
  MRNUMBER = {1113223},
MRREVIEWER = {T.\ J.\ Lyons},
       DOI = {10.1007/BF01198791},
       URL = {https://doi-org.kyoto-u.idm.oclc.org/10.1007/BF01198791},
}

@book {MR92,
    AUTHOR = {Ma, Zhi Ming and R\"ockner, Michael},
     TITLE = {Introduction to the theory of (nonsymmetric) {D}irichlet
              forms},
    SERIES = {Universitext},
 PUBLISHER = {Springer-Verlag, Berlin},
      YEAR = {1992},
     PAGES = {vi+209},
      ISBN = {3-540-55848-9},
   MRCLASS = {60J40 (31C25 60J45)},
  MRNUMBER = {1214375},
MRREVIEWER = {Tu\ Sheng\ Zhang},
       DOI = {10.1007/978-3-642-77739-4},
       URL = {https://doi-org.kyoto-u.idm.oclc.org/10.1007/978-3-642-77739-4},
}

@book {S74,
    AUTHOR = {Simon, Barry},
     TITLE = {The {$P(\phi )\sb{2}$} {E}uclidean (quantum) field theory},
    SERIES = {Princeton Series in Physics},
 PUBLISHER = {Princeton University Press, Princeton, NJ},
      YEAR = {1974},
     PAGES = {xx+392},
   MRCLASS = {81.60},
  MRNUMBER = {489552},
MRREVIEWER = {Raymond\ Streater},
}

@book {Jan97,
    AUTHOR = {Janson, Svante},
     TITLE = {Gaussian {H}ilbert spaces},
    SERIES = {Cambridge Tracts in Mathematics},
    VOLUME = {129},
 PUBLISHER = {Cambridge University Press, Cambridge},
      YEAR = {1997},
     PAGES = {x+340},
      ISBN = {0-521-56128-0},
   MRCLASS = {60G35 (60H15 62M20 81T08)},
  MRNUMBER = {1474726},
MRREVIEWER = {Amarjit\ Budhiraja},
       DOI = {10.1017/CBO9780511526169},
       URL = {https://doi-org.kyoto-u.idm.oclc.org/10.1017/CBO9780511526169},
}

@article{MW17-phi43,
  author = {Mourrat, Jean-Christophe and Weber, Hendrik},
  title = {The dynamic {$\Phi^4_3$} model comes down from infinity},
  journal = {Comm. Math. Phys.},
  fjournal = {Communications in Mathematical Physics},
  volume = {356},
  year = {2017},
  number = {3},
  pages = {673--753},
  issn = {0010-3616,1432-0916},
  mrclass = {60H15 (35K15 35K91 35R60)},
  mrnumber = {3719541},
  doi = {10.1007/s00220-017-2997-4},
  url = {https://doi.org/10.1007/s00220-017-2997-4},
}

@book{BCD11,
  author = {Bahouri, Hajer and Chemin, Jean-Yves and Danchin, Rapha\"{e}l},
  title = {Fourier analysis and nonlinear partial differential equations},
  series = {Grundlehren der mathematischen Wissenschaften [Fundamental Principles of Mathematical Sciences]},
  volume = {343},
  publisher = {Springer, Heidelberg},
  year = {2011},
  pages = {xvi+523},
  isbn = {978-3-642-16829-1},
  mrclass = {35-02 (35L72 35Q30 42-02 42B37 76B03 76D03 76N10)},
  mrnumber = {2768550},
  mrreviewer = {Peter¥ R.¥ Massopust},
  doi = {10.1007/978-3-642-16830-7},
  url = {https://doi.org/10.1007/978-3-642-16830-7},
}

@article {DPD03,
    AUTHOR = {Da Prato, Giuseppe and Debussche, Arnaud},
     TITLE = {Strong solutions to the stochastic quantization equations},
   JOURNAL = {Ann. Probab.},
  FJOURNAL = {The Annals of Probability},
    VOLUME = {31},
      YEAR = {2003},
    NUMBER = {4},
     PAGES = {1900--1916},
      ISSN = {0091-1798,2168-894X},
   MRCLASS = {81S20 (28C20 42B05)},
  MRNUMBER = {2016604},
       DOI = {10.1214/aop/1068646370},
       URL = {https://doi.org/10.1214/aop/1068646370},
}

@article {KT25,
    AUTHOR = {Kunick, Florian and Tsatsoulis, Pavlos},
     TITLE = {Gradient-type estimates for the dynamic {$\varphi
              _2^4$}-model},
   JOURNAL = {Stochastic Process. Appl.},
  FJOURNAL = {Stochastic Processes and their Applications},
    VOLUME = {181},
      YEAR = {2025},
     PAGES = {Paper No. 104548, 21},
      ISSN = {0304-4149,1879-209X},
   MRCLASS = {60H17 (47D07)},
  MRNUMBER = {4843971},
MRREVIEWER = {Le\ Chen},
       DOI = {10.1016/j.spa.2024.104548},
       URL = {https://doi.org/10.1016/j.spa.2024.104548},
}

@article {AKMR23,
    AUTHOR = {Albeverio, Sergio and Kawabi, Hiroshi and Mihalache,
              Stefan-Radu and R\"ockner, Michael},
     TITLE = {Strong uniqueness for {D}irichlet operators related to
              stochastic quantization under exponential/trigonometric
              interactions on the two-dimensional torus},
   JOURNAL = {Ann. Sc. Norm. Super. Pisa Cl. Sci. (5)},
  FJOURNAL = {Annali della Scuola Normale Superiore di Pisa. Classe di
              Scienze. Serie V},
    VOLUME = {24},
      YEAR = {2023},
    NUMBER = {1},
     PAGES = {33--69},
      ISSN = {0391-173X,2036-2145},
   MRCLASS = {81S20 (60H17 60J46)},
  MRNUMBER = {4587740},
MRREVIEWER = {Nizar\ Demni},
       DOI = {10.2422/2036-2145.202105\_106},
       URL = {https://doi.org/10.2422/2036-2145.202105_106},
}

@article {KR07,
    AUTHOR = {Kawabi, Hiroshi and R\"ockner, Michael},
     TITLE = {Essential self-adjointness of {D}irichlet operators on a path
              space with {G}ibbs measures via an {SPDE} approach},
   JOURNAL = {J. Funct. Anal.},
  FJOURNAL = {Journal of Functional Analysis},
    VOLUME = {242},
      YEAR = {2007},
    NUMBER = {2},
     PAGES = {486--518},
      ISSN = {0022-1236,1096-0783},
   MRCLASS = {81T08 (47B25 47N50 60H15 60J35 60J45 81S20)},
  MRNUMBER = {2274818},
MRREVIEWER = {Wilhelm\ Stannat},
       DOI = {10.1016/j.jfa.2006.06.017},
       URL = {https://doi.org/10.1016/j.jfa.2006.06.017},
}

@article {DD20,
    AUTHOR = {Da Prato, Giuseppe and Debussche, Arnaud},
     TITLE = {Gradient estimates and maximal dissipativity for the
              {K}olmogorov operator in {$\Phi^4_2$}},
   JOURNAL = {Electron. Commun. Probab.},
  FJOURNAL = {Electronic Communications in Probability},
    VOLUME = {25},
      YEAR = {2020},
     PAGES = {Paper No. 9, 16},
      ISSN = {1083-589X},
   MRCLASS = {60H15 (35R60 81S20)},
  MRNUMBER = {4066302},
       DOI = {10.1214/20-ecp294},
       URL = {https://doi.org/10.1214/20-ecp294},
}

@incollection {AMR15,
    AUTHOR = {Albeverio, Sergio and Ma, Zhi-Ming and R\"ockner, Michael},
     TITLE = {Quasi regular {D}irichlet forms and the stochastic
              quantization problem},
 BOOKTITLE = {Festschrift {M}asatoshi {F}ukushima},
    SERIES = {Interdiscip. Math. Sci.},
    VOLUME = {17},
     PAGES = {27--58},
 PUBLISHER = {World Sci. Publ., Hackensack, NJ},
      YEAR = {2015},
      ISBN = {978-981-4596-52-7},
   MRCLASS = {31C25 (47D07 60J60 81S20 81T70)},
  MRNUMBER = {3379334},
MRREVIEWER = {Gerald\ Trutnau},
       DOI = {10.1142/9789814596534\_0003},
       URL = {https://doi.org/10.1142/9789814596534_0003},
}

@article {RZZ17,
    AUTHOR = {R\"ockner, Michael and Zhu, Rongchan and Zhu, Xiangchan},
     TITLE = {Restricted {M}arkov uniqueness for the stochastic quantization
              of {$P(\Phi)_2$} and its applications},
   JOURNAL = {J. Funct. Anal.},
  FJOURNAL = {Journal of Functional Analysis},
    VOLUME = {272},
      YEAR = {2017},
    NUMBER = {10},
     PAGES = {4263--4303},
      ISSN = {0022-1236,1096-0783},
   MRCLASS = {60H15 (82C28)},
  MRNUMBER = {3626040},
       DOI = {10.1016/j.jfa.2017.01.023},
       URL = {https://doi.org/10.1016/j.jfa.2017.01.023},
}

@article {LR98,
    AUTHOR = {Liskevich, Vitali and R\"ockner, Michael},
     TITLE = {Strong uniqueness for certain infinite-dimensional {D}irichlet
              operators and applications to stochastic quantization},
   JOURNAL = {Ann. Scuola Norm. Sup. Pisa Cl. Sci. (4)},
  FJOURNAL = {Annali della Scuola Normale Superiore di Pisa. Classe di
              Scienze. Serie IV},
    VOLUME = {27},
      YEAR = {1998},
    NUMBER = {1},
     PAGES = {69--91},
      ISSN = {0391-173X,2036-2145},
   MRCLASS = {47B25 (31C25 47D07 81S20)},
  MRNUMBER = {1658889},
MRREVIEWER = {J.\ A.\ van Casteren},
       URL = {http://www.numdam.org/item?id=ASNSP_1998_4_27_1_69_0},
}

@article {DT00,
    AUTHOR = {Da Prato, Giuseppe and Tubaro, Luciano},
     TITLE = {Self-adjointness of some infinite-dimensional elliptic
              operators and application to stochastic quantization},
   JOURNAL = {Probab. Theory Related Fields},
  FJOURNAL = {Probability Theory and Related Fields},
    VOLUME = {118},
      YEAR = {2000},
    NUMBER = {1},
     PAGES = {131--145},
      ISSN = {0178-8051,1432-2064},
   MRCLASS = {47B25 (81S20)},
  MRNUMBER = {1785456},
MRREVIEWER = {Vladimir\ I.\ Bogachev},
       DOI = {10.1007/PL00008739},
       URL = {https://doi.org/10.1007/PL00008739},
}

@misc{DGT22,
      title={A singular integration by parts formula for the exponential Euclidean {QFT} on the plane}, 
      author = {De Vecchi, Francesco C. and Gubinelli, Massimiliano and Turra, Mattia},
      year={2022},
      eprint={2212.05584},
      archivePrefix={arXiv},
      primaryClass={math.PR},
      url={https://arxiv.org/abs/2212.05584}, 
}

@misc{DHYZ25,
      title={Ergodicity of infinite volume $\Phi^4_3$ at high temperature}, 
      author={Paweł Duch and Martin Hairer and Jaeyun Yi and Wenhao Zhao},
      year={2025},
      eprint={2508.07776},
      archivePrefix={arXiv},
      primaryClass={math.PR},
      url={https://arxiv.org/abs/2508.07776}, 
}

@article {MW17,
    AUTHOR = {Mourrat, Jean-Christophe and Weber, Hendrik},
     TITLE = {Global well-posedness of the dynamic {$\Phi^4$} model in the
              plane},
   JOURNAL = {Ann. Probab.},
  FJOURNAL = {The Annals of Probability},
    VOLUME = {45},
      YEAR = {2017},
    NUMBER = {4},
     PAGES = {2398--2476},
      ISSN = {0091-1798,2168-894X},
   MRCLASS = {81T27 (35B30 35K58 35R60 60H15 81T40)},
  MRNUMBER = {3693966},
MRREVIEWER = {Dimitri\ Petritis},
       DOI = {10.1214/16-AOP1116},
       URL = {https://doi.org/10.1214/16-AOP1116},
}

@book {Ebe99,
    AUTHOR = {Eberle, Andreas},
     TITLE = {Uniqueness and non-uniqueness of semigroups generated by
              singular diffusion operators},
    SERIES = {Lecture Notes in Mathematics},
    VOLUME = {1718},
 PUBLISHER = {Springer-Verlag, Berlin},
      YEAR = {1999},
     PAGES = {viii+262},
      ISBN = {3-540-66628-1},
   MRCLASS = {60J60 (31C25 47D07 60H10 60J35)},
  MRNUMBER = {1734956},
MRREVIEWER = {Niels\ Jacob},
       DOI = {10.1007/BFb0103045},
       URL = {https://doi.org/10.1007/BFb0103045},
}

@book {Cer01,
    AUTHOR = {Cerrai, Sandra},
     TITLE = {Second order {PDE}'s in finite and infinite dimension},
    SERIES = {Lecture Notes in Mathematics},
    VOLUME = {1762},
      NOTE = {A probabilistic approach},
 PUBLISHER = {Springer-Verlag, Berlin},
      YEAR = {2001},
     PAGES = {x+330},
      ISBN = {3-540-42136-X},
   MRCLASS = {35R60 (35K57 47H20 47N20 60H10 60H15 60J35)},
  MRNUMBER = {1840644},
MRREVIEWER = {Krystyna\ Twardowska},
       DOI = {10.1007/b80743},
       URL = {https://doi.org/10.1007/b80743},
}

@article {OTWZ26,
    AUTHOR = {Oh, Tadahiro and Tolomeo, Leonardo and Wang, Yuzhao and Zheng,
              Guangqu},
     TITLE = {Hyperbolic {$P(\Phi)_2$}-model on the plane},
   JOURNAL = {Comm. Math. Phys.},
  FJOURNAL = {Communications in Mathematical Physics},
    VOLUME = {407},
      YEAR = {2026},
    NUMBER = {2},
     PAGES = {Paper No. 34, 84},
      ISSN = {0010-3616,1432-0916},
   MRCLASS = {35L05 (35R60 60H15 81S20)},
  MRNUMBER = {5016812},
       DOI = {10.1007/s00220-025-05486-0},
       URL = {https://doi.org/10.1007/s00220-025-05486-0},
}

@book {Tri06,
    AUTHOR = {Triebel, Hans},
     TITLE = {Theory of function spaces. {III}},
    SERIES = {Monographs in Mathematics},
    VOLUME = {100},
 PUBLISHER = {Birkh\"auser Verlag, Basel},
      YEAR = {2006},
     PAGES = {xii+426},
      ISBN = {978-3-7643-7581-2; 3-7643-7581-7},
   MRCLASS = {46E35 (46-02)},
  MRNUMBER = {2250142},
MRREVIEWER = {Lubo\v s\ Pick},
}

@book {Tri83,
    AUTHOR = {Triebel, Hans},
     TITLE = {Theory of function spaces},
    SERIES = {Monographs in Mathematics},
    VOLUME = {78},
 PUBLISHER = {Birkh\"auser Verlag, Basel},
      YEAR = {1983},
     PAGES = {284},
      ISBN = {3-7643-1381-1},
   MRCLASS = {46Exx},
  MRNUMBER = {781540},
       DOI = {10.1007/978-3-0346-0416-1},
       URL = {https://doi.org/10.1007/978-3-0346-0416-1},
}

@article{TW18,
  author = {Tsatsoulis, Pavlos and Weber, Hendrik},
  title = {Spectral gap for the stochastic quantization equation on the 2-dimensional torus},
  journal = {Ann. Inst. Henri Poincar\'e{} Probab. Stat.},
  year = {2018},
  fjournal = {Annales de l'Institut Henri Poincar\'e{} Probabilit\'es et Statistiques},
  volume = {54},
  number = {3},
  pages = {1204--1249},
  issn = {0246-0203,1778-7017},
  mrclass = {60H15 (47A10 81T08)},
  mrnumber = {3825880},
  mrreviewer = {Jiang\ Lun\ Wu},
  doi = {10.1214/17-AIHP837},
  url = {https://doi.org/10.1214/17-AIHP837},
}

@book{Str11,
  author = {Stroock, Daniel W.},
  title = {Probability theory},
  edition = {Second},
  note = {An analytic view},
  publisher = {Cambridge University Press, Cambridge},
  year = {2011},
  pages = {xxii+527},
  isbn = {978-0-521-13250-3},
  mrclass = {60-02 (60B11 60F05 60G44 60G51 60J35)},
  mrnumber = {2760872},
  mrreviewer = {Ofer\ Zeitouni},
}

@article{PW81,
  author = {Parisi, G. and Wu, Yong Shi},
  title = {Perturbation theory without gauge fixing},
  journal = {Sci. Sinica},
  fjournal = {Scientia Sinica. Zhongguo Kexue},
  volume = {24},
  year = {1981},
  number = {4},
  pages = {483--496},
  issn = {0582-236x},
  mrclass = {81E99},
  mrnumber = {626795},
}

@book{FOT11,
  author = {Fukushima, Masatoshi and Oshima, Yoichi and Takeda, Masayoshi},
  title = {Dirichlet forms and symmetric {M}arkov processes},
  series = {De Gruyter Studies in Mathematics},
  volume = {19},
  edition = {extended},
  publisher = {Walter de Gruyter \& Co., Berlin},
  year = {2011},
  pages = {x+489},
  isbn = {978-3-11-021808-4},
  mrclass = {60J25 (28A12 31C45 60F10 60J40 60J45 60J55)},
  mrnumber = {2778606},
}

@book{DPZ14,
  author = {Da Prato, Giuseppe and Zabczyk, Jerzy},
  title = {Stochastic equations in infinite dimensions},
  series = {Encyclopedia of Mathematics and its Applications},
  volume = {152},
  edition = {Second},
  publisher = {Cambridge University Press, Cambridge},
  year = {2014},
  pages = {xviii+493},
  isbn = {978-1-107-05584-1},
  mrclass = {60H15 (34F05 34Gxx)},
  mrnumber = {3236753},
  mrreviewer = {David\ Nualart},
  doi = {10.1017/CBO9781107295513},
  url = {https://doi.org/10.1017/CBO9781107295513},
}

@article{AKR12,
  author = {Albeverio, Sergio and Kawabi, Hiroshi and R\"ockner, Michael},
  title = {Strong uniqueness for both {D}irichlet operators and stochastic dynamics to {G}ibbs measures on a path space with exponential interactions},
  journal = {J. Funct. Anal.},
  fjournal = {Journal of Functional Analysis},
  volume = {262},
  year = {2012},
  number = {2},
  pages = {602--638},
  issn = {0022-1236,1096-0783},
  mrclass = {81T08 (47D06 60H15 81S20)},
  mrnumber = {2854715},
  mrreviewer = {J.\ A.\ van Casteren},
  doi = {10.1016/j.jfa.2011.09.023},
  url = {https://doi.org/10.1016/j.jfa.2011.09.023},
}

@article{Hai14,
  author = {Hairer, M.},
  title = {A theory of regularity structures},
  journal = {Invent. Math.},
  fjournal = {Inventiones Mathematicae},
  volume = {198},
  year = {2014},
  number = {2},
  pages = {269--504},
  issn = {0020-9910,1432-1297},
  mrclass = {60H15 (35R60 60H40 81S20 82C28)},
  mrnumber = {3274562},
  mrreviewer = {Dora\ Sele\v{s}i},
  doi = {10.1007/s00222-014-0505-4},
  url = {https://doi.org/10.1007/s00222-014-0505-4},
}

@article{GIP15,
  author = {Gubinelli, Massimiliano and Imkeller, Peter and Perkowski, Nicolas},
  title = {Paracontrolled distributions and singular {PDE}s},
  journal = {Forum Math. Pi},
  fjournal = {Forum of Mathematics. Pi},
  volume = {3},
  year = {2015},
  pages = {e6, 75},
  issn = {2050-5086},
  mrclass = {60H15 (35S50)},
  mrnumber = {3406823},
  doi = {10.1017/fmp.2015.2},
  url = {https://doi.org/10.1017/fmp.2015.2},
}

@article{ALZ06,
  author = {Albeverio, Sergio and Liang, Song and Zegarlinski, Boguslav},
  title = {Remark on the integration by parts formula for the {$\phi^4_3$}-quantum field model},
  journal = {Infin. Dimens. Anal. Quantum Probab. Relat. Top.},
  fjournal = {Infinite Dimensional Analysis, Quantum Probability and Related Topics},
  volume = {9},
  year = {2006},
  number = {1},
  pages = {149--154},
  issn = {0219-0257,1793-6306},
  mrclass = {81T08 (46N50)},
  mrnumber = {2214506},
  mrreviewer = {Shigeaki\ Nagamachi},
  doi = {10.1142/S0219025706002275},
  url = {https://doi.org/10.1142/S0219025706002275},
}

@article{ZZ18,
  author = {Zhu, Rongchan and Zhu, Xiangchan},
  title = {Dirichlet form associated with the {$\Phi^4_3$} model},
  journal = {Electron. J. Probab.},
  fjournal = {Electronic Journal of Probability},
  volume = {23},
  year = {2018},
  pages = {Paper No. 78, 31},
  issn = {1083-6489},
  mrclass = {60H15 (60J35 82C28)},
  mrnumber = {3858906},
  doi = {10.1214/18-EJP207},
  url = {https://doi.org/10.1214/18-EJP207},
}

@article{AK20,
  author = {Albeverio, Sergio and Kusuoka, Seiichiro},
  title = {The invariant measure and the flow associated to the {$\Phi^4_3$}-quantum field model},
  journal = {Ann. Sc. Norm. Super. Pisa Cl. Sci. (5)},
  fjournal = {Annali della Scuola Normale Superiore di Pisa. Classe di Scienze. Serie V},
  volume = {20},
  year = {2020},
  number = {4},
  pages = {1359--1427},
  issn = {0391-173X,2036-2145},
  mrclass = {81S20 (35Q40 37H10 60H15 60L50 81T08)},
  mrnumber = {4201185},
}

@article{GH21,
  author = {Gubinelli, Massimiliano and Hofmanov\'{a}, Martina},
  title = {A {PDE} construction of the {E}uclidean {$\phi_3^4$} quantum field theory},
  journal = {Comm. Math. Phys.},
  fjournal = {Communications in Mathematical Physics},
  volume = {384},
  year = {2021},
  number = {1},
  pages = {1--75},
  issn = {0010-3616,1432-0916},
  mrclass = {81S20 (35Q40 37D35 37E20 60G60 81T08 81T27 82B31)},
  mrnumber = {4252872},
  mrreviewer = {Peter\ Bernard\ Weichman},
  doi = {10.1007/s00220-021-04022-0},
  url = {https://doi.org/10.1007/s00220-021-04022-0},
}

@article{AK22,
  author = {Albeverio, Sergio and Kusuoka, Seiichiro},
  title = {Construction of a non-{G}aussian and rotation-invariant {$\Phi^4$}-measure and associated flow on {$\Bbb R^3$} through stochastic quantization},
  journal = {Mem. Amer. Math. Soc.},
  fjournal = {Memoirs of the American Mathematical Society},
  volume = {308},
  year = {2025},
  number = {1558},
  pages = {viii+114},
  issn = {0065-9266,1947-6221},
  isbn = {978-1-4704-7318-1; 978-1-4704-8145-2},
  mrclass = {81S20 (35K58 35Q40 35R60 60H15 81T08)},
  mrnumber = {4896167},
  doi = {10.1090/memo/1558},
  url = {https://doi.org/10.1090/memo/1558},
}

@article{GH19,
  author = {Gubinelli, Massimiliano and Hofmanov\'{a}, Martina},
  title = {Global solutions to elliptic and parabolic {$\Phi^4$} models in {E}uclidean space},
  journal = {Comm. Math. Phys.},
  fjournal = {Communications in Mathematical Physics},
  volume = {368},
  year = {2019},
  number = {3},
  pages = {1201--1266},
  issn = {0010-3616,1432-0916},
  mrclass = {35Q40 (35J91 35K91 35R60 60H15 60H99 81T17 82C31)},
  mrnumber = {3951704},
  mrreviewer = {Peter\ Bernard\ Weichman},
  doi = {10.1007/s00220-019-03398-4},
  url = {https://doi.org/10.1007/s00220-019-03398-4},
}

@article{CC18,
  author = {Catellier, R\'emi and Chouk, Khalil},
  title = {Paracontrolled distributions and the 3-dimensional stochastic quantization equation},
  journal = {Ann. Probab.},
  fjournal = {The Annals of Probability},
  volume = {46},
  year = {2018},
  number = {5},
  pages = {2621--2679},
  issn = {0091-1798,2168-894X},
  mrclass = {60H15 (35K15 35R60 42B25 60H40)},
  mrnumber = {3846835},
  mrreviewer = {Dora\ Sele\v si},
  doi = {10.1214/17-AOP1235},
  url = {https://doi.org/10.1214/17-AOP1235},
}

@article{HM18,
  author = {Hairer, M. and Matetski, K.},
  title = {Discretisations of rough stochastic {PDE}s},
  journal = {Ann. Probab.},
  year = {2018},
  fjournal = {The Annals of Probability},
  volume = {46},
  number = {3},
  pages = {1651--1709},
  issn = {0091-1798,2168-894X},
  mrclass = {60H15 (35R60 65M75)},
  mrnumber = {3785597},
  mrreviewer = {Le\ Chen},
  doi = {10.1214/17-AOP1212},
  url = {https://doi.org/10.1214/17-AOP1212},
}

\end{document}